\documentclass[a4paper,10pt]{amsart}
\usepackage{amsmath,amsfonts,amssymb, amsthm, xypic}
\usepackage{epsfig,psfrag}
\usepackage[all]{xy}

\def\Mc{{\mathcal M}}
\def\Oc{{\mathcal O}}
\def\Ic{{\mathcal I}}
\def\Nc{{\mathcal N}}

\def\Z{\mathbb{Z}}
\def\Q{\mathbb{Q}}
\def\E{\mathcal{X}}
\def\Ext{\mathrm{Ext}}
\def\H{\mathbf{H}}

\def\Fb{\mathbf{F}}
\def\R{\mathbf{R}}
\def\U{\mathbf{U}}
\def\SH{\mathbf{S}\ddot{\mathbf{H}}}
\def\qed{\hfill $\sqcap \hskip-6.5pt \sqcup$}
\def\N{\mathbb{N}}
\def\tto{\twoheadrightarrow}

\def\a{\alpha}

\def\qlb{\overline{\mathbb{Q}}_p}
\def\cc{\mathbf{c}}
\def\C{\mathbb{C}}

\def\x{\mathbf{x}}
\def\y{\mathbf{y}}
\def\eE{\mathbf{E}}
\def\z{\mathbf{z}}
\def\p{\mathbf{p}}

\def\aeE{\hspace{.01in}_{\mathcal{A}}\eE}

\newtheorem{theo}{\bf{Theorem}}[section]
\newtheorem{lem}[theo]{Lemma}
\newtheorem{slem}[theo]{Sublemma}
\newtheorem{cor}[theo]{Corollary}
\newtheorem{prop}[theo]{Proposition}

\newtheorem*{theor}{\bf{Theorem}}

\numberwithin{equation}{section}

\title[Hall, Cherednik, Eisenstein, Macdonald]{The elliptic Hall algebra, Cherednik Hecke algebras and Macdonald polynomials.}
\author{O. Schiffmann, E. Vasserot}

\begin{document}
\maketitle

\tableofcontents

\section{Introduction}

\vspace{.2in}

\paragraph{} In this paper we continue the study of 
the Hall algebra $\H_{\E}$ of an
elliptic curve $\E$ defined over a finite field $\mathbb{F}_l$, 
started in \cite{BS}, \cite{Scano}.
Here we exhibit some strong link between the Hall algebra $\H_{\E}$, 
or more precisely its composition subalgebra $\U_{\E}$, and
Cherednik's double affine Hecke algebras $\ddot{\H}_n$ of type $GL(n)$, 
for all $n$.
This allows us to obtain a geometric construction of 
Macdonald polynomials $P_{\lambda}(q,t)$ in terms of certain
functions (Eisenstein series) on the moduli space of 
semistable vector bundles on the elliptic curve $\E$.

\vspace{.1in}

Let us describe our results in more details. The spherical affine
Hecke algebra $\mathbf{S}\dot{\H}_G$ of a reductive algebraic group
$G$ is the convolution algebra of $G(\mathcal{O})$-invariant
functions on the affine Grassmanian $\widehat{Gr}=
G(\mathcal{K})/G(\mathcal{O})$, where
$\mathcal{K}=\mathbb{F}_l((z))$ and $\mathcal{O}=\mathbb{F}_l[[z]]$,
see \cite{IM}. The Satake isomorphism identifies
$\mathbf{S}\dot{\H}_G$ with the representation ring $Rep(G^L)$ of
the dual group of $G$. Now let us assume that $G=G^L=GL(n)$, so that
the set of $\mathbb{F}_l$-points of $\widehat{Gr}$ is equal to
$$\{L \subset \mathbb{F}_l^n((z)); L\;\text{is\;a\;free\;}
\mathbb{F}_l[[z]]-\text{module\;of\;rank\;}n\}$$ and $Rep(G)\simeq
\C[x_1^{\pm 1}, \ldots, x_n^{\pm 1}]^{\mathfrak{S}_n}$. Lusztig
\cite{LusztigGreen} embeds the nilpotent cones $\mathcal{N}_k
\subset \mathfrak{gl}(k)$, $k\ge 1$, into the positive Schubert
variety
$$\widehat{Gr}^+=\{L \subset \mathbb{F}_l^n[[z]]; L\;\text{is\;a\;free\;} \mathbb{F}_l[[z]]-\text{module\;of\;rank\;}n\}$$
of $\widehat{Gr}$.
This yields a surjective algebra homomorphism
\begin{equation}\label{E:intro1}
\Theta^+_n~: \H_{cl}=\bigoplus_{k \geq 0} \C_{GL(k)}[\mathcal{N}_k]
\tto \mathbf{S}\dot{\H}^+_n\simeq \C[v^{\pm 1}][x_1, \ldots,
x_n]^{\mathfrak{S}_n}.
\end{equation}
See \cite[Chap. II]{Mac} or \cite{LusztigGreen}.
Here $\H_{cl}$ is the classical Hall algebra,
and $v=l^{-1/2}$.
Since the dependence on $v$ is polynomial,
we may treat it as a formal parameter.
Letting $n$ tend to infinity in (\ref{E:intro1}) yields an isomorphism in the stable limit
\begin{equation}\label{E:intro2}
\Theta^+_{\infty}~: \H_{cl} \stackrel{\sim}{\to} \mathbf{S}\dot{\H}_{\infty}^+=\varprojlim \mathbf{S}\dot{\H}_{n}^+\simeq\C[v^{\pm 1}][x_1,x_2, \ldots]^{\mathfrak{S}_{\infty}}.
\end{equation}
The first main result of this paper provides an affine version of (\ref{E:intro1}) and (\ref{E:intro2}).

In \cite{BS} it was found that the Hall algebra $\H_{\E}$ of the category of coherent sheaves on an elliptic curve $\E$ defined over $\mathbb{F}_l$ contains a natural ``composition'' subalgebra $\U_{\E}$ which is a two-parameter
deformation of the ring of diagonal invariants
$$\R^+_n=
\C[x_1, \ldots, x_n, y_1^{\pm 1}, \ldots, y_n^{\pm 1}]^{\mathfrak{S}_n}$$
where $\mathfrak{S}_n$ acts simultaneously on the $x$-variables
and the $y$-variables.
The two deformation parameters are
$v=l^{-1/2}$ and $t=\sigma v$,
where $\sigma$ is a Frobenius eigenvalue of $H^1(\E,\overline{\Q}_p)$
(viewed as a complex number).
The dependence on $v$, $t$ is polynomial
and we may treat them as formal variables.

Let $\ddot{\H}_n$ denote Cherednik's double affine Hecke algebra of
type $GL(n)$, and let $\SH_n=S \cdot \ddot{\H}_n \cdot S$ stand for
its spherical subalgebra. Here $S$ is the complete idempotent
associated to the \textit{finite} Hecke algebra $\H_n \subset
\ddot{\H}_n$. The algebra $\SH_n$ is a deformation of the ring
$$\R_n=\C[x_1^{\pm 1}, \ldots, x_n^{\pm 1}, y_1^{\pm 1},
\ldots, y_n^{\pm 1}]^{\mathfrak{S}_n}$$ depending on two parameters
$v$ and $q$. Let $\SH^+_n$ be the positive part of $\SH_n$, see
Section~2.1. In Theorem~\ref{T:1} we prove the following.

\vspace{.05in}

\begin{theor} 
If $q=vt$ then for any $n$ there exists a surjective algebra homomorphism
$\Psi^+_n: \U_{\E} \tto \SH^+_n.$ This map extends to a surjective
algebra homomorphism $\Psi_n: \mathbf{D}\U_{\E} \tto \SH_n.$ 
\end{theor}

\vspace{.05in}

Here $\mathbf{D}\U_{\E}$ is the Drinfeld double of $\U_{\E}$.
It is equipped with an action of
$SL(2,\Z)$ coming from the group of derived autoequivalences of
$D^b(Coh(\E))$. Cherednik has defined an action of $SL(2,\Z)$ on
$\SH_n$, see \cite{Cherednik}, \cite{Ion}. The map $\Psi_n$ is
defined so as to intertwine these two actions. The maps $\Psi^+_n$
behave well with respect to the stable limit, see
Theorem~\ref{T:infty}.

\begin{theor} The maps $\Psi^+_n$ induce an algebra isomomorphism
$$\Psi^+_{\infty}: \U_{\E} \stackrel{\sim}{\to} \SH^+_\infty=\varprojlim \SH^+_n.$$
\end{theor}

Note that our approach doesn't use any affinization of the affine Grassmanian.

\vspace{.1in}

One of the essential features of the construction of the spherical affine Hecke
algebras as convolution algebra of functions
(on the affine Grassmanian or on the nilpotent cones)
is that it lifts to a tensor category of perverse sheaves
(see e.g. \cite{Ginzburg}, \cite{MV}).
Such a geometric lift also makes sense here, and fits in Laumon's theory
of automorphic sheaves.
We refer to \cite{Scano} and Section~4.3. for more details.

\vspace{.1in}

In the second part of this paper, we give an application of the
above geometric construction of $\SH_n$ to Macdonald polynomials.

\vspace{.1in}

The Hall algebra $\H^{vec}_\E$ of the category of vector bundles on
$\E$ (or on any smooth projective curve) can be viewed as the
algebra of (unramified) automorphic forms for $GL(n)$, for all $n
\geq 1$, over the function field of $\E$. The product is given by
the functor of parabolic induction, see \cite{Kap}. To obtain the
whole Hall algebra $\H_\E$ one needs to take into account the
torsion sheaves as well. The Hall algebra $\H^{tor}_\E$ of the
category of torsion sheaves on $\E$ acts on $\H^{vec}_\E$ by the
adjoint action and $\H_\E$ is isomorphic to the semi-direct product
$\H^{vec}_\E \rtimes \H^{tor}_\E$. The action of torsion sheaves is
interpreted in the language of automorphic forms as \textit{Hecke
operators}. For instance, the skyscraper sheaf $\mathcal{O}_x$ at a
point $x \in \E$ corresponds to the elementary modification at $x$.

\vspace{.1in}

Under the map $\Psi^+_{\infty}$, the element
$\mathbf{1}_{(0,1)} \in \U^{tor}_{\E}$ responsible for the
\textit{Hecke operator} of rank one is sent to
\textit{Macdonald's element} $\Delta_1=S \sum_i Y_i S \in \SH^+_{\infty}$,
see Section~2. The importance of this element stems from the fact that,
in the polynomial representation of $\SH^+_{\infty}$, the operator
$\Delta_1$ has distinct eigenvalues and the corresponding eigenvectors are
the Macdonald polynomials $P_{\lambda}(q,v^2)$.
Thus the map $\Psi^+_{\infty}$ allows us to relate Hecke eigenvectors
on the Hall or automorphic side to Macdonald polynomials on the Hecke
algebra side. In particular, we are naturally led to find Hecke eigenvectors in $\U_{\E}$ whose eigenvalues match those of the $P_{\lambda}(q,v^2)$.

\vspace{.1in}

Eisenstein series yield a way to produce new Hecke eigenvectors from
old ones via parabolic induction. In the present situation, it so
happens that the simplest Eisenstein series, i.e., those induced
from trivial characters of parabolic subgroups, already have the
good eigenvalues under the Hecke operator. Unfortunately, we are
unable to construct the polynomial representation of
$\SH^+_{\infty}$ in a geometric way, see Remark 5.1, and thus we
cannot obtain directly a geometric construction of
$P_{\lambda}(q,v^2)$. To remedy this, we manage to lift the
Macdonald polynomials from the polynomial representation and view
them inside the Hecke (or Hall) algebra itself. More precisely, it
is shown in \cite{BS} that the subalgebra of $\U_{\E}$ consisting of
functions supported on the set $Coh(\E)^{(0)}$ of semistable sheaves
of zero slope is canonically isomorphic to the algebra
$$\Lambda^+_{v,t}= \C[v^{\pm 1}, t^{\pm 1}][x_1,
x_2,\ldots]^{\mathfrak{S}_{\infty}}.$$ In a few words, under
Fourier-Mukai transform the set $Coh(\E)^{(0)}$ is identified with
the set of torsion sheaves on $\E$, and any function on the set of
torsion sheaves with a fixed ponctual support in $\E$ can be viewed
as an element of the classical Hall algebra. See \cite{P}, Theorem
14.7 for details. If $f$ is any function in $\U_{\E}$ we let $f_{(0)}$ be its
restriction to $Coh(\E)^{(0)}$, viewed as an element of
$\Lambda^+_{v,t}$. For any $l \in \N^+$ put
\begin{equation}
\eE_l(z)=\sum_{d \in \Z} \mathbf{1}_{(l,d)}v^{d(l-1)}z^d \in
\widehat{\U}_{\E}[[z,z^{-1}]]
\end{equation}
where $\mathbf{1}_{(l,d)}$ is the characteristic function of the set
of all coherent sheaves on $\E$ of rank $l$ and degree $d$, and
$\widehat{\U}_{\E}$ is a certain completion of $\U_{\E}$. For $(l_1,
\ldots, l_n) \in \N^n$ we form the Eisenstein series
\begin{equation}
\eE_{l_1, \ldots, l_n}(z_1, \ldots, z_n) =\eE_{l_1}(z_1) \cdot
\eE_{l_2}(z_2) \cdots \eE_{l_n}(z_n)\in \widehat{\U}_{\E}[[z_1^{\pm
1}, \ldots, z_n^{\pm 1}]].
\end{equation}
By a theorem of Harder, this is a rational function in $z_1, \ldots,
z_n$. Our second main result (Theorem~\ref{T:MacEis}) reads as
follows.

\vspace{.05in}

\begin{theor} Let $(l_1, \ldots, l_n) \in \N^n$. We have
$\eE_{l_1, \ldots, l_n}(z,q^{-1}z, \ldots, q^{1-n}z)=0$ unless
$(l_1, \ldots, l_n)$ is dominant, i.e., unless $(l_1, \ldots, l_n)$
is a partition. In that case we have
$$\eE_{l_1, \ldots, l_n}(z,q^{-1}z, \ldots, 
q^{1-n}z)_{(0)}=\omega P_{\lambda}(q,v^2)$$
where $\lambda=(l_1, \ldots, l_n)'$ is the conjugate partition and where $\omega$ stands
for the standard involution on symmetric functions.
\end{theor}

We also give a similar construction of skew Macdonald polynomials
$P_{\lambda/\mu}(q,v^2)$. Note that the above Eisenstein series can
be lifted to some constructible sheaves via the theory of
\textit{Eisenstein sheaves}, see \cite{Laumon} and \cite{Scano}.
Hence the Macdonald polynomials $P_{\lambda}(q,v^2)$ may be realized
as Frobenius traces of certain canonical constructible sheaves on
the moduli stack of semistable sheaves of zero slope on $\E$. We
hope to come back to this point in the future.

\vspace{.15in}

There is a well-known and important geometric approach to Macdonald
polynomials, which is based on the equivariant K-theory of the
Hilbert schemes $(\C^2)^{[n]}$ of points on $\C^2$. There the
polynomials $P_{\lambda}(q,v^2)$ are realized as the classes of
certain canonical coherent sheaves on $(\C^2)^{[n]}$, see
\cite{Haiman}. It would be interesting to relate this ``coherent
sheaf'' picture with our ``constructible functions'' (or ``perverse
sheaf'') picture in a precise fashion, and to understand this
relation in the framework of Langlands duality, see e.g. \cite{Bez}.
Note that the Hall algebra $\H_{\E}$ is already involved in a
Langlands type duality (the geometric Langlands duality for the
elliptic curve $\E$), with the category of coherent sheaves on the
moduli stack of local systems on $\E$.

\vspace{.15in}

The structure of the paper is as follows : Sections~1 and 2 contain
some recollections on elliptic Hall algebras $\H_{\E}$ and
$\U_{\E}$, taken from \cite{BS}, and Cherednik double affine Hecke
algebra $\ddot{\H}_n$ and $\SH_n$ respectively; in Section~3 we
construct the algebra morphism $\Psi_n: \mathbf{D}\U_{\E} \tto
\SH_n$; in Section~4 we study and define the stable limit
$\SH^+_{\infty}$ of the spherical Cherednik algebra and establish
the isomorphism $\Psi^+_{\infty}: \U_{\E} \stackrel{\sim}{\to}
\SH^+_{\infty}$. This is the first main result of this paper. A
comparison table with the picture of the classical Hall algebra and
the ``finite'' spherical affine Hecke algebra is found in
Section~4.3. Section~5 deals with Macdonald polynomials~: we recall
the definition and provide a characterization of the family of all
(possibly skew) Macdonald polynomials which we use later. In
Section~6 we introduce the Eisenstein series which are relevant to
us, and study some of their specializations. Finally, our second
main theorem, which gives a geometric construction of (possibly
skew) Macdonald polynomials from Eisenstein series, is given in
Section~7. Several proofs in the paper require some lengthy
computations; these are written up in details in Appendices A
through D.

\vspace{.1in}

A final word of warning concerning the notation.
There is an unfortunate clash between the conventional notations
used in the quantum group/Hall algebra literature and those used in the
Macdonald polynomials literature~: the letter $q$ generally
denotes the size of the finite field in the first case whereas
it denotes the modular parameter in the second case.
We have opted to comply with the Macdonald polynomials conventions.

\newpage

\vspace{.2in}

\section{Elliptic Hall algebra}

\vspace{.2in}

\paragraph{\textbf{1.1.}} We will use the standard $v$-integers and $v$-factorials
$$[i]=[i]_v=\frac{v^i-v^{-i}}{v-v^{-1}}, \qquad [i]!=[2] \cdots [i], $$
as well as some positive and negative variants
$$ [i]^+=\frac{v^{2i}-1}{v^{2}-1}, \quad [i]^+!=[2]^+ \cdots [i]^+, \qquad
[i]^-=\frac{v^{-2i}-1}{v^{-2}-1}, \quad [i]^-!=[2]^- \cdots [i]^-.$$
Let us denote by $\Lambda^+_v$ Macdonald's ring of symmetric functions (\cite{Mac})
$$\Lambda^+_v=\C[v^{\pm 1}] [x_1, x_2, \ldots]^{\mathfrak{S}_{\infty}}$$
defined over $\C[v^{\pm 1}]$. 
We will denote by $e_{\lambda}$, $p_{\lambda}$, $m_\lambda$
the elementary, the power-sum, and the monomial symmetric functions
respectively.
This ring is equipped with a natural bialgebra structure
$\Delta: \Lambda^+_v \to \Lambda^+_v \otimes \Lambda^+_v$ defined by
$\Delta(p_r)=p_r \otimes 1 + 1 \otimes p_r$ for $r \geq 1$.

\vspace{.2in}

\paragraph{\textbf{1.2}} Let $\E$ be a smooth elliptic curve over some finite field $\mathbb{F}_l$,
and let $Coh(\E)$ stand for the category of coherent sheaves on
$\E$. If $F$ is a sheaf sheaf in $Coh(\E)$ we call the pair
$\overline{F}=(rk(F),deg(F))$ the \textit{class} of $F$. The set of
possible classes of sheaves in $Coh(\E)$ is equal to
$\Z^{2,+}=\{(r,d) \in \Z^2; r \geq 1\; \text{or}\; r=0, d \geq
0\}$. We briefly recall the definition of the Hall algebra of
$Coh(\E)$. See \cite{BS} for more details.

\vspace{.1in}

Let $\Ic(\E)$ stand for the set of isomorphism classes of
objects in $Coh(\E)$.  Following Ringel \cite{Ringel}, the
$\C$-vector space of finitely supported functions
$$\mathbf{H}_\E=\{f: \Ic(\E) \to \C; |supp(f)| < \infty\}$$
may be equipped with the convolution product
$$(f \cdot g)(M)=\sum_{N \subseteq M} \nu^{-\langle M/N,N \rangle} f(M/N)g(N),$$
where $\nu=l^{-\frac{1}{2}}$ and $\langle P,Q
\rangle=dim\;\text{Hom}(P,Q)-dim\; \text{Ext}(P,Q)$ is the Euler form. 
Here we write $\text{Ext}(P,Q)$ for $\text{Ext}^1(P,Q)$. 
The sum on
the right hand side is finite for any $M$ since $f$ and $g$ have
finite support and, for any $N,M \in Coh(\E)$, the group $\text{Hom}(N,M)$
is finite. The above formula indeed defines an element in
$\mathbf{H}_{\E}$ as for any $P,Q \in Coh(\E)$, the group $\text{Ext}(P,Q)$
is also finite. By the Riemann-Roch theorem, we have
\begin{equation}\label{E:Euler}
\langle P, Q \rangle=rk(P)deg(Q)-deg(P) rk(Q).
\end{equation}

\vspace{.05in}

By \cite{Green} the algebra $\mathbf{H}_{\E}$ also has the structure
of a bialgebra, with coproduct
$$(\Delta(f))(P,Q)=\frac{1}{|\text{Ext}(P,Q)|}\sum_{\xi \in \Ext(P,Q)} f(M_{\xi})$$
where $M_{\xi}$ is the extension of $P$ by $Q$ corresponding to
$\xi$. The product and the coproduct are related by the pairing
$$\mathbf{H}_{\E} \otimes \mathbf{H}_{\E} \to \C, \; (f,g)=\sum_{M} \frac{f(M)g(M)}{|\text{Aut}(M)|}$$
which is a Hopf pairing, i.e., which satisfies the identity
$(fg,h)=(f \otimes g, \Delta(h))$ for any $f,g,h$.

\vspace{.05in}

\paragraph{\textbf{Remarks.}} i) In our situation, as opposed to \cite{Green},
it is not necessary to twist the product in $\mathbf{H}_{\E} \otimes
\mathbf{H}_{\E}$ in order to obtain a bialgebra, because the Euler form $\langle \;,\;\rangle$ is antisymmetric.\\
ii) The coproduct $\Delta$ only takes values in a certain formal
completion of $\mathbf{H}_{\E} \otimes \mathbf{H}_{\E}$, see
\cite[Section~2.2]{BS} for details.

\vspace{.1in}

The characteristic functions $\{\mathbf{1}_{M}; M \in
\Ic(\E)\}$ form a basis $\mathbf{H}_{\E}$. Assigning the degree
$(rk(M), deg(M))$ to the element $\mathbf{1}_{M}$ yields a
$\Z^2$-grading on $\mathbf{H}_{\E}$ which is compatible with the
(co)multiplication.

\vspace{.2in}

\paragraph{\textbf{1.3.}} Let $\mu(M)=deg(M)/rk(M) \in \mathbb{Q} \cup \{\infty\}$
be the slope of a sheaf $M \in Coh(\E)$, and for $\mu\in
\mathbb{Q} \cup \{\infty\}$ let $\mathcal{C}_\mu$ stand for the
category of semistable sheaves of slope $\mu$. For instance,
$\mathcal{C}_{\infty}$ is the category of torsion sheaves on $\E$.
The following fundamental result on the structure of $Coh(\E)$ is
due to Atiyah.

\begin{theo}[Atiyah, \cite{At}]\label{T:At} The following hold~:\\
i) for any $\mu$, $\mu'$ there is an equivalence of abelian categories
$\epsilon_{\mu,\mu'} : 
\mathcal{C}_{\mu'} \stackrel{\sim}{\to} \mathcal{C}_{\mu},$
\\
ii) any coherent sheaf $\mathcal{F}$ decomposes uniquely as a direct sum
$\mathcal{F}=\mathcal{F}_1 \oplus \cdots \oplus \mathcal{F}_s$
of semistable sheaves $\mathcal{F}_i \in \mathcal{C}_{\mu_i}$ with 
$\mu_1 < \cdots < \mu_s$.
\end{theo}

\vspace{.05in}

By a standard property of semistable sheaves we have
$\text{Hom}(\mathcal{C}_{\mu},\mathcal{C}_{\mu'})=\{0\}$ for 
$\mu > \mu'$. By Serre duality, this implies that
$\text{Ext}(\mathcal{C}_{\mu'},\mathcal{C}_{\mu})=\{0\}$ whenever
$\mu >\mu'$. Hence any extension $0 \to \mathcal{F} \to
\mathcal{G} \to \mathcal{H} \to 0$ with $\mathcal{F} \in
\mathcal{C}_{\mu}, \mathcal{H} \in \mathcal{C}_{\mu'}$ is split.
From the above two facts, it follows that in $\mathbf{H}_{\E}$ we
have
\begin{equation}\label{E:131}
\mathbf{1}_{\mathcal{H}} \cdot \mathbf{1}_{\mathcal{F}}=
\nu^{-\langle \mathcal{H},\mathcal{F}\rangle} \mathbf{1}_{\mathcal{F} \oplus \mathcal{H}}
\end{equation}
{if} $\mathcal{F} \in \mathcal{C}_{\mu},\mathcal{H}\in \mathcal{C}_{\mu'}$ 
{and} $\mu>\mu'$.

\vspace{.1in}

For $\mu\in \Q \cup \{\infty\}$ let $\mathbf{H}_{\E}^{(\mu)}$
stand for the subspace consisting of functions supported on the set
of semistables sheaves of slope $\mu$. Since $\mathcal{C}_{\mu}$
is stable under extensions, $\mathbf{H}_{\E}^{(\mu)}$ is a
subalgebra of $\mathbf{H}_{\E}$. By Theorem~\ref{T:At}~i), all these
subalgebras are isomorphic. Let $\vec{\bigotimes}_{\mu}
\mathbf{H}_{\E}^{(\mu)}$ denote the ordered tensor product of
spaces $\mathbf{H}^{(\mu)}_{\E}$ with $\mu\in \mathbb{Q} \cup
\{\infty\}$, i.e., the vector space spanned by elements of the form
$a_{\mu_1} \otimes \cdots \otimes a_{\mu_r}$ with $a_{\mu_i} \in
\mathbf{H}_{\E}^{(\mu_i)}$ and $\mu_1 < \cdots < \mu_r$. From
(\ref{E:131}) and Theorem~\ref{T:At}~ii) we deduce the following,
see \cite[Lemma~2.4.]{BS}).

\begin{cor} The multiplication map induces an isomorphism of vector spaces
$\vec{\bigotimes}_\mu\mathbf{H}_{\E}^{(\mu)}\stackrel\sim\to\mathbf{H}_{\E}$.
\end{cor}

\vspace{.2in}

\paragraph{\textbf{1.4.}} We will mainly be interested in a certain subalgebra
$\mathbf{U}^+_{\E} \subset \mathbf{H}_{\E}$ which we now define. For
any class $\a \in \Z^{2,+}$ we set
$$\mathbf{1}^{ss}_{\a}=\sum_{\overline{{F}}=\a; {F} \in \mathcal{C}_{\mu(\a)}}
\mathbf{1}_{F} \in \mathbf{H}_{\E}.$$ This sum is finite.
Indeed, by
Theorem~\ref{T:At}~i) it is enough to check this for
$\mu(\a)=\infty$. Then, this follows from the fact that there are
only finitely many closed points on $\E$ which are rational over
a fixed finite extension of $\mathbb{F}_l$. Let
$\mathbf{U}^+_{\E}$ be the subalgebra generated by
$\mathbf{1}^{ss}_{\a}$ for $\a \in \Z^{2,+}$. It will be useful to
consider a different set of generators $T_{\a}$ of
$\mathbf{U}^+_{\E}$, uniquely determined by the collection of formal
relations
\begin{equation}\label{E:Talpha}
1+\sum_{l \geq 1} \mathbf{1}^{ss}_{l\a}s^l=exp\left(\sum_{l \geq 1} \frac{T_{l\a}}{[l]}s^l\right),
\end{equation}
for any $\a=(r,d)$ with $r$ and $d$ relatively prime.

\vspace{.1in}

To a slope $\mu\in \mathbb{Q} \cup \{\infty\}$ 
is naturally associated the subalgebra
$\mathbf{U}_{\E}^{+,(\mu)} \subset \mathbf{U}^+_{\E}$ 
generated by $\{\mathbf{1}^{ss}_{\a}; \mu(\a)=\mu\}$.
Of course, we have $\mathbf{U}_{\E}^{+,(\mu)} \subset 
\mathbf{H}^{(\mu)}_{\E}$.

\begin{prop}[{[BS, Theorem~4.1]}] \label{P:Macd}The following hold~:
\begin{enumerate}
\item[i)] The multiplication map induces an isomorphism 
$\vec{\bigotimes}_{\mu} \mathbf{U}^{+,(\mu)}_{\E} 
\stackrel{\sim}{\to} \mathbf{U}^+_{\E}$,
\item[ii)] For any $\x=(r,d) \in \Z^{2,+}$ 
for which $r$ and $d$ are relatively prime, the assignment
$T_{l\x}/[l]\mapsto p_l/l$ extends to an isomorphism of algebras
$\mathbf{U}^{+,(\mu(\x))}_{\E} \stackrel{\sim}{\to}
(\Lambda^+_v)_{|v=\nu}$. In particular,
$\mathbf{U}^{+,(\mu(\x))}_{\E}$ is a free commutative polynomial
algebra in the generators $\{T_{l\x};l \geq 1\}$.
\end{enumerate}
\end{prop}

\vspace{.2in}

\paragraph{\textbf{1.5.}} We now wish to give a presentation of $\mathbf{U}^+_{\E}$ by generators and relations.
In fact, we will give such a presentation for the \textit{Drinfeld
double} of $\mathbf{U}^+_{\E}$, which is a more symmetric object.
Recall that if $H$ is a bialgebra equipped with a Hopf pairing
$(\;,\;)$ then its Drinfeld double $DH$ is the algebra generated by
two copies $H^+$ and $H^-$ of $H$ subject to the collection of
relations
\begin{equation}
\sum_{i,j} (h^+)_i^{(1)} (g^-)_j^{(2)} (h_i^{(2)},g_j^{(1)})=
\sum_{i,j} (g^-)_j^{(1)}(h^+)_i^{(2)}(h_i^{(1)},g_j^{(2)})
\end{equation}
for any $g,h \in H$, where we write $h^+ \in H^+$ and $g^- \in H^-$
for the corresponding elements, and where use Sweedler's notation
${\Delta}(x)=\sum_i x_i^{(1)} \otimes x_i^{(2)}$.

\vspace{.05in}

By \cite[Proposition 4.2]{BS}, the algebra $\mathbf{U}^+_{\E}$ is a
subbialgebra of $\mathbf{H}_{\E}$, and we denote by
$\mathbf{U}_{\E}$ its Drinfeld double.

Set $\Z^{2,*}$ equal to $\Z^2 \setminus \{(0,0)\}$. For $\x,\y \in
\Z^{2,*}$ let $\boldsymbol{\Delta}_{\x,\y}$ stand for the triangle
with vertices $\mathbf{o}, \mathbf{x}, \mathbf{x+y}$, where
$\mathbf{o}$ denotes the origin in $\Z^2$.  If $\x=(r,d) \in
\Z^{2,*}$ we write $d(\x)=gcd(r,d)$. For a pair of
non-colinear vectors $(\x,\y) \in \Z^{2,*}$ we set
$\epsilon_{\x,\y}$ equal to $sign(det(\x,\y))$.

We set $\mathcal{A}=\C[v^{\pm 1}, t^{\pm 1}]$ and
$\mathcal{K}=\C(v,t)$.

\vspace{.1in}

\paragraph{\textbf{Definition.}} For $i \in \N$, put
$\cc_i=(v^i+v^{-i}-t^{i}-t^{-i})[i]/i \in \mathcal{A}$. Let
$\mathbf{A}_{\mathcal{K}}$ be the unital $\mathcal{K}$-algebra
generated by elements ${t}_{\x}$ for $\x \in \Z^{2,*}$ bound by the
following set of relations.
\begin{enumerate}
\item[i)] If $\x,\x'$ belong to the same line in $\Z^2$ then
$[{t}_\x,{t}_{\x'}]=0$.
\item[ii)] Assume that $\x,\y$ are such that $d(\x)=1$ and that
$\boldsymbol{\Delta}_{\x,\y}$ has no interior lattice point. Then
$$[{t}_\y,{t}_{\x}]=\epsilon_{\x,\y}\cc_{d(\y)}\frac{{\theta}_{\x+\y}}{v^{-1}-v}$$
where the elements ${\theta}_{\z}$, $\z \in \Z^2$ are obtained by equating the Fourier coefficients of the collection of relations
$$
\sum_i {\theta}_{i\x_0}s^i=exp((v^{-1}-v)\sum_{i \geq 1}{t}_{i\x_0}s^i),
$$
for any $\x_0 \in \Z^2$ such that $d(\x_0)=1$.
\end{enumerate}

\vspace{.1in}

The algebra $\mathbf{A}_{\mathcal{K}}$ is $\Z^2$-graded by
$deg(t_{\x})=\x$. Put $\tilde{t}_{\x}=t_{\x}/[d(\x)]$ and let
$\mathbf{A}_{\mathcal{A}}$ be the unital $\mathcal{A}$-subalgebra of
$\mathbf{A}_{\mathcal{K}}$ generated by $\{\tilde{t}_{\x};\x \in
\Z^{2,*}\}$. We will write $\mathbf{A}^\pm_{\mathcal{A}}$ for the
subalgebra of $\mathbf{A}$ generated by $\{\tilde{t}_{\x}; \x
\in \pm \Z^{2,+}\}$.  By \cite{BS}, Proposition~5.1, the
multiplication yields an isomorphism
$\mathbf{A}^-_{\mathcal{A}} \otimes_{\mathcal{A}}
\mathbf{A}^+_{\mathcal{A}} \simeq \mathbf{A}_{\mathcal{A}}$. Let
$\mathbf{A}^{++}_{\mathcal{A}}$ be the subalgebra generated by
$\{\tilde{t}_{\x}; \x \in \Z^{2,++}\}$, where
$\Z^{2,++}=\{(r,d);r \geq 0, d \geq 0\}$. We have
$$\mathbf{A}^{+}_{\mathcal{A}}=\bigoplus_{\x \in \Z^{2,+}} \mathbf{A}_{\mathcal{A}}[\x], \qquad
\mathbf{A}^{++}_{\mathcal{A}}=\bigoplus_{\x \in \Z^{2,++}}
\mathbf{A}_{\mathcal{A}}[\x].$$
The algebra $\mathbf{A}_{\mathcal{A}}$ has an obvious symmetry~: the
group $SL(2, \Z)$ acts by automorphisms such that 
$g\cdot\tilde{t}_{\x}=\tilde{t}_{g(\x)}$.

\vspace{.2in}

Let $\sigma$, $\overline{\sigma}$ be the two eigenvalues of the
Frobenius endomorphism acting on the vector space
$H^1(\E\otimes\overline{\mathbb{F}}_l,\qlb)$, with $p$ prime to $l$.
We'll fix once for all a field isomorphism
$\mathbb{C}\simeq\qlb$. This allows us to view
$\sigma$, $\overline{\sigma}$ as complex numbers. 
Let $\mathbf{A}_{\E}$ stand
for the specialization of $\mathbf{A}_{\mathcal{A}}$ at
$v=\nu=l^{-1/2}$ and $t=\sigma\nu$.

\begin{theo}[{[BS, Theorem~5.1]}]\label{T:BS}  
The assignment $\tilde{t}_{\x} \mapsto T_{\x}/[deg(\x)]$,
$\x \in \Z^{2,*}$, extends to an isomorphism $\mathbf{A}_{\E}
\stackrel{\sim}{\to} \mathbf{U}_{\E}$. It restricts to an
isomorphism $\mathbf{A}^+_{\E}  \stackrel{\sim}{\to}
\mathbf{U}^+_{\E}$.
\end{theo}

\vspace{.1in}

To a slope $\mu\in \mathbb{Q} \cup \{\infty\}$ 
is naturally associated the subalgebra
$\mathbf{A}_{\mathcal{A}}^{+,(\mu)} \subset \mathbf{A}^+_{\mathcal{A}}$ 
generated by $\{T_{\a}; \mu(\a)=\mu\}$.

\vspace{.2in}

\vspace{.2in}

\section{Double affine Hecke algebras}

\vspace{.2in}

\paragraph{\textbf{2.1.}} 
We set $\mathcal{A'}=\C[v^{\pm 1}, q^{\pm 1}]$ and
$\mathcal{K'}=\C(v,q)$.
The double affine Hecke algebra $\ddot{\H}_n$ of $GL(n)$, 
abreviated DAHA, is the
$\mathcal{K'}$-algebra generated by elements 
$T_i^{\pm 1}$, $X_j^{\pm 1}$ and $Y_j^{\pm 1}$
for $1 \leq i \leq n-1$ and $1 \leq j \leq n$, 
subject to the following relations~:

\begin{equation}\label{E:finiteHecke}
(T_i+v^{-1})(T_i-v)=0, \qquad T_i T_{i+1}T_i=T_{i+1}T_iT_{i+1}
\end{equation}
\begin{equation}
T_iT_k=T_kT_i\;\qquad  \text{if}\; |i-k| >1
\end{equation}
\begin{equation}\label{E:def-Hecke1}
X_jX_k=X_kX_j, \qquad Y_jY_k=Y_kY_j
\end{equation}
\begin{equation}\label{E:def-Hecke2}
T_iX_iT_i=X_{i+1}, \qquad T_i^{-1}Y_iT_i^{-1}=Y_{i+1}
\end{equation}
\begin{equation}\label{E:def-Hecke3}
T_iX_k=X_kT_i, \qquad T_i Y_k =Y_kT_i \qquad \text{if\;} |i-k|>1
\end{equation}
\begin{equation}\label{E:def-Hecke4}
Y_1 X_1 \cdots X_n=qX_1\cdots X_nY_1
\end{equation}
\begin{equation}\label{E:def-Hecke5}
X_1^{-1}Y_2=Y_2X_1^{-1}T_1^{-2}
\end{equation}

The subalgebra $\H_n$ generated by $\{T_i\}$ is the usual Hecke
algebra of the symmetric group $\mathfrak{S}_n$, while the
subalgebras $\dot{\H}_{n, X}$ and $\dot{\H}_{n,Y}$ respectively
generated by $\H_n$ and $\{X_i^{\pm 1}\}$, and $\H_n$ and
$\{Y_i^{\pm 1}\}$, are both isomorphic to the Hecke algebra of the
affine Weyl group $\widehat{\mathfrak{S}}_n \simeq \mathfrak{S}_n
\ltimes \Z^n$.  We define a $\Z^2$-grading on $\ddot{\H}_n$ by
giving $T_i$, $X_i$ and $Y_i$ degrees $0, (1,0)$ and $(0,1)$
respectively.

\vspace{.1in}

Let $s_i \in \mathfrak{S}_n$ denote the transposition $(i,i+1)$, and
let $l: \mathfrak{S}_n \to \mathbb{N}$ be the standard length
function. If $w =s_{i_1} \cdots s_{i_r}$ is a reduced decomposition
of $w \in \mathfrak{S}_n$, we set $T_w=T_{i_1} \cdots T_{i_r}$. We
put $\tilde{S}=\sum_{w \in \mathfrak{S}_n} v^{l(w)} T_w$. We have
$\tilde{S}^2=[n]^+! \tilde{S}$, so that the element
${S}=\tilde{S}/[n]^+!$ is idempotent. For any $i$ we have
$T_iS=S T_i=vS$.

\vspace{.1in}

We will mainly be interested in the \textit{spherical} DAHA
of $\ddot{\H}_n$ equal to $\mathbf{S}\ddot{\H}_n={S}\ddot{\H}_n{S}$.
Before we can give some bases for $\mathbf{S}\ddot{\H}_n$
we need a few notations.
Set $\R_n$, $\Fb_n$ equal to the algebras $$\C[ x_1^{\pm 1},
\ldots,x_n^{\pm 1}, y_1^{\pm 1}, \ldots y_n^{\pm
1}]^{\mathfrak{S}_n},\quad\C\langle x_1^{\pm 1}, \ldots,x_n^{\pm 1},
y_1^{\pm 1}, \ldots y_n^{\pm 1}\rangle^{\mathfrak{S}_n}.$$ The
algebra $\R_n$ consists of the symmetric Laurent polynomials where
the $x$'s and the $y$'s commute between themselves and with each
other. The algebra $\Fb_n$ consists of the symmetric Laurent
polynomials where the $x$'s and the $y$'s commute between
themselves, but not with each other. Here the symmetric group
$\mathfrak{S}_n$ acts by simultaneous permutation on the $x$'s and
the $y$'s. There is an obvious projection map
$Com~:\Fb_n\to\R_n.$
The size of the spherical DAHA is described in the following result.

\vspace{.1in}

\begin{prop}\label{P:Weyl?} Let $\{E_t\}$ be any collection of elements of
$\Fb_n$ such that $\{Com(E_t)\}$ forms a basis of $\R_n$. Then
$\{SE_tS\}$ is a $\mathcal{K'}$-basis of $\SH_n$.\end{prop}

\noindent \textit{Proof.}  Let $\ddot{\H}_{\mathcal{A},n}$ stand for
the $\mathcal{A'}$-subalgebra of $\ddot{\H}_n$ generated
by $T_i$ for $i=1, \ldots, n-1$ and $X_j^{\pm 1}, Y_j^{\pm 1}$ for
$j=1, \ldots, n$. We also set $\SH_{\mathcal{A},n}=S
\ddot{\H}_{\mathcal{A},n} S$. Both $\ddot{\H}_{\mathcal{A},n}$ and
$\SH_{\mathcal{A},n}$ are free $\mathcal{A}$-modules, see
\cite{Cherednik}. Write $\C$ for the one-dimensional complex
representation of $\mathcal{A'}$ in which $v$ and $q$
act as $1$. The assignment
$$S P(X_i,Y_i)S\otimes 1 \mapsto \frac{1}{n!} 
\sum_{\sigma \in \mathfrak{S}_n} \sigma \cdot P(x_i,y_i)$$
give rise to a $\C$-algebra isomorphism
$$\pi:~\SH_{\mathcal{A},n} \otimes \C \stackrel{\sim}{\to} \R_n.$$
Now let $\{E_t\}$ be as in the hypothesis of the proposition. Then
$SE_tS \in \SH_{\mathcal{A},n}$ for all $t$ and $\{\pi(SE_tS)\}$
forms a $\C$-basis of $ \R_n$. It follows that $\{SE_tS\}$ are
$\mathcal{K'}$-linearly independent and
generate $\SH_{n}$ over $\mathcal{K'}$. \qed

\vspace{.1in}

There is an action of the braid group $B_3$ on three strands by
automorphisms on $\ddot{\H}_n$, explicitly given by the following
operators.

\begin{equation*}
\rho_1: \begin{cases}
T_i \mapsto T_i,&\\
X_i \mapsto X_iY_i(T_{i-1} \cdots T_i)(T_i \cdots T_{i-1}), & \\
Y_i  \mapsto Y_i ,
\end{cases}
\end{equation*}

\begin{equation*}
\rho_2: \begin{cases}
T_i \mapsto T_i,&\\
Y_i \mapsto Y_iX_i(T^{-1}_{i-1} \cdots T^{-1}_i)(T^{-1}_i \cdots T^{-1}_{i-1}),& \\
X_i  \mapsto X_i.&
\end{cases}
\end{equation*}
These operators preserve $\mathbf{S}\ddot{\H}_n$, and the
corresponding $B_3$-action factors through an $SL(2,\Z)$-action
$\rho: SL(2,\Z) \to \mathrm{Aut}(\mathbf{S}\ddot{\H}_n)$ satisfying
$\rho(A_1)=\rho_1$, $\rho(A_2)=\rho_2$, where $A_1=\tiny{\begin{pmatrix} 1
& 0 \\1 & 1
\end{pmatrix}}$ and $A_2=\tiny{\begin{pmatrix} 1 & 1 \\
0 & 1 \end{pmatrix}}.$

\vspace{.1in}

The following technical lemma will be often used.

\vspace{.1in}

\begin{lem}\label{L:1} For each $l \geq 2$ we put 
$\a_l=T_{l-1}^{-1} \cdots T_2^{-1}T_1^{-2} T_2 \cdots T_{l-1}$.
The following relations hold
\begin{equation}\label{E:L11}
X_l^{-1}Y_1X_l=\a_lY_1,
\end{equation}
\begin{equation}\label{E:L12}
Y_lX_1=X_1Y_l+(v^{-1}-v)T_{l-1}^{-1} \cdots T_2^{-1}T_1^{-1}T_2^{-1} \cdots T_{l-1}^{-1}Y_1X_1,
\end{equation}
\begin{equation}\label{E:L13}
qX_1Y_1=T_1^{-1} \cdots T_{n-2}^{-1}T_{n-1}^{-2} T_{n-2}^{-1} \cdots T_1^{-1}Y_1X_1.
\end{equation}
\begin{equation}\label{E:L14}
\a_2 \cdots \a_l =T_1^{-1} \cdots T_{l-2}^{-1}T_{l-1}^{-2} T_{l-2}^{-1} \cdots T_1^{-1},
\end{equation}
\end{lem}
\noindent
\textit{Proof.} By definition, we have $T_1^{-2}Y_1X_2^{-1}=X_2^{-1}Y_1$, which is relation (\ref{E:L11}) for $l=2$. Multiplying on the left and on the right by $T_2^{-1}$ and using the fact that $[T_2,Y_1]=0$ we obtain
$$T_2^{-1}T_1^{-2} T_2Y_1 \cdot T_2^{-1}X_2^{-1}T_2^{-1} =T_2^{-1} X_2^{-1}T_2^{-1}Y_1=Y_1X_3^{-1}$$
Since $T_2^{-1}X_2^{-1}T_2^{-1}=X_3^{-1}$, we get
$$T_2^{-1}T_1^{-2}T_2 Y_1 X_3^{-1}=X_3^{-1}Y_1$$
which is (\ref{E:L11}) for $l=3$. A similar reasoning with multiplying on the left and on the right by
$T_3^{-1}$ yields (\ref{E:L11}) for $l=4$, etc.

We now prove (\ref{E:L12}). From the defining relations of $\ddot{\H}_n$ we have
\begin{equation*}
\begin{split}
Y_2X_1&=X_1 Y_2 X_1^{-1}T_1^{-2}X_1= X_1 Y_2 + (v^{-1}-v) X_1Y_2 X_1^{-1}T_1^{-1}X_1\\
&=X_1Y_2+(v^{-1}-v)X_1Y_2X_1^{-1}T_1^{-1}X_1 T_1 T_1^{-1}\\
&=X_1Y_2+(v^{-1}-v)X_1Y_2X_1^{-1}T_1^{-2}X_2 T_1^{-1}=X_1Y_2+(v^{-1}-v)Y_2X_2T_1^{-1}\\
&=X_1Y_2+(v^{-1}-v)T_1^{-1}Y_1X_1,
\end{split}
\end{equation*}
which is (\ref{E:L12}) for $l=2$. Now we multiply on the left and on the right by $T_2^{-1}$ and usi the fact that $[T_2,X_1]=[T_2,Y_1]=0$ to get
$$T_2^{-1}Y_2 T_2^{-1}X_1= X_1 T_2^{-1} Y_2 T_2^{-1} + (v^{-1}-v) T_2^{-1}T_1^{-1}T_2^{-1}Y_1X_1$$
which, by virtue of the relation $T_2^{-1}Y_2T_2^{-1}=Y_3$, gives (\ref{E:L12}) for $l=3$. To obtain (\ref{E:L12}) for $l=3,4$ etc., we successively multiply on the left and on the right by $T_3^{-1}, T_4^{-1}$, etc.

We turn to (\ref{E:L13}). Recall that by definition
$$(X_n^{-1} \cdots X_2^{-1} Y_1 X_2 \cdots X_n) X_1=qX_1Y_1.$$
By (\ref{E:L11}) above we have $X_2^{-1}Y_1X_2=\a_2Y_1$. Since $[\a_l,X_k]=0$ if $k>l$, we obtain after conjugation by $X_3$,
$$X_3^{-1}X_2^{-1}Y_1X_2X_3=\a_2X_3^{-1}Y_1X_3=\a_2\a_3Y_1.$$
Continuing in this manner yields in the end
$$X_n^{-1} \cdots X_2^{-1}Y_1 X_2 \cdots X_n=\a_2 \cdots \a_n Y_1.$$
Thus (\ref{E:L13}) will follow from (\ref{E:L14}) which we now prove. It is easy to check (\ref{E:L14})
for $l=2$ and $l=3$. We argue by induction on $l$. So we fix $l$ and assume that
$$\a_2 \cdots \a_l=T_1^{-1} \cdots T_{l-2}^{-1}T_{l-1}^{-2} T_{l-2}^{-1} \cdots T_1^{-1}.$$
Then
\begin{equation*}
\begin{split}
\a_2& \cdots \a_{l+1}=\\
&=T_1^{-1} \cdots T_{l-2}^{-1}T_{l-1}^{-2} T_{l-2}^{-1} \cdots T_1^{-1}\cdot T_{l}^{-1} \cdots T_2^{-1}T_1^{-2} T_2 \cdots T_l\\
&=T_1^{-1} \cdots T_{l-2}^{-1}T_{l-1}^{-2} T_{l-2}^{-1} \cdots T_2^{-1}\cdot T_{l}^{-1} \cdots T_3^{-1} T_1^{-1}T_2^{-1}T_1^{-2} T_2T_3 \cdots T_l\\
&=T_1^{-1} \cdots T_{l-2}^{-1}T_{l-1}^{-2} T_{l-2}^{-1} \cdots T_2^{-1}\cdot T_{l}^{-1} \cdots T_3^{-1} T_2^{-1}T_1^{-1}T_2^{-1}T_1^{-1} T_2T_3 \cdots T_l\\
&=T_1^{-1} \cdots T_{l-2}^{-1}T_{l-1}^{-2} T_{l-2}^{-1} \cdots T_2^{-1}\cdot T_{l}^{-1} \cdots T_3^{-1} T_2^{-2}T_1^{-1}T_3 \cdots T_l\\
&=T_1^{-1} \cdots T_{l-2}^{-1}T_{l-1}^{-2} T_{l-2}^{-1} \cdots T_2^{-1}\cdot T_{l}^{-1} \cdots T_3^{-1} T_2^{-2}T_3 \cdots T_lT_1^{-1}.
\end{split}
\end{equation*}
This last expression is of the form $T_1^{-1} Z T_1^{-1}$ where $Z$ is (using the induction hypothesis) equal to $\a_2 \ldots \a_{l}$ \textit{for the subalgebra of $\H_n$ generated by $T_2, \ldots, T_l$}. In particular, $T_1$ is not involved in $Z$. By our induction hypothesis again, we deduce that
$Z=T_2^{-1} \cdots T_{l-1}^{-1}T_l^{-2}T_{l-1}^{-1} \cdots T_2^{-1}$ from which
$$\a_2 \cdots \a_{l+1}=T_1^{-1} \cdots T_{l-1}^{-1}T_l^{-2}T_{l-1}^{-1} \cdots T_1^{-1}$$
of course follows. The lemma is proved.\qed

\vspace{.2in}

\paragraph{\textbf{2.2.}} For $e >0$ we set
$P^n_{(0,e)}={S} \sum_i Y_i^e {S}.$
More generally, if $(r,d) =g \cdot (0,e)$ 
we put $P^n_{(r,d)}=\rho(g) P^n_{(0,e)}$. If the element $g' \in
SL(2,\Z)$ fixes the couple $(0,e)$ then $\rho(g')=\rho_1^l$ for some
$l$ hence $\rho(g') P^n_{(0,e)}=P^n_{(0,e)}$. Therefore the above
definition makes sense, and it yields an element
$P^n_{\x}\in\mathbf{S}\ddot{\H}_n$ for each $\x \in \Z^{2,*}$, such that
$\rho(g) P^n_{\x}=P^n_{g(\x)}$
for any $g \in SL(2, \Z)$. \vspace{.1in}

As an illustration, let us give the expression for certain elements $P^n_{(r,d)}$ when $r,d$ are relatively prime.
To unburden the notation, we drop the exponent $n$ in $P^n_{(r,d)}$.

\vspace{.1in}

\begin{lem}\label{L:2} For any $l \in \Z$ we have
\begin{equation}\label{E:L21}
P_{(l,1)}=[n]^-SY_1X_1^lS,
\end{equation}
\begin{equation}\label{E:L22}
P_{(1,l)}=q[n]^+SX_1Y_1^lS,
\end{equation}
\begin{equation}\label{E:L23}
P_{(0,-1)}=q[n]^+SY_1^{-1}S,
\end{equation}
\begin{equation}\label{E:L24}
P_{(l,-1)}=q[n]^+SX_1^lY_1^{-1}S.
\end{equation}
\end{lem}
\noindent
\textit{Proof.} Observe that since $SY_{i+1}S=ST_{i}^{-1}Y_iT_{i}^{-1}S=v^{-2}SY_iS$, we have
$P_{(0,1)}=S \sum_i Y_i S=[n]^- SY_1S$. Equation (\ref{E:L21}) follows from this and an application
of $\rho_2^l$. In particular, we have $P_{(1,1)}=[n]^-SY_1X_1S$. Using Lemma~\ref{L:1}, (\ref{E:L13}), we obtain $P_{(1,1)}=qv^{2(n-1)}[n]^-SX_1Y_1S=q[n]^+SX_1Y_1S$. Application of $\rho_1^l$ yields (\ref{E:L22}). The last two equalities are proved using similar techniques.\qed

\vspace{.1in}

The values of $P_{(r,d)}$ when $r$ and $d$ are not relatively prime
is usually harder to compute. We give a few examples, which will be
important for us.

\vspace{.1in}

\begin{lem}\label{L:3} For any $l \geq 1$, we have
\begin{equation}\label{E:L31}
P_{(-l,0)}=S\sum_i X_i^{-l} S,
\end{equation}
\begin{equation}\label{E:L32}
P_{(0,-l)}=q^lS\sum_i Y_i^{-l} S.
\end{equation}
\begin{equation}\label{E:L33}
P_{(l,0)}=q^lS\sum_i X_i^{l} S,
\end{equation}
\end{lem}
\noindent \textit{Proof.} Set $A_3=\tiny{\begin{pmatrix}1 & -1\\1 &
0\end{pmatrix}}$ and $A_4=\tiny{\begin{pmatrix}-1 & 0\\2 &
-1\end{pmatrix}}$. We have
$$\rho(A_3)Y_1=\rho_1 \rho_2^{-1} \cdot Y_1=X_1^{-1}$$
and hence $\rho(A_3)Y_i=X_i^{-1}$ for all $i$. The first relation
(\ref{E:L31}) immediately follows. The proof of the second and the
third relations are identical, and we only treat (\ref{E:L32}). We
have, using Lemma~\ref{L:1}, (\ref{E:L13})
\begin{equation*}
\begin{split}
\rho(A_4) Y_1&=\rho_2^{-1}\rho_1^2 \rho_2^{-1} \cdot Y_1=X_1^{-1}Y_1^{-1}X_1\\
&=qY_1^{-1}T_1 \cdots T_{n-2}T_{n-1}^2 T_{n-2} \cdots T_1\\
&=q T_1^{-1} \cdots T_{n-1}^{-1} Y_n^{-1} T_{n-1} \cdots T_1,
\end{split}
\end{equation*}
where in the last equation we have the relations
$Y^{-1}_iT_i=T_i^{-1}Y_{i+1}$. It follows that
\begin{equation}
\rho(A_4) (Y^l_1)= q^lT_1^{-1} \cdots T_{n-1}^{-1} Y_n^{-l} T_{n-1}
\cdots T_1
\end{equation}
for any $l$. Hence $\rho(A_4)(SY^l_1S)=q^lS Y_n^{-l} S$. It is
well-known and easy to show that the elements $SY_1^lS$ for $l=1,
\ldots, n$ freely generate the ring $S\C[Y_1, \ldots, Y_n]S$. Let
$$\theta~: \C[SY_1S,\ldots ,SY_1^nS] \stackrel{\sim}{\to} S\C[Y_1, \ldots, Y_n]S =
S \C[Y_1, \ldots, Y_n]^{\mathfrak{S}_n}S$$
and let
$$\theta'~: \C[SY_n^{-1}S,\ldots ,SY_n^{-n}S] \stackrel{\sim}{\to} S\C[Y_1^{-1}, \ldots, Y_n^{-1}]S =
S \C[Y_1^{-1}, \ldots, Y_n^{-1}]^{\mathfrak{S}_n}S$$
be defined in a similar fashion. Equation (\ref{E:L32}) is a
consequence of the following result.

\vspace{.05in}

\begin{slem} The composition $u=\theta'\circ\rho(A_4)\circ\theta^{-1}$ satisfies
$$u(SP(Y_1, \ldots, Y_n)S)=q^lSP(Y_1^{-1}, \ldots, Y_n^{-1})S$$
for any symmetric polynomial $P(t_1, \ldots, t_n)$.
\end{slem}
\noindent \textit{Proof of sublemma.} Let $\dot{\H}_{n,Y}^+$ (resp.
$\dot{\H}^-_{n,Y}$) be the subalgebras generated by $\H_n$ and the
elements $Y_1, \ldots, Y_n$ (resp. by $\H_n$ and the elements
$Y_1^{-1}, \ldots, Y_n^{-1}$). The assignment
$T_i\mapsto T_{n-i}$, $Y_i\mapsto Y_{n+1-i}^{-1}$
gives rise to an isomorphism of algebras
$\Theta~: \dot{\H}^+_{n,Y} \stackrel{\sim}{\to} \dot{\H}^-_{n,Y}$.
It restricts to an isomorphism of spherical algebras
$S\Theta~:\mathbf{S}\dot{\H}^+_{n,Y} \stackrel{\sim}{\to} 
\mathbf{S}\dot{\H}^-_{n,Y}$. This last map clearly satisfies
$$S\Theta ( SP(Y_1, \ldots, Y_n)S)=SP(Y_1^{-1}, \ldots, Y_n^{-1})S$$
for any symmetric polynomial $P(t_1, \ldots, t_n)$. 
It remains to observe that $u$
coincides with $q^lS\Theta$ on the elements $SY_1^lS$, 
and that these elements generate
$S\dot{\H}^+_{n,Y}S$, so that in fact $u=q^lS\Theta$. 
The sublemma is proved.\qed

\vspace{.2in}

\begin{prop} The elements $\{P^n_{\x}; \x \in \Z^{2,*}
\}$ generate $\SH_n$ as a $\mathcal{K'}$-algebra. \end{prop} \noindent
\textit{Proof.} As in the proof of Proposition~\ref{P:Weyl?}, let
$\SH_{\mathcal{A},n}$ be the integral form of $\SH_n$ and let $\pi:
\SH_{\mathcal{A},n} \to \R_n$ be the specialization at $q=v=1$. It
is easy to see that $\pi(P^n_{(r,d)})=\sum_i x_i^ry_i^d$. By Weyl's
theorem, the elements $\{\pi(P^n_{(r,d)}); (r,d) \in \Z^{2,*}\}$
generate the ring $\R_n$. It follows that $\{ P^n_{(r,d)}; (r,d)
\in \Z^{2,*}\}$ generate the algebra $\SH_n$ over the field
$\mathcal{K'}$. \qed

\vspace{.2in}

\section{The projection map}

\vspace{.1in}

Recall that $\mathcal{K}=\C(v,t)$,
$\mathcal{K'}=\C(v,q)$,
and $\mathbf{A}_{\mathcal{K}}=\mathbf{A}_{\mathcal{A}} \otimes \mathcal{K}$.
The first main Theorem of this paper is the following.

\begin{theo}\label{T:1}  For any $n >0$, the assignment
$$t \mapsto qv^{-1},\quad
\tilde{t}_{\x} \mapsto
\frac{1}{d(\x)}(q^{-d(\x)}v^{d(\x)}-v^{-d(\x)}){P}^n_{\x}\
\text{for}\
 \x \in \Z^{2,*}$$
extends to a surjective $\C$-algebra homomorphism 
$\Psi_n:\mathbf{A}_{\mathcal{K}} \tto \mathbf{S}\ddot{\H}_n$.
\end{theo}
\noindent
\textit{Proof.} Fix an integer $n$.
For simplicity we drop the index $n$ in all notations.
We have to show that the elements
$[d(\x)](q^{-d(\x)}v^{d(\x)}-v^{-d(\x)}){P}^n_{\x}/d(\x)$
satisfy relations i) and ii) of Section~1.4.
Relation i) is clear for $\x=(0,r), \x'=(0,r')$ with $r,r'$ of the same sign,
and follows from Lemma~\ref{L:3}, (\ref{E:L32}) for $r,r'$ of different signs.
By applying a suitable automorphism $\rho_g$ for $g \in SL(2,\Z)$,
we deduce relation i) for any other line in $\Z^2$ through the origin.

The proof of relation ii) is much more involved.
We reduce it to two sets of equalities
(\ref{E:P1}) and (\ref{E:T28}) which are dealt with in
Appendices A and B respectively.
Note that the assignment $\Psi$ respects the $SL(2,\Z)$ action
on both sides.
Hence it is enough to check the relation ii) for one pair in each
orbit under this $SL(2,\Z)$-action. By the same argument (based on
Pick's formula) as in \cite{BS}, Theorem~5.1, we may reduce
ourselves to the cases where $\x=(1,0), \y=(0,l)$ with $l \in \Z^*$,
or $\x=(0,1), \y=(l,-1)$ with $l \geq 0$.

\vspace{.1in}

\noindent
\textit{Case a1.} Assume that $\x=(1,0)$ and $\y=(0,l)$ with $l >0$. We have to show that
\begin{equation}\label{E:T21}
[\Psi({t}_{(1,0)}),\Psi({t}_{(0,l)})]=-\cc_l \Psi({t}_{(1,l)}),
\end{equation}
which we may rewrite as
\begin{equation}\label{E:T22}
\frac{[l]}{l}(q^{-1}v-v^{-1})(q^{-l}v^l-v^{-l})[P_{(1,0)},P_{(0,l)}]=-\cc_l(q^{-1}v-v^{-1})P_{(1,l)}.
\end{equation}
By Lemma~\ref{L:2}, we have
$$P_{(1,0)}=q[n]^+SX_1S, \qquad P_{(0,l)}=\sum_iSY_i^lS,\qquad P_{(1,l)}=q[n]^+SX_1Y_1^lS.$$
Using this and the identity
$[l](1-q^l)(q^{-l}v^l-v^{-l})/l=-\cc_l$, we see that
(\ref{E:T22}) is equivalent to the following proposition.

\vspace{.1in}

\begin{prop}\label{P:1} For any $l >0$ we have
\begin{equation}\label{E:P1}
\bigg[SX_1S,\sum_iSY_i^lS\bigg]=S\bigg[X_1,\sum_i Y_i^l\bigg]S=
(1-q^l)SX_1Y_1^lS.
\end{equation}
\end{prop}

\vspace{.1in}

\noindent
\textit{Case a2.} Let us now assume that $\x=(1,0)$ and $\y=(0,-l)$
with $l >0$. We have to show that
\begin{equation}\label{E:T23}
[\Psi({t}_{(1,0)}),\Psi({t}_{(0,-l)})]=\cc_{l} \Psi({t}_{(1,-l)}),
\end{equation}
which, after using the definitions and Lemmas~\ref{L:2} and~\ref{L:3},
reduces to
\begin{equation}\label{E:T24}
\bigg[SX_1S,\sum_iSY_i^{-l}S\bigg]=(1-q^{-l})SX_1Y_1^{-l}S.
\end{equation}
Consider the $\C$-algebra isomorphism
$\sigma:\ddot{\H}_n\to\ddot{\H}_n$ given by
$$T_i \mapsto T_i^{-1}, \qquad X_i \mapsto Y_i,
\qquad Y_i \mapsto X_i, \qquad v \mapsto v^{-1}, \qquad q \mapsto q^{-1},$$
see \cite{Cherednik}. Applying $\sigma$ to (\ref{E:T24}) gives the equation
$$\bigg[SY_1S,\sum_iSX_i^{-l}S\bigg]=(1-q^{l})SY_1X_1^{-l}S,$$
which, once transformed by the automorphism $\rho(A_5)$,
$A_5=\tiny{\begin{pmatrix} 0 & 1\\ -1 & 0 \end{pmatrix}}$,
is nothing else than (\ref{E:P1}).
Thus this case also follows from Proposition~\ref{P:1}.

\vspace{.1in}

\noindent
\textit{Case b.} The final case to consider is that of $\x=(0,1)$ and $\y=(-1,l)$ with $l >0$.
Here we have to prove that
\begin{equation}\label{E:T25}
[\Psi(t_{(1,0)}),\Psi_(t_{(-1,l)})]=\cc_1 \frac{\Psi(\theta_{(0,l)})}{v^{-1}-v},
\end{equation}
which reduces to
\begin{equation}\label{E:T26}
\frac{(1-v^{2n})(1-qv^{-2})}{q-1}\bigg[\sum_iSY_i S,SX_1^lY_1^{-1}S\bigg]=\Psi(\theta_{(0,l)}).
\end{equation}
Forming a generating series, we may write this as
\begin{equation}\label{E:T27}
1+\sum_{l \geq 1}\frac{(1-v^{2n})(1-qv^{-2})}{q-1}S\bigg[ \sum_i Y_i ,X_1^lY_1^{-1}\bigg]Ss^l=
1+\sum_{l \geq 1}\Psi(\theta_{(0,l)})s^l.
\end{equation}
Given the definition of $\theta_{(0,l)}$, we finally obtain that (\ref{E:T25}) for all $l>0$
 is equivalent to the following assertion.
\vspace{.1in}

\begin{prop}\label{P:2} The following holds~:
\begin{equation}\label{E:T28}
\begin{split}
exp &\left(\sum_{l \geq 1}\frac{(v^{-l}-v^l)(v^l-q^lv^{-l})}{l}\sum_iSX_i^lSs^l \right)=\\
&\qquad \qquad \qquad =1+\sum_{l \geq 1} \frac{(1-qv^{-2})(1-v^{2n})}{q-1} S \bigg[\sum_i Y_i,X_1Y_1^{-1} \bigg]S s^l.
\end{split}
\end{equation}
\end{prop}

This completes the proof of the Theorem.\qed

\vspace{.2in}

\section{Stable limits of DAHA's}

\vspace{.1in}

\paragraph{\textbf{4.1.}}
In considering stable limits of DAHA's we will be concerned
with the graded subalgebras $\mathbf{S}\ddot{\H}_m^{+}$,
$\mathbf{S}\ddot{\H}_m^{++}$ of $\mathbf{S}\ddot{\H}_m$ generated by
the elements $P^m_{\x}$ for $\x \in \Z^{2,+}$, $\Z^{2,++}$
respectively. We have
$$\mathbf{S}\ddot{\H}_m^{+}=
\bigoplus_{\x \in \Z^{2,+}}\mathbf{S} \ddot{\H}^+_m[\x], \qquad
\mathbf{S}\ddot{\H}_m^{++}=
\bigoplus_{\x \in \Z^{2,++}}\mathbf{S} \ddot{\H}^{++}_m[\x].$$

\vspace{.1in}

\begin{prop}\label{P:3} The assignment
$P_{\x}^m \mapsto P_{\x}^{m-1}$ for each $\x \in \Z^{2,+}$
extends to a unique $\mathcal{K'}$-algebra morphism
$\Phi_m: \mathbf{S}\ddot{\H}^+_{m} \tto \mathbf{S}\ddot{\H}^+_{m-1}$.
A similar statement holds for $\mathbf{S}\ddot{\H}^{++}_{m}$.\end{prop}
\noindent
\textit{Proof.} The proof is based on the realization of Cherednik algebras 
as certain algebras of difference operators. 
Let $\mathbf{D}_{m}$ stand for the algebra of $q$-difference operators
$\mathcal{K'}[x_1^{\pm 1}, \partial_i^{\pm 1}, \ldots, x_m^{\pm 1}, 
\partial_m^{\pm 1}]$ with defining relations
$$[x_i,x_j]=[\partial_i,\partial_j]=0, \qquad 
\partial_ix_j=q^{\delta_{ij}}x_j\partial_i.$$
We also denote by $\mathbf{D}_{m,loc}$ the localization of
$\mathbf{D}_m$ with respect to the elements
$\{x_i-v^lq^nx_j;l, n \in \Z , \; i,j=1, \ldots, m, \}$.
The symmetric group $\mathfrak{S}_m$ acts in an obvious fashion on
$\mathbf{D}_{m,loc}$ and we may form the semidirect product
$\mathbf{D}_{m,loc}\rtimes \mathfrak{S}_m$.
The following lemma is due to Cherednik.

\vspace{.1in}

\begin{lem}[Cherednik, \cite{Cherednik}]
Set $\omega=s_{m-1} \cdots s_1 \partial_1$.
There is a unique embedding of algebras
$\varphi_m~:\ddot{\H}_m \hookrightarrow \mathbf{D}_{m,loc}\ltimes \mathfrak{S}_m$
satisfying
\begin{align*}
\varphi_m(X_i)&=x_i,\\
\varphi_m(T_i)&=vs_i + \frac{v-v^{-1}}{x_i/x_{i+1}-1}(s_i-1),\\
\varphi_m(Y_i)&=\varphi_m(T_i) \cdots \varphi_m(T_{m-1})\omega
\varphi_m(T_1^{-1}) \cdots \varphi_m(T_{i-1}^{-1}).
\end{align*}
\end{lem}

\vspace{.15in}

It is known that
$\varphi_m(\mathbf{S}\ddot{\H}_m)\subset\mathbf{D}_{m,loc}^{\mathfrak{S}_m}
\rtimes \mathfrak{S}_m$.
Composing the restriction of $\varphi_m$ to $\mathbf{S}\ddot{\H}_m$
with the projection
\begin{align*}
\mathbf{D}_{m,loc}^{\mathfrak{S}_m}\rtimes \mathfrak{S}_m\to
\mathbf{D}_{m,loc}^{\mathfrak{S}_m},\quad
P(x_i^{\pm 1}, \partial_i^{\pm 1})\sigma\mapsto
P(x_i^{\pm 1}, \partial_i^{\pm 1})
\end{align*}
provides us with an embedding
$\psi_m: \SH_m \hookrightarrow \mathbf{D}_{m,loc}^{\mathfrak{S}_m}$.
Set $\mathbf{D}^{++}_{m,loc}$ equal to the algebra
$\mathcal{K'}[x_1,\partial_1, \ldots, x_m, \partial_m]_{loc}$.

\begin{lem} We have
$\psi_m(\SH^{++}_m)\subset (\mathbf{D}_{m,loc}^{++})^{\mathfrak{S}_m}$.
\end{lem}
\noindent
\textit{Proof.} It is easy to see that $\mathbf{A}^{++}_{\mathcal{K}}$ is generated by
$\{t_{(0,l)}, t_{(l,0)}; l \geq 1\}$. Hence, by theorem~\ref{T:1},
$\SH_m^{++}$ is generated by $\{P^m_{(0,l)}, P^m_{(l,0)}; l \geq 1\}$,
and it suffices to check the veracity of the lemma for these elements, for which it is obvious.\qed

\vspace{.1in}

We now consider the map
$$\pi_m:(\mathbf{D}_{m,loc}^{++})^{\mathfrak{S}_m} \to (\mathbf{D}_{m-1,loc})^{\mathfrak{S}_{m-1}}$$
defined by sending $x_l$, $\partial_l$ to $x_l$,
$v^{-1} \partial_l$ if $l<m$, and $x_m$, $\partial_m$ to zero.
This is a well-defined algebra homomorphism.
We may summarize the situation in the following diagram
of algebra homomorphisms

$$\xymatrix{
\SH^{++}_m \ar[r]^-{\psi_m} & (\mathbf{D}_{m,loc}^{++})^{\mathfrak{S}_m} \ar[d]^-{\pi_m}\\
\SH^{++}_{m-1} \ar[r]^-{\psi_{m-1}} &
(\mathbf{D}_{m-1,loc})^{\mathfrak{S}_{m-1}}}$$ in which $\psi_m$ and
$\psi_{m-1}$ are embeddings. Therefore, Proposition~\ref{P:3} will
be proved for the algebra $\SH^{++}_m$ once we show that
\begin{equation}\label{E:P31}
\pi_m \circ \psi_m(P^m_{\x})=\psi_{m-1}(P^{m-1}_{\x}),\quad
\x \in \Z^{2,++}.
\end{equation}

\vspace{.1in}

\begin{lem}\label{L:P31} For any $\x \in \Z^{2,++}$
there exists a polynomial $Q_{\x}\in\Fb_\infty$ such that
the following formula holds in $\mathbf{S}\ddot{\H}_m$ for any $m$
$$Q_{\x}(P^m_{(0,1)}, P^m_{(0,2)}, \ldots, P^m_{(1,0)}, P^m_{(2,0)}, \ldots)=
P^m_{\x}.$$
\end{lem}
\noindent
\textit{Proof.} Since $\mathbf{A}^{++}_{\mathcal{K}}$ is generated by the elements $t_{(0,l)}, t_{(l,0)}$ for $l \geq 1$ there exists,
for any $\x \in \Z^{2,++}$, a polynomial $R_{\x}$ such that
$$R_{\x}(t_{(0,1)},t_{(0,2)}, \ldots, t_{(1,0)},t_{(2,0)}, \ldots)=t_{\x}.$$
By theorem~\ref{T:1}, we may take as $Q_{\x}$ the polynomial defined by
\begin{equation*}
\begin{split}
Q_{\x}(u_{(0,1)},& u_{(0,2)}, \ldots, u_{(1,0)},u_{(2,0)}, \ldots)\\
&=\frac{1}{o(d(\x))}R_{\x}(o(1)u_{(0,1)},o(2)u_{(0,2)}, \ldots, o(1)u_{(1,0)}, o(2)u_{(2,0)}, \ldots)
\end{split}
\end{equation*}
where we have set $o(l)=(q^{-l}v^l-v^{-l})/l$. \qed

\vspace{.1in}

Lemma~\ref{L:P31} implies that it is enough to show that (\ref{E:P31}) holds for $\x=(l,0)$ or $\x=(0,l)$.
This is obvious by definition for $\x=(l,0)$. We deal with the second case. It is of course enough to prove
that $\pi_m \circ \psi_m ( Sf_r(Y_1, \ldots, Y_m)S)=
\psi_{m-1} (Sf_r(Y_1, \ldots, Y_{m-1},0)S)$ for any family of symmetric polynomials $\{f_r\}$ which generates the ring $\C[Y_1, \ldots, Y_m]^{\mathfrak{S}_m}$. We may in particular take for $f_r$ the monomial symmetric function
$$m_r(Y_1, \ldots, Y_m)=\sum_{1 \leq i_1 < \cdots < i_r \leq m} Y_{i_1} \cdots Y_{i_r}.$$

We now use the following explicit computation of
$\psi_n(Sm_r(Y_1, \ldots, Y_n) S)$. See \cite[Chap. VI.5]{Mac} for a proof.

\vspace{.1in}

\begin{lem}[Macdonald] For any $n$ and any $l \geq 1$ we have
$$\psi_n(Sm_r(Y_1, \ldots, Y_n) S)=
\sum_{I \subset \{1, \ldots, n\}} A_I(x_1, \ldots, x_n) \prod_{i \in I} \partial_i,$$
where $I$ runs among all subsets of $\{1, \ldots, n\}$ of size $r$ and where
$$A_I(x_1, \ldots, x_n)=\prod_{i \in I, j \not\in I} \frac{vx_i-v^{-1}x_j}{x_i-x_j}.$$
\end{lem}

By the above lemma, we have
\begin{equation*}
\begin{split}
\pi_m \circ \psi_m(Sm_r(Y_1, \ldots, Y_m)S)&= \sum_{I \subset \{1, \ldots, m-1\}} A_I(x_1, \ldots, x_{m-1},0) \prod_{i \in I}v^{-1} \partial_i\\
&=\sum_{I \subset \{1, \ldots, m-1\}} A_I(x_1, \ldots, x_{m-1}) \prod_{i \in I} \partial_i\\
&=\psi_{m-1}(Sm_r(Y_1, \ldots , Y_{m-1})S).
\end{split}
\end{equation*}

\vspace{.1in}

We have proved that the assignment $P^m_{\x} \mapsto P^{m-1}_{\x}$
for each $\x \in \Z^{2,++}$ extends to a surjective $\mathcal{K'}$-algebra
homomorphism $\Phi_m:\SH_m^{++} \tto \SH_{m-1}^{++}$. Applying the
operator $\rho(A_1^{-k})$ and using the fact that $\rho(g) \cdot
P^n_{\x}=P^n_{g \cdot \x}$ for any $g \in SL(2,\Z)$ and $\x \in
\Z^2$, we deduce that the map $\Phi_m$ extends to an algebra
homomorphism
$$\Phi_m: \SH_{m}^{\geq -k} \tto \SH_{m-1}^{\geq -k}.$$
Here, for any $n$, we have writen $\SH_n^{\geq -k}$ for the subalgebra of
$\SH_n^+$ generated by elements
$P^n_{\x}$ where $\x=(r,d) \in \Z^{2,+}$ satisfies $d/r \geq -k$.
Letting $k$ tend to infinity, we finally obtain that the map
$\Phi_m$ extends to a surjective algebra homomorphism
$$\Phi_m: \SH_m^+ \tto \SH_{m-1}^{+}$$
such that $\Phi_m(P^m_{\x})=P^{m-1}_{\x}$ for all $\x \in \Z^{2,+}$.
This completes the proof of Proposition~\ref{P:3}. \qed

\vspace{.2in}

\paragraph{\textbf{4.2.}} Proposition~\ref{P:3} allows us to define the projective
limits $\varprojlim \mathbf{S}\ddot{\H}^+_m$ and $\varprojlim
\mathbf{S}\ddot{\H}^{++}_m$. By construction, the collection of
generators $P^m_{\x}$ for $m \geq 1$, give rise in the limit to some
elements $P_{\x}$ of these projective limits. Let $\mathbf{S}
\ddot{\H}^{+}_{\infty}$, $\mathbf{S} \ddot{\H}^{++}_{\infty}$ stand
for the subalgebras generated by $P_{\x}$ for $\x \in \Z^{2,+}$,
$\Z^{2,++}$ respectively. We may view $\varprojlim
\mathbf{S}\ddot{\H}^+_m$ and $\varprojlim
\mathbf{S}\ddot{\H}^{++}_m$ as some completions of $\mathbf{S}
\ddot{\H}^{+}_{\infty}$ and $\mathbf{S} \ddot{\H}^{++}_{\infty}$.

\vspace{.2in}

By construction the map $\Psi_m: \mathbf{A}_{\mathcal{K}} \tto
\mathbf{S}\ddot{\H}_m$ sends $\mathbf{A}^+_{\mathcal{K}}$ and
$\mathbf{A}^{++}_{\mathcal{K}}$ onto $\mathbf{S}\ddot{\H}_m^{+}$ and
$\mathbf{S}\ddot{\H}_m^{++}$ respectively. Let us call $\Psi_m^+$
and $\Psi_m^{++}$ the restrictions of $\Psi_m$ to
$\mathbf{A}^+_{\mathcal{K}}$ and $\mathbf{A}^{++}_{\mathcal{K}}$.
The collection of maps $\Psi^+_{m}$ and $\Psi^{++}_m$ give rise, in
the limit, to algebra homomorphisms
$$\Psi^+_{\infty}: \mathbf{A}^+_{\mathcal{K}} {\to}
\mathbf{S}\ddot{\H}_{\infty}^+,\quad \Psi^{++}_{\infty}:
\mathbf{A}^{++}_{\mathcal{K}} {\to}
\mathbf{S}\ddot{\H}_{\infty}^{++}.$$

\vspace{.1in}

\begin{theo}\label{T:infty} The maps $\Psi^{+}_{\infty}$ and $\Psi^{++}_{\infty}$ are isomorphisms.
\end{theo}
\noindent
\textit{Proof.} Both $\Psi^+_{\infty}$ and $\Psi^{++}_{\infty}$ are surjective by construction and we have
to show their injectivity. The subgroup
$$\mathbf{G}=\left\{\begin{pmatrix} 1 & n\\ 0 & 1 \end{pmatrix} ;n \in \Z \right\} \subset SL(2,\Z)$$
preserves $\Z^{2,+}$ and for any $\x \in \Z^{2,+}$ there exists
$g \in \mathbf{G}$ such that $g \cdot \x \in \Z^{2,++}$.
Since the map $\Psi_{\infty}^+$ is clearly compatible with the action of
$\mathbf{G}$ on $\mathbf{A}_{\mathcal{K}}^+$ and $\SH^+_{\infty}$,
we see that it is in fact enough to prove the injectivity of
$\Psi_{\infty}^{++}$.

Now fix $(r,d) \in \Z^{2,++}$. By Section~1.5.,
the dimension of the weight space $\mathbf{A}^{++}_{\mathcal{K}}[r,d]$
is equal to the number of convex paths $\mathbf{p}=(\x_1, \ldots, \x_r)$
for which $\x_i \in \Z^{2,++}$ for all $i$ and for which
$\sum_i \x_i=(r,d)$. By Proposition~\ref{P:Weyl?},
the dimension of the weight space $\SH^{++}_n[r,d]$
is equal to the dimension of the space of polynomial diagonal invariants
$$\R_n^{++}=\C[x_1, \ldots, x_n,y_1, \ldots, y_n]^{\mathfrak{S}_n}$$
of $x$-degree $r$ and $y$-degree $d$.

The latter dimension is equal to the number of orbits under $\mathfrak{S}_n$
of monomials $x_1^{g_1} \cdots x_n^{g_n}y_1^{h_1} \cdots y_n^{h_n}$
with $g_i, h_i \in \N$ satisfying $\sum_i g_i=r$ and $\sum_i h_i=d$;
equivalently, it is equal to the number of $n$-tuples of pairs
$\{(g_1,h_1), \ldots, (g_n,h_n)\}$ satisfying again
$\sum_i g_i=r$ and $\sum_i h_i=d$, or to the number of convex paths
$\p=(\x_1, \ldots, \x_r)$ in $\Z^{2,++}$ of length $r \leq n$ satisfying
$\sum_i \x_i=(r,d)$. It remains to observe that for any given $(r,d)$,
the length of convex paths $\p=(\x_1, \ldots, \x_r)$ in $\Z^{2,++}$
for which $\sum_i \x_i=(r,d)$ is bounded above by, say, $n(r,d)$. Hence
$$\text{dim\;}\mathbf{A}^{++}_{\mathcal{K}}[r,d]=\text{dim\;}\SH^{++}_n[r,d]$$
whenever $n \geq n(r,d)$, and finally
$$\text{dim\;}\mathbf{A}^{++}_{\mathcal{K}}[r,d]=\text{dim\;}\SH^{++}_\infty[r,d].$$
The injectivity of the map $\Psi_{\infty}^{++}$ follows. Theorem~\ref{T:infty} is proved. \qed

\vspace{.2in}

\paragraph{\textbf{Remark.}} 
Theorem~\ref{T:infty} allows us to transport the PBW basis
$\{\beta_{\p}; \p \in \textbf{Conv}^+\}$ and the canonical basis
$\{\mathbf{b}_{\p}; \p \in \textbf{Conv}^+\}$ 
of $\mathbf{A}^+_{\mathcal{K}}$ defined in \cite[Section~2.3]{Scano}
to bases $\{\gamma_{\p}; \p \in \textbf{Conv}^+\}$
and $\{\mathbf{c}_{\p}; \p \in \textbf{Conv}^+\}$ of $\SH^+_{\infty}$
such that
$\gamma_{\p}=\Psi_{\infty}^+(\beta_{\p})$ and
$\mathbf{c}_{\p}=\Psi^+_{\infty}(\mathbf{b}_{\p})$.
The element $\mathbf{c}_{\p}$ belongs to the completion 
$\widehat{\SH}^+_{\infty}$ of $\SH^+_{\infty}$
equal to the sum $\bigoplus_{(r,d)}\widehat{\SH}^+_{\infty}[r,d]$
over all couples $(r,d)\in\Z^{2,+}$ of the vector spaces
$$\widehat{\SH}^+_{\infty}[r,d]=\prod_{\p} \mathcal{K'}\gamma_{\p}.$$
Here $\p$ runs among all paths $\p=(\x_1, \ldots, \x_r)$ in $\textbf{Conv}^+$ 
satisfying $\sum_i \x_i=(r,d)$.






\vspace{.2in}

\paragraph{\textbf{4.3.}} 
Theorem~\ref{T:infty} should be put
in perspective with the theory of the classical Hall algebra 
$\H_{cl}$ of a discrete
valuation ring $\mathcal{O}$.
See \cite[Chap.II]{Mac}. 
Recall that $\H_{cl}$ is canonically isomorphic to the algebra
$\Lambda^+_v$,
and that this isomorphism naturally fits in a chain
\begin{equation}\label{E:bla21}
\H_{cl} \simeq \mathbf{S}\dot{\H}^+_{\infty} \simeq \Lambda^+_v
\end{equation}
where $\mathbf{S}\dot{\H}^+_{\infty}$ is the 
stable limit of the positive spherical affine Hecke algebra
of type $GL(n)$ as $n$ tends to infinity. 
Hence Theorem~\ref{T:infty} may be interpreted as an 
affine version of (\ref{E:bla21}). 
Observe that 
$\mathbf{S}\dot{\H}^+_{\infty}$ is a \textit{trivial} one-parameter 
deformation $\Lambda^+_v$ of $\Lambda^+$, 
while $\SH^+_{\infty}$ is a \textit{nontrivial} two-parameter deformation of
the ring
$$\R^+=
\C[x_1,x_2,\ldots,y_1^{\pm 1},y_2^{\pm 1},\ldots]^{\mathfrak{S}_{\infty}}.$$
The analogy may be summarized in the following table.

$$
\begin{tabular}{c|c}
Classical Hall algebra $\H_{cl}$ & Elliptic Hall algebra $\H_{el}=\mathbf{A}^+_{\mathcal{K}}$\\
&\\
\hline
&\\
$\mathcal{O}$-\textit{Mod} & $Coh(\E)$\\
&\\
$\Lambda^+=\C[x_1, x_2, \ldots]^{\mathfrak{S}_{\infty}}$ &$
\R^+=\C[x_1, x_2, \ldots, y_1^{\pm 1}, y_2^{\pm 1}, \ldots]^{\mathfrak{S}_{\infty}}$\\
&\\
$\Theta^+_{\infty}: \H_{cl} \stackrel{\sim}{\to} \mathbf{S}\dot{\H}^+_{\infty}$ &$\Psi^+_{\infty}: \H_{el} \stackrel{\sim}{\to} \SH^+_{\infty}$\\
&\\
$\Pi=\bigsqcup_n (\Z^+)^n/\mathfrak{S}_n$ & $\textbf{Conv}^+=
\bigsqcup_n (\Z^{2,+})^n/\mathfrak{S}_n$\\
&\\
$\mathbf{1}_{\mathcal{O}_{\lambda}}=v^{-2n(\lambda)}P_{\lambda}(v^2)$ &PBW-basis $\beta_{\mathbf{p}}$\\
&\\
\hline
&\\
$ \mathcal{N}_n, \ n \in \N$ & $\underline{Coh}^{r,d}(\E),
\ (r,d) \in \Z^{2,+}$\\
&\\
$\bigsqcup_{n}\; \mathcal{P}erv_{GL(n)}(\mathcal{N}_n)$ & $\bigsqcup_{r,d}\;\mathbb{U}^{r,d}$\\
&\\
$IC(\mathcal{O}_{\lambda}),\ \lambda \in \Pi$ & $\mathbb{P}_{\mathbf{p}},\
\mathbf{p} \in \textbf{Conv}^+$\\
&\\
$\Theta^+_{\infty}(tr(IC(\mathcal{O}_{\lambda}))=s_{\lambda}$& $\Psi^+_{\infty}(tr(\mathbb{P}_{\mathbf{p}}))=\mathbf{c}_{\mathbf{p}}$\\
&\\
$K_{\lambda,\mu}(v) \in \N[v]$ & $\daleth_{\mathbf{p},\mathbf{q}} \in \N[v,-t^{\pm 1}]$\\
&\\
\hline
&\\
Affine Grassmanian $\widehat{Gr}$ & ??\\
&\\
Geometric Satake isomorphism & ??\\
$\bigsqcup_n \mathcal{P}erv_{GL(n)}(\mathcal{N}_n) \simeq Rep^+GL(\infty)$&\\
&\\
\hline
\end{tabular}
$$

\vspace{.15in}

Here $P_\lambda$ is the Hall-Littlewood polynomial, 
$s_\lambda$ is the Schur polynomial,
and $K_{\lambda,\mu}$ is the Kostka polynomial.

The second part of the table is based on the geometric version of the 
elliptic Hall algebra which involves the theory of 
\textit{automorphic sheaves} defined in \cite{Laumon} and studied in 
details for an elliptic curve in \cite{Scano}. 
We refer to that paper for notations.

Lastly, in the third part of the table we mention two important features of the classical picture, for which we don't know of any analog in the setting of the elliptic Hall algebra~: functions on the nilpotent cone $\mathcal{N}_n$ may be lifted to functions on some Schubert variety of the affine 
Grassmanian~$\widehat{Gr}$ of type $GL(n)$, 
and the category of perverse sheaves 
$\bigsqcup_n\;\mathcal{P}erv_{GL(n)}(\mathcal{N}_n)$ 
is equivalent to the category $Rep^+(GL(\infty))$ 
of finite-dimensional \textit{polynomial} representations of $GL(\infty)$,
see \cite{Ginzburg}, \cite{MV}.

\vspace{.2in}

\section{Macdonald Polynomials}

\vspace{.1in}

\paragraph{\textbf{5.1.}} Macdonald discovered in the late 80's in 
\cite{MacMac} a remarkable family of symmetric polynomials $P_{\lambda}(q,t)$ 
depending on two parameters, and from which many of the classical symmetric 
functions may be obtained by specializations. 
We will use the variable $v^2$ rather than the conventional $t$ 
to comply with the notation in force in the rest of this paper.

\vspace{.1in}

The Macdonald polynomials are defined as eigenfunctions of certain difference operators acting on the spaces of symmetric functions 
$$\Lambda^m_{(q,v)}=\mathcal{K'}[x_1, \ldots, x_m]^{\mathfrak{S}_m}. $$
Recall the embedding 
$\psi_m: \SH^+_m \hookrightarrow \mathbf{D}_{m,loc}^{\mathfrak{S}_m}$,
which gives rise to an action $\rho_m$ of $\SH^+_m$ on $\Lambda^m_{(q,v)}$. 
Consider the following linear operator on
$\Lambda^m_{(q,v)}$
$$D_m=\rho_m(S (Y_1+ \cdots + Y_m)S)=
\sum_{i=1}^m \bigg(\prod_{j \neq i} 
\frac{vx_i-v^{-1}x_j}{x_i-x_j}\bigg)\partial_i.$$
By \cite{Mac}, VI, (3.10) the operator $D_m$ is upper triangular with 
respect to the basis $\{m_{\lambda}\}$ of monomial symmetric functions 
and has distinct eigenvalues.

We are interested in the stable limit as $m$ goes to infinity of the 
corresponding eigenfunctions. 
Let $\theta_m: \Lambda^m_{(q,v)}\to \Lambda^{m-1}_{(q,v)}$ 
be the specialization $x_m=0$. 
It is not true that $\theta_m \circ D_m= D_{m-1} \circ \theta_m$. 
However, the operator $E_m=v^{1-m}(D_m-[m])$ does satisfies 
$\theta_m \circ E_m= E_{m-1} \circ \theta_m$.
Recall that the space 
$$\Lambda_{(q,v)}=\mathcal{K'}[x_1, x_2, \ldots]^{\mathfrak{S}_{\infty}}$$
of symmetric functions is the projective limit of
$(\Lambda_{(q,v)}^m,\theta_m)$ in the category of graded rings.
See \cite{Mac}, Remark 1.2.1.
Hence the operators $E_m$, $m\geq 1$, 
give rise to a linear operator $E$ on the space 
$\Lambda_{(q,v)}$.
This operator is still upper triangular with respect to the basis 
$\{m_{\lambda}\}$ and has distinct eigenvalues $\{\alpha_{\lambda}\}$ given by
\begin{equation}\label{E:Maceig}
\alpha_{\lambda}=\sum_{i \geq 1} (q^{\lambda_i}-1)v^{-2(i-1)}.
\end{equation}
The Macdonald polynomial is defined to be the unique 
$\alpha_{\lambda}$-eigenvector of $E$ such that
$$P_{\lambda}(q,v^2) \in m_{\lambda} \oplus \bigoplus_{\mu<\lambda} 
\mathcal{K'}m_{\mu}.$$
For a pair of partitions $\mu \subset \lambda$, 
the skew Macdonald polynomial $P_{\lambda / \mu}(q,v^2)$ 
is determined by the coproduct formula
$$\Delta(P_{\lambda}(q,v^2))=
\sum_{\mu \subset \lambda} P_{\mu}(q,v^2) \otimes P_{\lambda/\mu}(q,v^2).$$

\vspace{.1in}

\noindent
\textbf{Examples.} i) We have $P_{(1^r)}(q,v^2)=e_r$.\\
ii) We have
$$P_{(r)}(q,v^2)=\prod_{l=0}^{r-1}\frac{1-q^{l+1}}{1-v^2q^l}\cdot\sum_{\lambda \vdash r} z_{\lambda}(q,v^2)^{-1}p_{\lambda},$$
where
$$z_{\lambda}(q,v^2)=z_{\lambda}\prod_{i=1}^{l(\lambda)}\frac{1-q^{\lambda_i}}{1-v^{2\lambda_i}},\qquad z_{(1^{m_1}2^{m_2}\cdots)}=\prod_i i^{m_i}m_i!.$$
In particular,
$$P_{(2)}(q,v^2)=\frac{(1-q)(1+v^2)}{2(1-qv^2)}p_2+\frac{(1+q)(1-v^2)}{2(1-qv^2)}p_1^2.$$

\vspace{.2in}

\noindent
\textbf{Remark.} The representations 
$\rho_m: \SH^+_m \to \text{End}(\Lambda^m_{(q,v)})$ lift, after a suitable renormalization, to a stable limit representation $\rho_{\infty}: \SH_{\infty} \to \text{End}(\Lambda_{(q,v)})$ in which $P_{(0,1)}=S (\sum_i Y_i) S$ acts as the operator $E$.
Composing with the isomorphism 
$\Psi^+_{\infty}: \mathbf{A}^+_{\mathcal{K}} \simeq \SH^+_{\infty}$ 
we obtain a representation of the Hall algebra $\mathbf{A}^+_{\mathcal{K}}$ 
on $\Lambda_{(q,v)}$ in which the element $t_{(0,1)}/\cc_1$, i.e.,
the so-called \textit{Hecke operator}, acts as Macdonald's operator $E/(q-1)$.
We will not need this representation here.

\vspace{.2in}

\paragraph{\textbf{5.2.}} There are many different characterizations of 
Macdonald polynomials. See \cite{Haiman} for instance. 
The one which fits our needs best treats the polynomials 
$P_{\lambda}(q,v^2)$ and $P_{\lambda/\mu}(q,v^2)$ at the same time. 
We first recall some standard notations from \cite{Mac}.

Let $\mu \subset \lambda$ be two partitions. 
Put $|\lambda/\mu|=|\lambda|-|\mu|$.
The skew partition $\lambda/\mu$ is said to be a \textit{vertical strip} 
if $\lambda_i-\mu_i \leq 1$ for all $i$, i.e., if the corresponding diagram 
contains at most one box per row. A skew partition $\lambda/\mu$ is a 
\textit{horizontal strip} if its conjugate $\lambda'/\mu'$ is a vertical strip.
If $\lambda/\mu$ is a horizontal strip, we put
$$\psi_{\lambda/\mu}(q,v^2)=
\prod \frac{(1-v^{2(\mu'_i-\mu'_j)}q^{j-i-1})(1-v^{2(\mu'_i-\mu'_j-1)}
q^{j-i+1})}{(1-v^{2(\mu'_i-\mu'_j)}q^{j-i})
(1-v^{2(\mu'_i-\mu'_j-1)}q^{j-i})},$$
where the sum ranges over all pairs $(i,j)$ with $i <j$ 
such that $\mu'_i=\lambda'_i$ but $\mu'_i=\lambda'_i-1$. 
In particular, we have $\psi_{\lambda/\mu}(q,v^2)=1$ if $\lambda/\mu$ 
is a horizontal strip containing no empty columns.

\vspace{.1in}

\begin{prop}\label{P:Mac51} The family 
$\{P_{\lambda/\mu}(q,v^2); \mu \subset \lambda\}$ 
is uniquely determined by the following set of properties~:
\begin{enumerate}
\item[i)] $P_{\lambda/\mu}(q,v^2)$ is homogeneous of degree $|\lambda/\mu|$,
\item[ii)] we have
$$\Delta(P_{\lambda/\mu}(q,v^2))=
\sum_{\mu \subseteq \nu \subseteq \lambda} 
P_{\nu/\mu}(q,v^2) \otimes P_{\lambda/\nu}(q,v^2),$$
\item[iii)] if $\lambda/\mu$ is not a horizontal strip then
$$P_{\lambda/\mu}(q,v^2) \in \bigoplus_{\nu < (r)} \mathcal{K'}m_{\nu},$$
where $r=|\lambda/\mu|$,
\item[iv)] if $\lambda/\mu$ is a horizontal strip then
$$P_{\lambda/\mu}(q,v^2) \in 
\psi_{\lambda/\mu}(q,v^2)m_{r}\oplus \bigoplus_{\nu < (r)} 
\mathcal{K'}m_{\nu},$$
where $r=|\lambda/\mu|$.
\end{enumerate}
\end{prop}

\noindent
\textit{Proof.} Properties i) through iv) are all known to hold for Macdonald 
polynomials: statement ii) follows from \cite[VI.7, (7.9')]{Mac}, 
while statements iii) and iv) are consequences of \cite[VI.7, (7.13')]{Mac}). 
We now prove the unicity of polynomials satisfying i) through iv). 
Let $Q_{\lambda/\mu}(q,v^2)$ be such a family. 
When $|\lambda/\mu|=1$ we have, by iv)
$$Q_{\lambda/\mu}(q^2,v)=\psi_{\lambda/\mu}(q,v^2)m_{1}=P_{\lambda/\mu}(q,v^2).$$
Let $r>1$ and assume that $Q_{\eta/\nu}(q,v^2)=P_{\eta/\nu}(q,v^2)$ 
for all $\eta/\nu$ satisfying
$|\eta/\nu| <r$. Let $\lambda/\mu$ be a skew partition with $|\lambda/\mu|=r$. By ii) and the induction hypothesis
\begin{equation*}
\begin{split}
\Delta(Q_{\lambda/\mu}&(q,v^2))=\\
&=Q_{\lambda/\mu}(q,v^2) \otimes 1 + 1 \otimes Q_{\lambda/\mu}(q,v^2) + \sum_{\mu \subset \nu \subset \lambda} Q_{\nu/\mu}(q,v^2) \otimes Q_{\lambda/\nu}(q,v^2)\\
&=Q_{\lambda/\mu}(q,v^2) \otimes 1 + 1 \otimes Q_{\lambda/\mu}(q,v^2) + \sum_{\mu \subset \nu \subset \lambda} P_{\nu/\mu}(q,v^2) \otimes P_{\lambda/\nu}(q,v^2).
\end{split}
\end{equation*}
It follows that $Q_{\lambda/\mu}(q,v^2)-P_{\lambda/\mu}(q,v^2)$ 
is contained in
$$\text{Ker}\big(\Delta-Id \otimes 1 - 1 \otimes Id\big)=
\mathcal{K'}p_{|\lambda/\mu|}.$$
But then the coefficient of $p_{|\lambda/\mu|}$ in 
$Q_{\lambda/\mu}(q,v^2)$ is uniquely determined by iii) or iv), 
and $Q_{\lambda/\mu}(q,v^2)=P_{\lambda/\mu}(q,v^2)$.\qed

\vspace{.2in}

\section{Eisenstein series}

\vspace{.1in}

\paragraph{\textbf{6.1.}} We return to the setting of Section~1, i.e.,
$\E$ is a smooth elliptic curve over $\mathbb{F}_l$, 
$\H_{\E}$ is its Hall algebra,
and $\U_{\E} \subset \H_{\E}$ is the subalgebra introduced in Section~1.4. 
For simplicity, we drop the exponent in $\U^+_{\E}$. 
Recall that $\H_\E$ and $\U_{\E}$ are
$\Z^2$-graded in the following way~:
$$\H_{\E}=\bigoplus_{(r,d)}\H_{\E}[r,d], \qquad 
\U_{\E}=\bigoplus_{(r,d)}\U_{\E}[r,d].$$
The Eisenstein series which we will need to consider are certain elements 
of a completion of the Hall algebra, which we now define in details. 
Let $\widehat{\H}_{\E}[r,d]$ stand for the space of \textit{all} functions 
$f: \Ic(\E)_{r,d} \to \C$ on the set of coherent sheaves of rank $r$ 
and degree $d$, and put 
$\widehat{\H}_{\E}=\bigoplus_{(r,d)}\widehat{\H}_{\E}[r,d]$. 
By \cite[Proposition~2.1]{BS}, the space
$\widehat{\H}_{\E}$ is still a bialgebra. 
Recall, see Section~1.3, that as a vector space we have
$$\H_{\E}[r,d]=\bigoplus_{\a_1, \ldots, \a_n} 
\H_{\E}^{(\mu(\a_1))}[\a_1] \otimes \cdots  \H_{\E}^{(\mu(\a_n))}[\a_n],$$
where the sum ranges over all tuples $(\a_1, \ldots, \a_n)$ of elements in
$\Z^{2,+}$ satisfying $\mu(\a_1) < \cdots < \mu(\a_n)$ and $\sum \a_i=(r,d)$. 
Then $\widehat{\H}_{\E}[r,d]$ is simply
$$\widehat{\H}_{\E}[r,d]=
\prod_{\a_1, \ldots, \a_n} \H_{\E}^{(\mu(\a_1))}[\a_1] \otimes 
\cdots  \H_{\E}^{(\mu(\a_n))}[\a_n].$$
In a similar fashion, we define the subalgebra 
$\widehat{\U}_{\E}$ of $\widehat{\H}_{\E}$ as
$\widehat{\U}_{\E}=\bigoplus_{(r,d)}\widehat{\U}_{\E}[r,d]$ where
$$\widehat{\U}_{\E}[r,d]=\prod_{\a_1, \ldots, \a_n} \U_{\E}^{(\mu(\a_1))}[\a_1] \otimes \cdots  \U_{\E}^{(\mu(\a_n))}[\a_n].$$
For instance, for any $(r,d)$ the element
$$\mathbf{1}_{(r,d)}=\sum_{\substack{\mathcal{F}\\ \overline{\mathcal{F}}=(r,d)}} \mathbf{1}_{\mathcal{F}}$$
is a function with infinite support belonging to $\widehat{\U}_{\E}[r,d]$ since it may be written as the infinite sum (see \cite[Equation~(4.4)]{BS})
\begin{equation}\label{E:Talpha2}
\mathbf{1}_{r,d}=\mathbf{1}^{ss}_{r,d} +  \sum_{\substack{\mu(\a_1)< \cdots < \mu(\a_n)\\\a_1 + \cdots + \a_n=(r,d)}} \nu^{\sum_{i<j}\langle \a_i,\a_j\rangle}{\mathbf{1}}_{\a_1}^{ss} \cdots {\mathbf{1}}_{\a_n}^{ss}.
\end{equation}

\vspace{.2in}

\paragraph{\textbf{6.2.}} Consider the generating series
$$\eE_0(z)=1+\sum_{d \geq 1} \mathbf{1}_{(0,d)}\nu^{-d}z^d,$$
and for $r \geq 1$
$$\eE_{r}(z)=\sum_{d \in \Z} \mathbf{1}_{(r,d)}\nu^{(r-1)d}z^d.$$
These take values in the space $\widehat{\U}_{\E}[[z,z^{-1}]]$ 
of Laurent series infinite in both directions.
We will be interested in products
$$\eE_{r_1, \ldots, r_n}(z_1, \ldots, z_n)=
\eE_{r_1}(z_1) \cdots \eE_{r_n}(z_n)\in\widehat{\U}_{\E}
[[z_1^{\pm 1}, \ldots, z_n^{\pm 1}]]$$
where $r_1, \ldots, r_n$ are nonnegative integers. 
The value of such a series at a coherent sheaf of rank $r=\sum r_i$ 
and degree $d$ is equal to the infinite sum
\begin{equation}\label{E:DefEis}
\eE_{r_1,\ldots, r_n}(\mathcal{F})=\nu^{-(r+1)d}\sum_{\substack{\mathcal{F}_1 \subset \cdots \subset \mathcal{F}_n=\mathcal{F}\\rk(\mathcal{F}_i/\mathcal{F}_{i-1})=r_i}} \nu^{2\sum_i r_i deg(\mathcal{F}_i)}z_1^{deg(\mathcal{F}_1)}\cdots z_n^{deg(\mathcal{F}_n/\mathcal{F}_{n-1})}.
\end{equation}

The following fundamental result is due to Harder.

\begin{theo}[Harder, \cite{Harder}] The series $\eE_{r_1, \ldots, r_n}(z_1, \ldots, z_n)$ converges in the region
$|z_1| \ll \cdots \ll |z_n|$ to a rational function in $\widehat{\U}_{\E}(z_1, \ldots, z_n)$ with at most simple poles along the hyperplanes $z_i/z_j \in \{1,\nu^{2}, \ldots, \nu^{2r}\}$ where $r=\sum r_i$.
\end{theo}

\vspace{.1in}

In other words, for each $\mathcal{F}$ the series (\ref{E:DefEis}) 
is the expansion in the region $|z_1| \ll \cdots \ll |z_n|$ 
of some rational function in the variables $z_1, \ldots, z_n$. 
When $r_1= \cdots =r_n=1$ the series $\eE_{1,\ldots, 1}(z_1, \ldots, z_n)$ 
is the Eisenstein series attached to the cusp form of rank one 
corresponding to the trivial character $Pic(\E) \to \C^*$ taken $n$ times.
See \cite[Section~(2.4)]{Kap} for details. 
For other values of $r_1, \ldots, r_n$ the series 
$\eE_{r_1, \ldots, r_n}(z_1, \ldots, z_n)$ is the Eisenstein series attached 
to the trivial character of the parabolic subgroup 
$GL_{r_1}(\mathbf{k}_{\E}) \times \cdots \times  
GL_{r_n}(\mathbf{k}_{\E})$ of $GL_{\sum r_i}(\mathbf{k}_{\E})$, 
where $\mathbf{k}_{\E}$ is the function field of $\E$.

\vspace{.1in}

The Eisenstein series behave well with respect to the coproduct.

\begin{prop}\label{P:Eiscop} For nonnegative integers $r_1, \ldots, r_n$ 
we have
\begin{equation*}
\begin{split}
\Delta(\eE_{r_1, \ldots, r_n}&(z_1, \ldots, z_n))=\\
&=\sum_{0 \leq s_i \leq r_i} \eE_{s_1, \ldots, s_n}(z_1, \ldots, z_n) \otimes \eE_{r_1-s_1, \ldots, r_n-s_n}(\nu^{2s_1}z_1, \ldots, \nu^{2s_n}z_n).
\end{split}
\end{equation*}
In particular, we have
$$\Delta(\eE_r(z))=\sum_{s=0}^r \eE_s(z) \otimes \eE_{r-s}(\nu^{2s}z).$$
\end{prop}
\noindent
\textit{Proof.} This is a consequence of the fact that $\widehat{\U}_{\E}$ 
is a bialgebra and that, by \cite[Equation~(4.5)]{BS} we have
$$\Delta(\mathbf{1}_{(r,d)})=\sum_{\substack{r_1+r_2=r\\d_1+d_2=d}} 
\nu^{r_1d_2-r_2d_1}\mathbf{1}_{(r_1,d_1)}\otimes \mathbf{1}_{(r_2,d_2)}.$$
\qed

\vspace{.2in}

One of the most crucial properties of Eisenstein series for us is the fact 
that they are eigenvectors for the adjoint action of the element 
$T_{(0,1)}=\sum_{x \in \E(\mathbb{F}_l)}\mathbf{1}_{\mathcal{O}_{x}}$, 
and more generally of the elements $T_{(0,d)}$ for $d \geq 1$. 
These are the so-called \textit{Hecke operators} in the theory of automorphic 
forms on function fields. Let
$$\zeta(z)=\frac{(1-\sigma z)(1-\overline{\sigma}z)}{(1-z)(1-\nu^{-2}z)}$$
be the zeta function of $\E$.

\begin{theo}\label{T:EisHeck} For any $r \geq 0$ the following holds~:
\begin{equation}\label{E:EisHeck1}
[T_{(0,1)},\eE_r(z)]=\nu \# \E(\mathbb{F}_l)\frac{\nu^{-2r}-1}{\nu^{-2}-1}z^{-1}\eE_r(z),
\end{equation}
\begin{equation}\label{E:EisHeck2}
\eE_0(z_1)\eE_r(z_2)=\prod_{i=0}^{r-1}\zeta\left(\nu^{-2i}\frac{z_1}{z_2}\right)\cdot \eE_{r}(z_2)\eE_0(z_1).
\end{equation}
In particular, we have
$\eE_0(z_1)\eE_1(z_2)=\zeta\left(\frac{z_1}{z_2}\right)\eE_1(z_2)\eE_0(z_1)$.
\end{theo}
\noindent
\textit{Proof.} Both statements are well-known (maybe in a different form) in the theory of automorphic forms. For the reader's convenience, we have included a proof in the spirit of Hall algebras in the Appendix C. \qed

\vspace{.1in}

We finish with the so-called \textit{functional equation} 
for rank one Eisenstein series.

\begin{theo}[Harder, \cite{Harder}] The rational function $\eE_{1,\ldots,1}(z_1, \ldots, z_n)$ is symmetric in variables $z_1, \ldots, z_n$.
\end{theo}

\vspace{.2in}

\paragraph{\textbf{Remark.}} Strictly speaking, the Eisenstein series most 
often considered in the theory of automorphic forms are given by expressions 
like (\ref{E:DefEis}) but in which one requires in addition each factor 
$\mathcal{F}_i/\mathcal{F}_{i-1}$ to be a vector bundle. 
In other words, if one sets
$$\mathbf{1}^{vec}_{(r,d)}=
\sum_{\substack{\mathcal{V}\;vec.\;bdle\\ \overline{\mathcal{V}}=(r,d)}} 
\mathbf{1}_{\mathcal{V}}, \qquad 
\eE_r^{vec}(z)=\sum_{d\in \Z}\mathbf{1}^{vec}_{(r,d)}\nu^{(r-1)d}z^d$$
then the corresponding product would be
$$\eE^{vec}_{r_1, \ldots, r_n}(z_1, \ldots, z_n)=\eE_{r_1}^{vec}(z_1) \cdots \eE_{r_n}^{vec}(z_n).$$
The two series, when restricted to vector bundles, 
are related by a global rational factor. 
It is the so-called \textit{$L$-factor}.
Indeed there is an obvious factorisation
$$\eE_{r}(z)=\eE_r^{vec}(z)\eE_0(\nu^{2r}z).$$
Therefore, by Theorem~\ref{T:EisHeck} we have, 
after restricting to the set of vector bundles
$$\eE_{r_1, \ldots, r_n}(z_1, \ldots, z_n)=
L_{r_1,\ldots,r_n}(z_1,\ldots,z_n)\eE^{vec}_{r_1,\ldots,r_n}(z_1,\ldots,z_n),$$
where
$$L_{r_1, \ldots, r_n}(z_1, \ldots, z_n)=\prod_{i<j}\prod_{k=0}^{r_j-1}\zeta\left(\nu^{2(r_i-k)}\frac{z_i}{z_j}\right).$$

\vspace{.2in}

\paragraph{\textbf{6.3. Example}.} 
To conclude this section, we give the example of the series 
$\eE_{1,1}(z_1,z_2)$. 
For simplicity, we will only compute the degree zero component
$$\eE_{1,1}(z_1,z_2)_0=
\sum_{d \in \Z} \left(\frac{z_1}{z_2}\right)^d 
\mathbf{1}_{(1,d)}\mathbf{1}_{(1,-d)}$$
and then only the values of $\eE_{1,1}(z_1,z_2)_0$ on vector bundles. 
So let $\mathcal{F}$ be a vector bundle of degree zero and rank two. 
Because any rank one subsheaf of $\mathcal{F}$ is a line bundle, 
and any nonzero map from a line bundle to $\mathcal{F}$ is injective, 
we have
$$\eE_{1,1}(z_1,z_2)(\mathcal{F})=
\sum_{d\in \Z}\left(\frac{z_1}{z_2}\right)^d\nu^{2d}\sum_{\mathcal{L}_{-d} 
\in Pic^{-d}(\E)} \frac{\#\text{Hom}(\mathcal{L}_{-d},\mathcal{F})-1}
{\nu^{-2}-1}.$$

If $\mathcal{F}$ is a stable bundle then
$$\text{Hom}(\mathcal{L}_{-d},\mathcal{F})=
\begin{cases} \mathbb{F}_l^{2d} & \qquad \text{if}\;d>0\\
\{0\} & \qquad \text{if} \; d \leq 0.\end{cases}$$
Hence
\begin{equation}\label{E:Eisex1}
\begin{split}
\eE_{1,1}(z_1,z_2)(\mathcal{F})&=
\frac{\# \E(\mathbb{F}_l)}{\nu^{-2}-1}
\sum_{d>0} \left(\frac{z_1}{z_2}\right)^d
\nu^{2d}(\nu^{-4d}-1)\\
&=\frac{z_1z_2(1+\nu^{-2}) 
\#\E(\mathbb{F}_l)}{(z_2-\nu^{-2}z_1)(\nu^{-2}z_2-z_1)}.
\end{split}
\end{equation}
If $\mathcal{F}=\mathcal{L}_0\oplus \mathcal{L}'_0$ 
is a direct sum of two distinct line bundles of degree zero then
$$\text{Hom}(\mathcal{L}_{-d},\mathcal{F})=
\begin{cases} \mathbb{F}_l^{2d} & \qquad \text{if}\;d>0\\
\mathbb{F}_l 
& \qquad \text{if}\;d=0\;
\text{and}\;\mathcal{L}_{-d} \in \{\mathcal{L}_0,\mathcal{L}'_0\}\\
\{0\} 
& \qquad \text{if}\;d=0\;\text{and}\;\mathcal{L}_{-d} 
\not\in \{\mathcal{L}_0,\mathcal{L}'_0\}\\
\{0\} & \qquad \text{if} \; d < 0.\end{cases}$$
Hence we get
\begin{equation}\label{E:EisEx2}
\eE_{1,1}(z_1,z_2)(\mathcal{F})=
\frac{z_1z_2 (1+\nu^{-2})\#\E(\mathbb{F}_l)}
{(z_2-\nu^{-2}z_1)(\nu^{-2}z_2-z_1)}+2.
\end{equation}
From (\ref{E:Eisex1}) and (\ref{E:EisEx2}) we deduce that the semistable 
component of $\eE_{1,1}(z_1,z_2)_0$ is equal to
\begin{equation*}
\eE_{1,1}(z_1,z_2)_{(0)}=
\frac{z_1z_2 (1+\nu^{-2})\#\E(\mathbb{F}_l)}
{(z_2-\nu^{-2}z_1)(\nu^{-2}z_2-z_1)} 
\left\{\frac{T_{(2,0)}}{[2]}+\frac{T_{(1,0)}^2}{2}\right\}
+ T_{(1,0)}^2.
\end{equation*}

Finally, to compute the unstable component of $\eE_{1,1}(z_1,z_2)_0$ 
we use the coproduct. 
Observe that since $\text{Ext}(\mathcal{L}_{-d},\mathcal{L}_d)=\{0\}$ the 
component of bidegree $(1,-d)$, $(1,d)$ of
$\Delta(\mathbf{1}_{\mathcal{L}_{-d}\oplus \mathcal{L}_{d}})$
is equal to 
$$\nu^{2d}\mathbf{1}_{\mathcal{L}_{-d}}\otimes \mathbf{1}_{\mathcal{L}_{d}},$$
and no other term may contribute to 
$\mathbf{1}_{\mathcal{L}_{-d}}\otimes \mathbf{1}_{\mathcal{L}_{d}}$. Hence
\begin{equation*}
\eE_{1,1}(z_1,z_2)(\mathcal{L}_{-d} \oplus \mathcal{L}_{d})=
\nu^{-2d}\Delta(\eE_{1,1}(z_1,z_2))(\mathcal{L}_{-d},\mathcal{L}_d)
\end{equation*}
By Proposition~\ref{P:Eiscop} and Theorem~\ref{T:EisHeck} we have
\begin{equation*}
\begin{split}
&\Delta_{1,1}(\eE_{1,1}(z_1,z_2))=\\
&=\eE_0(z_1)\eE_1(z_2)\otimes \eE_1(z_1)\eE_0(\nu^2z_2) + 
\eE_1(z_1)\eE_0(z_2)\otimes \eE_0(\nu^2z_1)\eE_1(z_2)\\
&=\zeta\left(\frac{z_1}{z_2}\right)\eE_1(z_2)\eE_0(z_1) \otimes 
\eE_1(z_1)\eE_0(\nu^2z_2) +\zeta\left(\frac{z_2}{z_1}\right)\eE_1(z_1)
\eE_0(z_2) \otimes \eE_1(z_2)\eE_0(\nu^2z_1)
\end{split}
\end{equation*}
from which we eventually obtain
$$\eE_{1,1}(z_1,z_2)(\mathcal{L}_{-d}\oplus \mathcal{L}_{d})=
\nu^{-2d}\left[\zeta\left(\frac{z_1}{z_2}\right)z_1^dz_2^{-d}+
\zeta\left(\frac{z_2}{z_1}\right)z_1^{-d}z_2^{d}\right].$$

\vspace{.2in}

\paragraph{\textbf{6.3.}} We have so far considered the Eisenstein series 
$\eE_{r_1, \ldots, r_n}(z_1, \ldots, z_n)$ for a 
\textit{fixed} elliptic curve $\E$ only. 
Recall from Section~{1.5} that there exists an algebra 
$\mathbf{A}^+_{\mathcal{A}}$ defined over the ring 
$\mathcal{A}=\C[v^{\pm 1}, t^{\pm 1}]$ 
whose specialization at $v=\nu=l^{-1/2}$ and $t=\sigma\nu$ 
for any $\E$ is isomorphic to $\U^+_{\E}$. 
Using the formulas (\ref{E:Talpha}) and (\ref{E:Talpha2}) 
we see that the generating series $\eE_r(z)$ 
and hence the Eisenstein series 
$\eE_{r_1, \ldots, r_n}(z_1, \ldots, z_n)$ 
may naturally be lifted to elements
$$\aeE_r(z) \in \widehat{\mathbf{A}}^+_{\mathcal{A}}[[z,z^{-1}]], 
\qquad \aeE_{r_1, \ldots, r_n}(z_1, \ldots, z_n)\in 
\widehat{\mathbf{A}}^+_{\mathcal{A}}[[z_1^{\pm 1}, \ldots, z_n^{\pm 1}]].$$

\vspace{.1in}

\begin{prop}\label{P:formalHarder} The series $\aeE_{r_1, \ldots, r_n}(z_1, \ldots, z_n)$ converges in the region
$|z_1| \ll \cdots \ll |z_n|$ to a rational function in 
$\widehat{\mathbf{A}}^+_{\mathcal{A}}(z_1, \ldots, z_n)$ 
with at most simple poles along the hyperplanes 
$z_i/z_j \in \{1, v^{2}, \ldots, v^{2r}\}$, where $r=\sum r_i$.
\end{prop}
\noindent
\textit{Proof.} The coefficient of $\aeE_{r_1, \ldots, r_n}(z_1, \ldots, z_n)$ 
on any basis element of
${\mathbf{A}}^+_{\mathcal{A}}$ is given by a Laurent series of the form
$$P(z_1, \ldots, z_n) \sum_{d_1, \ldots, d_n \geq 0} 
\alpha_{d_1, \ldots, d_n}\left(\frac{z_1}{z_2}\right)^{d_1}\cdots
\left(\frac{z_{n-1}}{z_n}\right)^{d_n},$$
where $P(z_1, \ldots, z_n) \in \mathcal{A}[z_1^{\pm 1}, \cdots, z_n^{\pm n}]$ 
and $\alpha_{d_1, \ldots, d_n} \in \mathcal{A}$. 
By Harder's Theorem, the evaluation at $v=\nu$ and $t=\sigma\nu$ 
for any elliptic curve $\E$ of the expression
$$\bigg(\prod_{l =1}^r \prod_{i,j} (z_i-v^{-2l}z_j)\bigg)\cdot 
P(z_1, \ldots, z_n) \sum_{d_1, \ldots, d_n}
\alpha_{d_1, \ldots, d_n}\left(\frac{z_1}{z_2}\right)^{d_1}\cdots
\left(\frac{z_{n-1}}{z_n}\right)^{d_n}$$
is a Laurent polynomial (of fixed degree). This is equivalent to the vanishing 
of certain $\mathcal{A}$-linear combinations of the 
$\alpha_{d_1, \ldots, d_n}$'s. Of course, if such a linear combination 
vanishes when evaluated at all (i.e., infinitely many) elliptic curves $\E$ 
then it must already vanish in $\mathcal{A}$. We are done.\qed

\vspace{.2in}

\paragraph{\textbf{6.4.}} Motivated by the analogy between the Hecke operator 
$T_{(0,1)}$ and Macdonald's operator (see Section~{5.1} and Remark~{5.1}) 
we introduce, for every partition $\lambda=(\lambda_1, \ldots, \lambda_n)$ 
the following specialization of Eisenstein series
$$\eE_{\lambda}(z)=\aeE_{\lambda_1, \ldots, \lambda_n}
(z,q^{-1}z,\ldots, q^{1-n}z)$$
where $q=vt$. By Proposition~\ref{P:formalHarder}, 
the line $(z,q^{-1}z,\ldots, q^{1-n}z)$ is not contained in the pole locus of 
$\aeE_{\lambda_1, \ldots, \lambda_n}(z_1, \ldots, z_n)$ and hence 
$\eE_{\lambda}(z)$ belongs to 
$\widehat{\mathbf{A}}_{\mathcal{A}}\otimes_{\mathcal{A}}\mathcal{K}(z)$. 
More generally, for any pairs of partitions
$\mu \subset \lambda$ we put
$$\eE_{\lambda/\mu}(z)=\aeE_{\lambda_1-\mu_1, \ldots, \lambda_n-\mu_n}
(v^{2\mu_1}z,v^{2\mu_2}q^{-1}z,\ldots, v^{2\mu_n}q^{1-n}z).$$
Observe that by Theorem~\ref{T:EisHeck} the series $\eE_{\lambda}(z)$ 
are eigenvectors for the adjoint action of the Hecke operator $T_{(0,1)}$, 
whose eigenvalues $\beta_{\lambda}$ are (up to a global factor) 
equal to that of the Macdonald polynomials, namely we have
$$\beta_{\lambda}=z^{-1}\mathbf{c}_1(v,t)\sum_i \frac{v^{-2\lambda_i}-1}
{v^{-2}-1}q^{i-1}=z^{-1}
\frac{\mathbf{c}_1(v,t)}{q-1}\alpha_{\lambda'}$$
where $\alpha_{\lambda'}$ is given by formula (\ref{E:Maceig}) 
and $\lambda'$ is the conjugate partition to $\lambda$.

\vspace{.2in}

It would seem natural to define more generally the specialization
$$\eE_{\underline{l}}(z)=\aeE_{l_1, \ldots, l_n}(z,q^{-1}z,\ldots, q^{1-n}z)$$
for \textit{any} sequence of nonnegative integers $l_1, \ldots, l_n$. 
However we have the following vanishing result.

\begin{lem}\label{L:Van} If $\underline{l}=(l_1, \ldots, l_n)$ 
is not dominant, i.e., if $l_k > l_{k-1}$ for some $k$, 
then $\eE_{\underline{l}}(z)=0$.
\end{lem}
\noindent
\textit{Proof.} One may check that the $L$-factor 
$L_{l_1, \ldots, l_n}(z_1, \ldots, z_n)$ vanishes on the line
$(z,q^{-1}z, \ldots, q^{1-n}z)$ whenever $\underline{l}$ is not dominant. 
Hence Lemma~\ref{L:Van} would follow from the fact that the 
unnormalized Eisenstein series $\eE^{vec}_{l_1, \ldots, l_n}(z_1, \ldots, z_n)$
is regular on that line. Rather than appealing to this fact, 
we provide a direct proof. To unburden the notation, we drop the subscript 
$\mathcal{A}$ throughout. By Proposition~\ref{P:Eiscop} we have
$$\Delta_{(1,\ldots, 1)}(\eE_r(z))=
\eE_1(z) \otimes \eE_1(v^2z)\otimes \cdots \otimes \eE_1(v^{2(r-1)}z),$$
and more generally given integers $\epsilon_k \in \{0,1\}$ with 
$\epsilon_{i_1}=\cdots =\epsilon_{i_r}=1$ while $\epsilon_k=0$ if 
$k \not\in\{i_1, \ldots, i_r\}$ we have
\begin{equation}\label{E:vanproof1}
\Delta_{(\epsilon_1, \ldots, \epsilon_n)}(\eE_r(z))=
\eE_{\epsilon_1}(z) \otimes \cdots \otimes \eE_{\epsilon_k}(v^{2s_k}z) \otimes 
\cdots \otimes \eE_{\epsilon_n}(v^{2s_n}z)
\end{equation}
where $s_k=\#\{l; i_l <k\}$. Now let $\underline{l}=(l_1, \ldots, l_n) \in \N^n$
and set $l=\sum l_i$. We may compute
$\Delta_{(1, \ldots, 1)}(\eE_{\underline{l}}(z))$ using (\ref{E:vanproof1}). 
It is equal to a sum, indexed by the set of maps
$\phi: \{1, \ldots, l\} \to \{1, \ldots, n\}$ of terms
$$a_{\phi}=\Delta_{(\epsilon_1^1, \ldots, \epsilon_n^1)}(\eE_{l_1}(z)) \cdots \Delta_{(\epsilon_1^k, \ldots, \epsilon_n^k)}(\eE_{l_k}(q^{1-k}z)) \cdots \Delta_{(\epsilon_1^n, \ldots, \epsilon_n^n)}(\eE_{l_n}(q^{1-n}z)),$$
where $\epsilon_i^k \in \{0,1\}$ is defined by
$$\epsilon_{i}^k=\begin{cases} 0\quad & if\;\phi(i) \neq k,\\ 
1 \quad & if\;\phi(i)=k.\end{cases}$$
In other terms, the map $\phi$ describes the way the coproducts 
(\ref{E:vanproof1}) of the $\eE_{l_k}(q^{1-k}z)$ have been distributed 
among the $l$ components of the tensor product.
We claim that if $\underline{l}$ is not dominant then each term $a_{\phi}$ 
vanishes. Indeed, suppose that $l_k >l_{k-1}$ for some $k$. 
Then $a_{\phi}$ is divisible by a term of the form
\begin{equation}\label{E:vanproof2}
\begin{split}
&\Delta_{(\epsilon_1, \ldots, \epsilon_n)}(\eE_{l_{k-1}}(q^{2-k}z)) \cdot 
\Delta_{(\epsilon'_1, \ldots, \epsilon'_n)}(\eE_{l_k}(q^{1-k}z))=\\
&\quad=\eE_{\epsilon_1}(q^{2-k}z)\eE_{\epsilon'_1}(q^{1-k}z) \otimes \cdots \otimes \eE_{\epsilon_n}(v^{2s_n}q^{2-k}z)\eE_{\epsilon'_n}(v^{2s'_n}q^{1-k}z).
\end{split}
\end{equation}
Of course, in the above if $\epsilon_i=1$ then $\epsilon'_i=0$ and vice versa.
As $s_1=s'_1=0$ while $s_n \in \{l_{k-1},l_{k-1}-1\}$ and 
$s'_n \in \{l_k,l_k-1\}$ it is easy to see that there exists an index $j$ 
for which $\epsilon_j=0, \epsilon'_j=1$ and $s_j=s'_j$. 
But then the $j$th component of (\ref{E:vanproof2}) is equal to
$$\eE_0(v^{2s_j}q^{2-k}z) \eE_1(v^{2s_j}q^{1-k}z)=
\zeta(q)\eE_1(v^{2s_j}q^{1-k}z)\eE_0(v^{2s_j}q^{2-k}z)=0$$
since $\zeta(q)=0$. Hence $a_{\phi}=0$ as wanted and 
$\Delta_{(1,\ldots,1)}(\eE_{\underline{l}}(z))=0$.
It remains to show that the map
$$\Delta_{(1,\ldots,1)}:~\widehat{\U}[r,d] \to \prod_{d_1+\cdots+d_r=d} 
\widehat{\U}[1,d_1] \otimes \cdots \otimes \widehat{\U}[1,d_r]$$
is injective. This in turn follows from the fact that $\widehat{\U}$ 
is equipped with a nondegenerate Hopf pairing, 
and that it is generated by elements of degree zero and one,
see \cite[Cor.~5.1]{BS}.\qed

\vspace{.1in}

\begin{prop}\label{P:Kiral} For any partition $\lambda$ we have
\begin{equation}\label{E:kiral1}
\Delta(\eE_{\lambda}(z))=\sum_{\mu \subseteq \lambda}\eE_{\mu}(z) \otimes \eE_{\lambda/\mu}(z).
\end{equation}
More generally, for any skew partition $\lambda/\mu$ we have
\begin{equation}\label{E:kiral2}
\Delta(\eE_{\lambda/\mu}(z))=\sum_{\mu \subseteq \nu \subseteq \lambda}\eE_{\nu/\mu}(z) \otimes \eE_{\lambda/\nu}(z).
\end{equation}
\end{prop}
\noindent
\textit{Proof.} We prove the first statement. 
By Proposition~\ref{P:Eiscop} it holds
\begin{equation}\label{E:kiral3}
\begin{split}
\Delta&(\eE_{\lambda}(z, \ldots, q^{1-n}z))=\\
&=\sum_{\substack{s_1, \ldots, s_n\\0 \leq s_i \leq \lambda_i}} 
\aeE_{s_1, \ldots, s_n}(z, \ldots, q^{1-n}z) \otimes 
\aeE_{\lambda_1-s_1, \ldots, \lambda_n-s_n}
(v^{2s_1}z, \ldots, v^{2s_n}q^{1-n}z).
\end{split}
\end{equation}
By Lemma~\ref{L:Van} we have
$\aeE_{s_1, \ldots, s_n}(z,q^{-1}z, \ldots, q^{1-n}z)=0$ 
if $(s_1, \ldots, s_n)$ is not a partition. 
Therefore the r.h.s. of (\ref{E:kiral3}) reduces to (\ref{E:kiral1}). 
The proof of the second statement of the Proposition is similar.\qed

\vspace{.2in}

\section{Geometric construction of Macdonald polynomials}

\vspace{.15in}

In this section, we make explicit the link between Macdonald polynomials 
$P_{\lambda}(q,v^2)$ and the  Eisenstein series $\eE_{\lambda}(z)$.

\vspace{.15in}

\paragraph{\textbf{7.1.}} For any skew partition $\lambda/\mu$ we denote by
$\eE_{\lambda/\mu}^{(0)}$ the restriction of $\eE_{\lambda/\mu}(z)$ 
to the set of \textit{semistable} vector bundles of degree zero. 
Notice that by homogeneity this is independent of $z$. 
This is therefore an element of the subalgebra
$\mathbf{A}^{+,(0)}_\mathcal{A}$ of the universal Hall algebra 
$\mathbf{A}^+_{\mathcal{A}}$ generated by elements $\widetilde{t}_{(r,0)}$ 
for $r \geq 0$. See Section~{1.5} for details. By Proposition~\ref{P:Macd} 
this last subalgebra is canonically identified with the algebra of 
symmetric functions $\Lambda_{(q,v)}$, 
where we have set $q=tv$.
Explicitly, the isomorphism is given by $\tilde{t}_{(r,0)}=p_r/r$.

\vspace{.1in}

For instance, from Example~{6.3} we see that
\begin{equation*}
\begin{split}
\eE_{1,1}^{(0)}&=\frac{q^{-1}(1+v^{-2})(1-v^{-2}q)(1-q^{-1})}
{(q^{-1}-v^{-2})(q^{-1}v^{-2}-1)}
\left(\frac{p_2}{2}+\frac{p_1^2}{2}\right) + p_1^2\\
&=\frac{(1+v^{-2})(q-1)}{v^{-2}-q}\frac{p_2}{2}+
\frac{(v^{-2}-1)(q+1)}{v^{-2}-q}\frac{p_1^2}{2}.
\end{split}
\end{equation*}

We let $\omega$ stand for the standard involution on symmetric functions, 
defined by $\omega(p_r)=(-1)^{r-1}p_r$.

\vspace{.2in}

\paragraph{\textbf{7.2.}} We are now ready to state the second main Theorem 
of this paper.

\vspace{.1in}

\begin{theo}\label{T:MacEis} For any partition $\lambda$ we have
$$\eE_{\lambda}^{(0)}=\omega P_{\lambda'}(q,v^2),$$
and for any skew partition $\lambda/\mu$ we have
$$\eE_{\lambda/\mu}^{(0)}=\omega P_{\lambda'/\mu'}(q,v^2).$$
\end{theo}

\vspace{.1in}

The rest of this section is devoted to the proof of this theorem. 
We will use the characterization of the polynomials $P_{\lambda/\mu}(q,v^2)$ 
given in Proposition~\ref{P:Mac51}. It is clear from the definitions that 
$\omega (\eE_{\lambda'/\mu'}^{(0)})$ is of degree $|\lambda/\mu|$. 
Property ii) of Proposition~\ref{P:Mac51} is shown for 
$\omega (\eE_{\lambda'/\mu'}^{(0)})$ in Proposition~\ref{P:Kiral}. 
Thus it only remains to check that the coefficient of $m_{r}$ 
in $\omega (\eE_{\lambda'/\mu'}^{(0)})$
for $r=|\lambda/\mu|$ is given by Proposition~\ref{P:Mac51}, iii) and iv).

\vspace{.1in}

To this aim we introduce the following family of elements in 
$\mathbf{A}^{+,(0)}$~:
$$g_r=
\sum_{\lambda \vdash r} z_{\lambda}^{-1}\prod_i\frac{v^{-\lambda_i}-v^{\lambda_i}}{\mathbf{c}_{\lambda_i}(v,t)}t_{\lambda_i}
=\sum_{\lambda \vdash r}z_{\lambda}^{-1}\prod_i\frac{v^{-2\lambda_i}-1}{(1-q^{-r})(1-(qv^{-2})^r)} p_{\lambda}.$$
Alternatively, these may be defined by the formula
$$1+\sum_{r >0}g_rs^r=exp\bigg((v^{-1}-v)\sum_{r \geq 1} 
\frac{t_{(r,0)}}{v^r+v^{-r}-t^r-t^{-r}}s^r\bigg).$$
Recall that $\mathbf{A}^+$ is equipped with a nondegenerate Hopf scalar product.
By \cite[Lemma~5.2]{BS} it satisfies
\begin{equation}\label{E:Hopscal}
\langle t_{(r,0)}, t_{(s,0)}\rangle=
\delta_{r,s}\frac{[r]^2\#\E(\mathbb{F}_{l^r})}{r(v^{-2r}-1)}=
\delta_{r,s}\frac{\mathbf{c}_r(v,t)}{v^{-1}-v},
\end{equation}
so that, after identification with $\Lambda_{(q,v)}$ it reads
\begin{equation}\label{E:Hopscal2}
\langle p_{r}, p_{s}\rangle=
\delta_{r,s}r\frac{(1-q^{-r})(1-(qv^{-2})^r)}{v^{-2r}-1}.
\end{equation}
Using \cite[Chap VI.2]{Mac} we deduce that $g_r$ is dual to $m_r$ 
with respect to the basis $\{m_{\lambda}\}$, i.e., that
$$\langle g_r,m_r\rangle=1,$$
$$\langle g_r, m_{\lambda} \rangle=0\;\quad{if}\;|\lambda|=r\;{and}\;\lambda < (r).$$
Therefore the proof of Theorem~\ref{T:MacEis} will be complete 
once we have shown that
\begin{equation}\label{E:ask1}
\langle g_r, \omega(\eE_{\lambda'/\mu'}^{(0)})\rangle=0
\end{equation}
if $\lambda/\mu$ is not a horizontal strip while
\begin{equation}\label{E:ask2}
\langle g_r, \omega(\eE_{\lambda'/\mu'}^{(0)})\rangle=
\psi_{\lambda/\mu}(q,v^2)
\end{equation}
if $\lambda'/\mu'$ is a horizontal strip. 
As we will see, these equations essentially amount to
certain relations between the factors $\psi_{\lambda/\mu}(q,v^2)$ and the 
$L$-factors appearing in the Eisenstein series.

Observe that as $g_r$ is itself semistable of degree zero, i.e., we have
$g_r \in \mathbf{A}^{+,(0)}$, and the subalgebras
$\mathbf{A}^{+,(\mu)}$ are all mutually orthogonal, 
we may as well replace $\eE_{\lambda'/\mu'}^{(0)}$ by $\eE_{\lambda'/\mu'}(z)$ 
in equations~(\ref{E:ask1}) and (\ref{E:ask2}). 
Note also that $\omega$ is an orthogonal involution for $\langle\;,\;\rangle$.

\vspace{.2in}

\paragraph{\textbf{7.3.}} The basic idea is to find a factorization of $g_r$ 
and to use the Hopf property of the scalar product  $\langle\;,\;\rangle$ 
to reduce (\ref{E:ask1}) and (\ref{E:ask2}) to a lower rank. Of course, 
since $g_r$ is dual to $m_r$ and $m_r$ is primitive, this is not directly 
feasable. However, it becomes possible as soon as we step out of the 
subalgebra $\mathbf{A}^{+,(0)}$. More precisely, put
$g_r^{(1)}=[t_{(0,1)},g_r]$ and $\omega g_r^{(1)}=[t_{(0,1)},\omega g_r].$

\begin{lem}\label{L:factgr} For any $r \geq 1$ we have
$$\omega g_{r+1}^{(1)}=\frac{v}{\mathbf{c}_1(v,t)}[t_{(1,0)},\omega g_r^{(1)}]+
(v^{-1}-v)\omega g_r t_{(1,1)}.$$
\end{lem}
\noindent
\textit{Proof.} An essentially direct computation, based on the relation 
$[t_{(s,1)},t_{(u,0)}]=\mathbf{c}_{u}t_{(s+u,1)}$ for any $s,u$, yields
\begin{equation}\label{E:factgr}
\omega g_r^{(1)}=(-1)^r\sum_{s=1}^r v^s(1-v^{-2})(-1)^s\omega g_{r-s}t_{(s,1)}.
\end{equation}
The recursion formula in the lemma is an easy 
consequence of (\ref{E:factgr}).\qed

\vspace{.1in}

\begin{lem}\label{L:eiggr1} For any skew partition $\lambda/\mu$ we have
$$\langle \omega g_r^{(1)}, \eE_{\lambda'/\mu'}(z)\rangle=\frac{\mathbf{c}_1(v,t)}{1-v^2}\big( \sum_i (v^{2\mu'_i}-v^{2\lambda'_i})q^{1-i}z
\big) \langle \omega g_r, \eE_{\lambda'/\mu'}(z)\rangle.$$
\end{lem}
\noindent
\textit{Proof.} Because $\langle\;,\;\rangle$ is a Hopf pairing we have
\begin{equation*}
\begin{split}
\langle \omega g_r^{(1)}, \eE_{\lambda'/\mu'}(z)\rangle&=\langle t_{(0,1)}\cdot  \omega g_r -\omega g_r \cdot t_{(0,1)}, \eE_{\lambda'/\mu'}(z)\rangle\\
&=\langle t_{(0,1)} \otimes \omega g_r, \Delta_{0,r}(\eE_{\lambda'/\mu'}(z))\rangle -
\langle \omega g_r\otimes t_{(0,1)}, \Delta_{r,0}(\eE_{\lambda'/\mu'}(z))\rangle .
\end{split}
\end{equation*}
Using Proposition~\ref{P:Eiscop} the coproducts are computed to be
\begin{equation*}
\begin{split}
\Delta_{0,r}(\eE_{\lambda'/\mu'}(z))&=\eE_{0}(v^{2\mu'_1}z) \cdots \eE_{0}(v^{2\mu'_n}q^{1-n}z)\otimes \eE_{\lambda'/\mu'}(z)\\
&=(1 + \sum_i v^{2\mu'_i-1}q^{1-i}z t_{(0,1)}+\cdots ) \otimes\eE_{\lambda'/\mu'}(z)
\end{split}
\end{equation*}
and
\begin{equation*}
\begin{split}
\Delta_{r,0}(\eE_{\lambda'/\mu'}(z))&=\eE_{\lambda'/\mu'}(z) \otimes \eE_{0}(v^{2\lambda'_1}z) \cdots \eE_{0}(v^{2\lambda'_n}q^{1-n}z) \\
&=\eE_{\lambda'/\mu'}(z) \otimes(1 + \sum_i v^{2\lambda'_i-1}q^{1-i}z t_{(0,1)}+\cdots ). \end{split}
\end{equation*}
The lemma follows.
\qed

\vspace{.2in}

\paragraph{\textbf{7.4.}} We now proceed with the proof of (\ref{E:ask1}) and 
(\ref{E:ask2}). We argue by induction on $|\lambda/\mu|$. 
Assume first that $|\lambda/\mu|=1$. This means that $\lambda'_i=\mu'_i$ 
for all $i$ except for one value, say $j$, for which $\lambda'_j=\mu'_j+1$. 
Then on the one hand
$$P_{\lambda/\mu}(q,v^2)=
\psi_{\lambda/\mu}m_{1}=
\prod_{i<j} \frac{(1-v^{2(\mu'_i-\mu'_j)}q^{j-i-1})
(1-v^{2(\mu'_i-\mu'_j-1)}q^{j-i+1})}
{(1-v^{2(\mu'_i-\mu'_j)}q^{j-i})(1-v^{2(\mu'_i-\mu'_j-1)}q^{j-i})}
\cdot m_{1}$$
while on the other hand
\begin{equation*}
\eE_{\lambda'/\mu'}(z)=
\eE_{0}(v^{2\mu'_1}z)\cdots \eE_0(v^{2\mu'_{j-1}}q^{2-j}z)
\eE_1(v^{2\mu'_j}q^{1-j}z)\cdots \eE_0(v^{2\mu'_n}q^{1-n}z)
\end{equation*}
so that by Theorem~\ref{T:EisHeck} we get
\begin{equation*}
\begin{split}
\eE_{\lambda'/\mu'}^{(0)}&=
\prod_{i<j}\zeta (v^{2(\mu'_i-\mu'_j)}q^{j-i})\cdot
\eE_1(v^{2\mu'_j}q^{1-j}z)^{(0)}\\
&=\prod_{i<j}\zeta (v^{2(\mu'_i-\mu'_j)}q^{j-i})\cdot m_{1}.
\end{split}
\end{equation*}
It remains to notice that
$$\frac{(1-v^{2(\mu'_i-\mu'_j)}q^{j-i-1})(1-v^{2(\mu'_i-\mu'_j-1)}q^{j-i+1})}
{(1-v^{2(\mu'_i-\mu'_j)}q^{j-i})(1-v^{2(\mu'_i-\mu'_j-1)}q^{j-i})}=
\zeta\left(v^{2(\mu'_i-\mu'_j)}q^{j-i}\right)$$
so that the $\psi$-factor and the $L$-factors indeed coincide
$$\psi_{\lambda/\mu}(q,v^2)=
\prod_{i<j}\zeta\left(v^{2(\mu'_i-\mu'_j)}q^{j-i}\right).$$

\vspace{.1in}

Next let us assume that equations (\ref{E:ask1}) 
and (\ref{E:ask2}) hold true for all skew partitions
$\nu/\eta$ for which $|\nu/\eta|<r$, 
and let $\lambda/\mu$ be a skew partition of size $r$. 
Combining Lemmas~\ref{L:eiggr1} and \ref{L:factgr} we have
\begin{equation}\label{E:proof1}
\begin{split}
&\frac{\mathbf{c}_1(v,t)}{1-v^2}\big(\sum_i (v^{2\mu'}-v^{2\lambda'})
q^{1-i})z\big)\langle \omega g_r,\eE_{\lambda'/\mu'}(z)\rangle=\\
\;&=\frac{v}{\mathbf{c}_1(v,t)}\left\{ \langle t_{(1,0)} \otimes \omega 
g_{r-1}^{(1)}, \Delta_{1,r-1}(\eE_{\lambda'/\mu'}(z))\rangle -
\langle  \omega g_{r-1}^{(1)}\otimes t_{(1,0)}, 
\Delta_{r-1,1}(\eE_{\lambda'/\mu'}(z))\rangle\right\}\\
&\qquad\qquad\qquad\qquad\qquad\qquad\qquad\qquad\;\; + (v^{-1}-v)\langle 
\omega g_{r-1} \otimes t_{(1,1)}, \Delta_{r-1,1}(\eE_{\lambda'/\mu'}(z))\rangle.
\end{split}
\end{equation}
Let us first consider the case of a vertical strip $\lambda'/\mu'$ 
(hence $\lambda/\mu$ is a horizontal strip). Hence
$$\eE_{\lambda'/\mu'}(z)=\eE_{\epsilon_1}(v^{2\mu'_1}z) \cdots \eE_{\epsilon_n}(v^{2\mu'_n}q^{1-n}z)$$ for some $\epsilon_i \in \{0,1\}$. Let $I$ (resp. $J$) be the set of $k \in \{1, \ldots, n\}$ for which $\epsilon_k=0$ (resp. $\epsilon_k=1$). Then
$$\eE_{\lambda'/\mu'}(z)=\prod_{\substack{i<j\\i \in I, j \in J}} \zeta (v^{2(\mu'_i-\mu'_j)}q^{j-i}) \cdot \vec{\prod}_{j \in J} \eE_1(v^{2\mu'_j}q^{1-j}z) \cdot \vec{\prod}_{i \in I}\eE_0(v^{2\mu'_i}q^{1-i}z).$$
As above, the $\psi$-factor and the $L$-factor coincide
$$\psi_{\lambda/\mu}(q,v^2)=\prod_{\substack{i<j\\i \in I, j \in J}} \zeta (v^{2(\mu'_i-\mu'_j)}q^{j-i}) $$
and therefore (\ref{E:ask2}) reduces to the simple relation
$$\langle \omega g_r,\eE_1(v^{2\mu'_{j_1}}q^{1-j_1}z) \cdots \eE_1(v^{2\mu'_{j_r}}q^{1-j_r}z)\rangle=1.$$
We claim that in fact $\langle \omega g_r,\eE_1(\a_1) \cdots \eE_1(\a_r)\rangle=1$ for \textit{any} $\a_1, \ldots, \a_r$. Developping (\ref{E:proof1}) one finds that this is equivalent to the following strange identity, which is proved in Appendix~D.

\vspace{.1in}

\begin{lem}\label{L:strangeid} For any $r \geq 1$, the following identity holds over the field of rational functions $\mathcal{K'}(\a_1, \ldots, \a_r)$~:
\begin{equation}\label{E:strangeid}
\begin{split}
\sum_{i=1}^r \a_i=\frac{1}{v^{-2}-1}&\sum_{j=1}^r \left[ \bigg( \prod_{l \neq j} \zeta \left( \frac{\a_l}{\a_j}\right) -\prod_{l \neq j} \zeta \left( \frac{\a_j}{\a_l}\right)\bigg) \cdot \bigg(\sum_{l \neq j} \a_l \bigg)\right]\\
&+\sum_{j=1}^r \bigg(\prod_{l \neq j} \zeta \left( \frac{\a_j}{\a_l}\right) \bigg) \a_j.
\end{split}
\end{equation}
\end{lem}

Next, let us assume that $\lambda'/\mu'$ has exactly one part of length two and $r-2$ parts of length one.
Arguing as above and cancelling the $L$-factor, we see that (\ref{E:ask1}) is equivalent to
$$\langle \omega g_r,\eE_1(v^{2\mu'_{j_1}}q^{1-j_1}z) \cdots \eE_2(v^{2\mu'_{k}}q^{1-k}z) \cdots \eE_1(v^{2\mu'_{j_{r-1}}}q^{1-j_{r-1}}z)\rangle=0.$$
Again, we claim that in fact $\langle \omega g_r,\eE_1(\a_1) \cdots \eE_2(\a_k) \cdots \eE_1(\a_{r-1})\rangle=0$ for any $\a_1, \ldots, \a_{r-1}$. This may be checked directly using (\ref{E:proof1}).

In all the remaining cases $\lambda'/\mu'$ has at least two parts of length at 
least two, i.e., if $\lambda'_i> \mu'_i+1$ for more than one value of $i$. 
But then no sub skew partition of size $r-1$ of $\lambda'/\mu'$ 
may be vertical. 
By the induction hypothesis this implies that all terms on the r.h.s.~
vanish and thus $\langle \omega g_r, \eE_{\lambda'/\mu'}(z)\rangle=0$ as wanted.
Theorem~\ref{T:MacEis} is proved.\qed

\vspace{.2in}

\paragraph{\textbf{Remarks.}} i) A factorization similar to (\ref{E:factgr}) 
involving rank one difference operators in the context of Pieri 
rules for skew Macdonald polynomials appears in \cite{Haiman2}. 
We thank Mark Haiman for this remark.\\

\noindent
ii) In addition to Macdonald's operator $\Delta_1$, 
one defines an operator $\nabla$ acting on symmetric polynomials in 
$\Lambda_{(q,v)}$ (see \cite{Haiman2}), which has distinct eigenvalues and 
whose eigenvectors are the Macdonald polynomials. 
Namely $\nabla$ is defined by
$$\nabla (P_{\lambda}(q,v^2))=
v^{-2n(\lambda)}q^{n(\lambda')}P_{\lambda}(q,v^2).$$
Our conventions, taken from \cite{Mac}, differ slightly from \cite{Haiman2}.
In our picture, this operator $\nabla$ is simply given 
by the action of the element $A_2\in SL(2,\Z)$ 
by automorphism on the Hall algebra, i.e.,
it is the tensor product with a line bundle of degree one.
Thus we have
$$\rho(A_2)( \mathbf{E}_{\lambda}(z))=
v^{-2n(\lambda')}q^{n(\lambda)}\mathbf{E}_{\lambda}(z).$$

\noindent
iii) Laumon defined and studied in \cite{Laumon} a ``geometric lift'' 
of Eisenstein series to certain perverse sheaves 
(or more precisely, constructible complexes) on the stacks 
$\underline{Coh}^{r,d}(\E)$ called \textit{Eisenstein sheaves}. 
The Eisenstein series themselves are recovered from the Eisenstein sheaves 
via the faisceaux-function correspondence. In the special case of an 
elliptic curve \textit{simple} Eisenstein sheaves are determined 
in \cite{Scano}.  The construction of the (non simple) 
Eisenstein sheaves relevant to Macdonald polynomials may be easily 
translated from Theorem~\ref{T:1}. Let us denote by 
$(\qlb)_{r,d}$ the trivial rank one constructible sheaf on 
$\underline{Coh}^{r,d}(\E)$, and let us consider the formal 
series whose coefficients are semisimple constructible complexes
$$\mathbb{E}_{r}(q^{-l}z)=
\bigoplus_{d \in \Z}(\qlb)_{r,d}[(r-1)d](ld)_{\E}z^d.$$
Here $[n]$ is the standard shift of complexes and 
$(m)_{\E}$ denotes the Tate twist \textit{by the Frobenius eigenvalue} 
$\sigma$ in $H^1(\E, \qlb)$. Note that there is a choice of one 
Frobenius eigenvalue $\sigma$ involved here, 
but of course choosing the other eigenvalue $\overline{\sigma}$ 
would give a similar result.
Using the induction functor of Laumon \cite{Laumon} we may form the product
$$\mathbb{E}_{\lambda}(z)=\mathbb{E}_{\lambda_1}(z) \star 
\mathbb{E}_{\lambda_2}(q^{-1}z) \cdots \star \mathbb{E}_{\lambda_l}(q^{1-l}z).$$
It is still a series with coefficients in semisimple 
constructible complexes. These will usually be of infinite rank. 
Restricting to the open substack parametrizing semistable 
sheaves of zero slope we finally obtain a semisimple constructible 
complex $\mathbb{E}_{\lambda}^{(0)}$. Using \cite{Scano}, 
Proposition~6.1, one can show that the Frobenius eigenvalues of 
$\mathbb{E}_{\lambda}(z)$ and $\mathbb{E}_{\lambda}^{(0)}$ 
all belong to $v^{\Z}q^{\Z}$. Hence the Frobenius trace 
$Tr(\mathbb{E}_{\lambda}^{(0)})$ is a Laurent series in $v$ and $q$. 
Recall that we have fixed an isomorphism
$\mathbb{C}\simeq\qlb$.
By Harder's theorem 
the series $Tr(\mathbb{E}_{\lambda}^{(0)})$  
converges (in a suitable domain) to 
$\mathbf{E}_{\lambda}^{(0)}$ and hence by Theorem~\ref{T:1} we have
$Tr(\mathbb{E}_{\lambda}^{(0)})=\omega P_{\lambda'}(q,v^2).$
\\

\noindent
iv) Pick a $\mathbb{F}_l$-rational
closed point $x\in\E(\mathbb{F}_l)$. Let $i:D_x\to\E$ be the
embedding of the
formal neighborhood of $x$ in $\E$. Given an \'etale coordinate at $x$ we get
an isomorphism $D_x\simeq\text{Spec}(\mathbb{F}_l((\varpi)))$, 
where $\varpi$ is a formal variable. 
Thus the set of isomorphism classes of torsion sheaves on $D_x$ 
is equal to the set of conjugacy classes of nilpotent matrices.
Invariant functions on the nilpotent cone $\Nc_d$ , $d\geq 1$,
are canonically identified with elements of the ring $\Lambda^+$
of symmetric functions.
The restriction of coherent sheaves on $\E$ to $D_x$ 
yields a map
$\Ic(\E)_{0,d}\to\coprod_{d'\le d}\Nc_{d'}$.
It factors to an algebra isomorphism
$\Lambda^+\simeq\U^{+,(\infty)}_\E$.
Fourier-Mukai transform 
yields an algebra isomorphism
$FM:\U^{+,(0)}_\E\to\U^{+,(\infty)}_\E$.
The composed map 
$\U^{+,(0)}_\E\to\Lambda^+$
coincides with the isomorphism 
in Proposition~\ref{P:Macd}.

The involution $\omega$ in Theorem~\ref{T:MacEis}
can be removed as follows.
We'll give another isomorphism $\U^{+,(\infty)}_\E\simeq\Lambda^+$ 
which takes the Laurent series
$FM(\mathbf{E}_{\lambda}^{(0)})$ to $P_{\lambda'}(q,v^2)$.
Let $\E^{(d)}$ be the $d$-th symmetric power of $E$,
and $\underline{\widetilde{Coh}}^{0,d}(\E)$ be the stacks of flags
$$\Mc_d\to\Mc_{d-1}\to\cdots\Mc_1,$$ 
where each $\Mc_i$ is a coherent sheaf on $\E$ of length $i$.
Consider the Cartesian square
$$\xymatrix{
\E^{d} \ar[r]^-{\iota_{d}} 
\ar[d]^-{r_d}
& \underline{\widetilde{Coh}}^{0,d}(\E)
\ar[d]^-{\pi_d}\\
\E^{(d)} \ar[r]^-{\iota_{(d)}} &
\underline{{Coh}}^{0,d}(\E)
}$$ 
in which $\pi_d$ is the Springer map,
$r_d$ is the ramified finite cover 
$(x_1,x_2,\dots x_d)\mapsto x_1+x_2+\cdots x_d$,
and $\iota_{(d)}$ takes a divisor $D$ to the sheaf $\Oc_D$.
According to Laumon, the complex $F=R(\pi_d)_*(\qlb)$
is the intermediate extension of its restriction
$F|_{U_d}$
to the dense open subset $U_d=\iota_{(d)}(\E^{(d)})$.
We have $\iota_{(d)}^*F=(r_d)_*(\qlb)$ by base change.
Thus the symmetric group 
$\mathfrak{S}_d$ acts on $F|_{U_d}$.
For each irreducible character $\phi$ of
$\mathfrak{S}_d$ let $F_\phi$ be the intermediate extension
of the constructible sheaf $Hom_{\mathfrak{S}_d}(\phi,F|_{U_d})$.
Each $F_\phi$ is a simple constructible complex 
on $\underline{{Coh}}^{0,d}(\E)$.
The representation ring of $\mathfrak{S}_d$ is canonically identified with a 
subring of $\Lambda^+$.
We claim that there is an unique isomorphism 
$\U^{+,(\infty)}_\E\simeq\Lambda^+$ 
taking $Tr(F_\phi)$ to the symmetric function associated to $\phi$.
This is the map we want.

\vspace{.4in}

\centerline{APPENDICES}

\vspace{.2in}

\appendix

\section{Proof of Proposition~\ref{P:1}.}

\vspace{.1in}

\paragraph{\textbf{A.1.}} We start the proof of Proposition~\ref{P:1} with a sequence of lemmas.

\vspace{.1in}

\begin{lem}\label{L:41}  We have
\begin{equation}\label{E:L41}
S \bigg[ X_1, \sum_i Y_i \bigg] S=(1-q)SX_1Y_1S.
\end{equation}
\end{lem}
\noindent
\textit{Proof.} Using Lemma~\ref{L:1}, equations (\ref{E:L11}) and 
(\ref{E:L12}), we get
\begin{equation*}
\begin{split}
S \sum_i& Y_i X_1S\\
&= S X_1 (\sum_{i \geq 2} Y_i)S + qv^{2(n-1)}SX_1Y_1S + (v^{-2}-1)SY_1X_1S \\
&\qquad+ (v^{-4}-v^{-2})SY_1X_1S + \cdots + (v^{-2(n-1)}-v^{-2(n-2)})SY_1X_1S\\
&=S X_1 (\sum_{i \geq 2} Y_i)S + qv^{2(n-1)}SX_1Y_1S+qv^{2(n-1)}(v^{-2(n-1)}-1)SX_1Y_1S\\
&=SX_1(\sum_i Y_iX_1)S + (q-1)SX_1Y_1S.
\end{split}
\end{equation*}
\qed

\vspace{.1in}

\begin{lem}\label{L:42} For any indices $2 \leq j_2 < j_3 < \cdots < j_l \leq n$ we have
\begin{equation}\label{E:L42}
SY_1Y_{j_2} \cdots Y_{j_l}X_1=qv^{2(n-l)}SX_1Y_1Y_{j_2} \cdots Y_{j_l}.
\end{equation}
\end{lem}
\noindent
\textit{Proof.} We have $SY_1X_1=qv^{2(n-1)}SX_1Y_1$ by Lemma~\ref{L:1}, 
equation~(\ref{E:L11}).
By (\ref{E:L12}) we have 
\begin{equation*}
Y_{j_2}X_1=X_1Y_{j_2}+(v^{-1}-v)T_{j_2-1}^{-1} \cdots T_1^{-1} \cdots T_{j_2-1}^{-1}Y_1X_1.
\end{equation*}
Multiplying by $Y_{j_3} \cdots Y_{j_l}$ and using the fact 
that $[T_k,Y_h]=0$ if $h >k-1$ we deduce
\begin{equation*}
\begin{split}
Y_{j_2}Y_{j_3} \cdots Y_{j_l}X_1&=Y_{j_3} \cdots Y_{j_l}X_1Y_{j_2}+\\
& \qquad +(v^{-1}-v)T_{j_2-1}^{-1} \cdots T_1^{-1} \cdots T_{j_2-1}^{-1}Y_1Y_{j_3} \cdots Y_{j_l}X_1.
\end{split}
\end{equation*}
Multiplying now by $Y_1$ and using the relation 
$Y_1T_1^{-1} \cdots Y_{j_2-1}^{-1}=T_1 \cdots T_{j_2-1}Y_{j_2}$ yields
\begin{equation*}
\begin{split}
Y_1Y_{j_2} \cdots Y_{j_l}X_1&=Y_1Y_{j_3} \cdots Y_{j_l}X_1Y_{j_2}+\\
& \qquad+(v^{-1}-v)T_{j_2-1}^{-1} \cdots T_2^{-1}T_1 \cdots T_{j_2-1}Y_1Y_{j_2} \cdots Y_{j_l}X_1,
\end{split}
\end{equation*}
from which it follows in turn that
$$SY_1Y_{j_2} \cdots Y_{j_l}X_1=SY_1Y_{j_3} \cdots Y_{j_l}X_1Y_{j_2}+(1-v^2)SY_1Y_{j_2} \cdots Y_{j_l}X_1$$
and thus that
$$v^2 SY_1 Y_{j_2} \cdots Y_{j_l}X_1=SY_1 Y_{j_3} \cdots Y_{j_l}X_1Y_{j_2}.$$
By the same argument,
$$v^2 S Y_1 Y_{j_3} \cdots Y_{j_l}X_1Y_{j_2}=SY_1Y_{j_4} \cdots Y_{j_l}X_1Y_{j_2}Y_{j_3},$$
and continuing in this manner we finally arrive at
$$SY_1 Y_{j_2} \cdots Y_{j_l}X_1=v^{-2(l-1)}SY_1X_1Y_{j_2} \cdots Y_{j_l}=qv^{2(n-l)}SX_1Y_1Y_{j_2} \cdots Y_{j_l}$$
as expected.\qed

\vspace{.1in}

\begin{lem}\label{L:43} For any indices $1< j_1 < j_2 \cdots < j_l \leq n$ we have
\begin{equation}\label{E:L43}
\begin{split}
SY_{j_1} &\cdots Y_{j_l}X_1=\\
&=SX_1Y_{j_1} \cdots Y_{j_l}+\sum_{u=1}^l q(1-v^2)v^{2(n-l-j_u+u)}SX_1Y_1 Y_{j_1} \cdots \widehat{Y_{j_u}} \cdots Y_{j_l}.
\end{split}
\end{equation}
Here the symbol $\widehat{x}$ means that we omit the term $x$ in the product.
\end{lem}
\noindent
\textit{Proof.} First of all, we have, again by Lemma~\ref{L:1}, (\ref{E:L12})
$$Y_jX_1=X_1Y_j+(v^{-1}-v)\beta_jY_jX_1,$$
for all $j >1$ where we have set
$$\beta_j=T_{j-1}^{-1} \cdots T_1{-1} \cdots T_{j-1}^{-1}.$$
Define elements
$$A(j_1, \ldots, j_l)=SY_{j_1} \cdots Y_{j_l}X_1, \qquad B(j_2, \ldots, j_l)=SY_1Y_{j_2} \cdots Y_{j_l}X_1.$$
We have, by the same arguments as in the previous lemma,
\begin{equation*}
\begin{split}
A(j_1, \ldots, j_l)&=A(j_2, \ldots, j_l)Y_{j_1} + (v^{-1}-v) S \beta_{j_1}Y_{1} Y_{j_2} \cdots Y_{j_l}X_1\\
&=A(j_2, \ldots, j_l)Y_{j_1} + (v^{-1}-v)v^{-2(j_1-1)+1}B(j_2, \ldots, j_l).
\end{split}
\end{equation*}
By Lemma~\ref{L:42}, $B(j_2, \ldots, j_l)=qv^{2(n-l)}SX_1Y_1Y_{j_2} \cdots Y_{j_l}$, therefore
\begin{equation*}
\begin{split}
A(j_1, \ldots, j_l)&=A(j_2, \ldots, j_l) Y_{j_1} + q(1-v^2)v^{2(n-l-j_1+1)}SX_1Y_1 Y_{j_2} \cdots Y_{j_l}\\
&=\left( A(j_3, \ldots, j_l)Y_{j_2}+q(1-v^2)v^{2(n-l-j_2+2)}SX_1Y_1Y_{j_3} \cdots Y_{j_l}\right)Y_{j_1}\\
&\qquad + q(1-v^2)v^{2(n-l-j_1+1)}SX_1Y_1 Y_{j_2} \cdots Y_{j_l}\\
&=\cdots\\
&=SX_1 Y_{j_1} \cdots Y_{j_l} +
\sum_{u=1}^l q(1-v^2)v^{2(n-l-j_u+u)}SX_1Y_{j_1} \cdots \widehat{Y_{j_u}} \cdots Y_{j_l}
\end{split}
\end{equation*}
which is what we wanted to prove.\qed

\vspace{.1in}

\begin{lem}\label{L:44} The following holds~:
\begin{equation}\label{E:L44}
S \bigg[ X_1, \sum_{j_1 < \cdots < j_l} Y_{j_1} \cdots Y_{j_l}\bigg] =(1-q)\sum_{1<j_2< \cdots <j_l} SX_1Y_1Y_{j_2} \cdots Y_{j_l}.
\end{equation}
\end{lem}
\noindent
\textit{Proof.} Using the previous two lemmas, we compute
\begin{equation*}
\begin{split}
\sum_{j_1< \cdots <j_l}& SY_{j_1} \cdots Y_{j_l}X_1\\
&=\sum_{1<j_2 < \cdots <j_l} SY_1 Y_{j_2} \cdots Y_{j_l}X_1 + \sum_{1<j_1 < \cdots < j_l} SY_{j_1} \cdots Y_{j_l}X_1\\
&=\sum_{1 < j_2 < \cdots < j_l} qv^{2(n-l)} SX_1Y_1 Y_{j_2} \cdots Y_{j_l} + \sum_{1<j_1< \cdots <j_l} SX_1 Y_{j_1} \cdots Y_{j_l} \\
& \qquad + \sum_{1<j_1 < \cdots <j_l} \left\{\sum_{u=1}^l q(1-v^2)v^{2(n-l-j_u+u)}SX_1Y_1 Y_{j_1} \cdots \widehat{Y_{j_u}} \cdots Y_{j_l}\right\}\\
&=\sum_{1 < j_2 < \cdots < j_l} qv^{2(n-l)} SX_1Y_1 Y_{j_2} \cdots Y_{j_l} + \sum_{1<j_1< \cdots <j_l} SX_1 Y_{j_1} \cdots Y_{j_l} \\
&\qquad+\sum_{1<k_1< \cdots <k_{l-1}}qv^{2(n-l)}(1-v^2)\sigma_{k_1, \ldots, k_{l-1}}SX_1Y_1Y_{k_1} \cdots Y_{k_{l-1}}
\end{split}
\end{equation*}
where
\begin{equation*}
\begin{split}
\sigma_{k_1, \ldots, k_{l-1}}&=\bigg\{(v^{-2} + \cdots +v^{2(1-(k_1-1))})+(v^{2(2-(k_1+1))}+ \cdots + v^{2(2-(k_2-1))}) \\
&\hspace{.4in} +\cdots + (v^{2(l-(k_{l-1}+1))} + \cdots + v^{2(l-n)})\bigg\}\\
&=\frac{v^{2(l-n)}-1}{1-v^2}
\end{split}
\end{equation*}
Hence,
\begin{equation*}
\begin{split}
\sum_{j_1< \cdots <j_l}& SY_{j_1} \cdots Y_{j_l}X_1\\
&=\sum_{1 < j_2 < \cdots < j_l} qv^{2(n-l)} SX_1Y_1 Y_{j_2} \cdots Y_{j_l} + \sum_{1<j_1< \cdots <j_l} SX_1 Y_{j_1} \cdots Y_{j_l} \\
&\qquad+\sum_{1<k_1 < \cdots < k_{l-1}} qv^{2(n-l)}(v^{2(l-n)}-1)SX_1Y_1Y_{k_1} \cdots Y_{k_{l-1}}\\
&=q\sum_{1 < j_2 < \cdots < j_l} qv^{2(n-l)} SX_1Y_1 Y_{j_2} \cdots Y_{j_l} + \sum_{1<j_1< \cdots <j_l} SX_1 Y_{j_1} \cdots Y_{j_l},
\end{split}
\end{equation*}
from which the lemma follows. \qed

\vspace{.2in}

We are finally ready to give the proof of Proposition~\ref{P:1}. We will argue by induction, with Lemma~\ref{L:41} being the case $l=1$. So we fix $l \in \N$ and assume that Proposition~\ref{P:1} has been proved for all $l'<l$. It is necessary to distinguish two cases~:

\vspace{.1in}

\noindent
\textit{Case 1.} Let us assume that $l \leq n$. We will use the formula
\begin{equation*}
\begin{split}
\sum_i Y_i^l&=\bigg(\sum_i Y_i^{l-1} \cdot \sum_i Y_i\bigg)-\bigg(\sum_i Y_i^{l-2}\cdot\sum_{i<j} Y_iY_j\bigg)\\
&\quad+ \cdots + (-1)^l\bigg(\sum_i Y_i\cdot\sum_{j_1 < \cdots < j_{l-1}} Y_{j_1} \cdots Y_{j_{l-1}}\bigg)\\
&\quad+(-1)^{l+1}l\sum_{j_1 < \cdots < j_l}Y_{j_1} \cdots Y_{j_l}.
\end{split}
\end{equation*}
According to the above, we have therefore
\begin{equation}\label{E:Case1}
\begin{split}
S\bigg[ X_1,&\sum_i Y_i^l \bigg]S\\
=& S \bigg[X_1, \sum_i Y_i \bigg]S\cdot\sum_i Y_i^{l-1} + \sum_i Y_i\cdot S \bigg[ X_1,\sum_i Y_i^{l-1}\bigg]S\\
&- S \bigg[X_1, \sum_{j_1<j_2} Y_{j_1}Y_{j_2} \bigg]S\cdot \sum_i Y_i^{l-2}  -\sum_{j_1<j_2} Y_{j_1}Y_{j_2}\cdot S \bigg[ X_1,\sum_i Y_i^{l-2}\bigg]S\\
&+\cdots\\
&+(-1)^l S\bigg[ X_1, \sum_{j_1< \cdots < j_{l-1}}Y_{j_1} \cdots Y_{j_{l-1}}\bigg]S\cdot  \sum_i Y_i\\
&+(-1)^l  \sum_{j_1< \cdots < j_{l-1}}Y_{j_1} \cdots Y_{j_{l-1}}\cdot S\bigg[ X_1, \sum_i Y_i \bigg]S\\
&+(-1)^{l+1} l S \bigg[ X_1, \sum_{j_1< \cdots < j_{l}}Y_{j_1} \cdots Y_{j_{l}}\bigg]S.
\end{split}
\end{equation}

By the induction hypothesis and Lemma~\ref{L:44},
$$S\bigg[X_1, \sum_i Y_i^t\bigg]S=(1-q^t)SY_1^tX_1S,$$
$$S \bigg[ X_1, \sum_{j_1 < \cdots < j_s} Y_{j_1} \cdots Y_{j_s}\bigg] S=(1-q)\sum_{1<j_2< \cdots <j_s} SX_1Y_1Y_{j_2} \cdots Y_{j_s}S $$
for all $t<l$ and for all $s$. Substituting in (\ref{E:Case1}), we deduce that
\begin{equation}\label{E:Case12}
\begin{split}
S\bigg[X_1&,\sum_i Y_i^l\bigg]S\\
=&(1-q)S X_1Y_1\sum_i Y_i^{l-1}S + (1-q^{l-1})S\sum_i Y_i X_1 Y_1^{l-1}S\\
&-(1-q)SX_1Y_1\sum_{1<j_2} Y_{j_2}\sum_i Y_i^{l-2}S -(1-q^{l-2})S\sum_{j_1<j_2} Y_{j_1}Y_{j_2}X_1Y_1S\\
&+\cdots\\
&+ (-1)^{l}(1-q)SX_1Y_1\sum_{1<j_2< \cdots < j_{l-1}}Y_{j_2} \cdots Y_{j_{l-1}}\sum_i Y_iS\\
&+(-1)^{l}(1-q)S\sum_{j_1< \cdots < j_{l-1}}Y_{j_1} \cdots Y_{j_{l-1}}X_1Y_1S\\
&+(-1)^{l+1}l(1-q)SX_1 Y_1 \sum_{1<j_2 < \cdots < j_l} Y_{j_2} \cdots Y_{j_l}S.
\end{split}
\end{equation}
This is where we use Lemma~\ref{L:44} again, in the form
\begin{equation*}
\begin{split}
S\sum_{j_1< \cdots < j_t}Y_{j_1} \cdots Y_{j_t}X_1= SX_1&\sum_{j_1< \cdots < j_t}Y_{j_1} \cdots Y_{j_t}
\\
&+(q-1)SX_1Y_1 \sum_{1<j_2< \cdots <j_t} Y_{j_2} \cdots Y_{j_t}
\end{split}
\end{equation*}
to obtain the following expression for the bracket $S\bigg[X_1,\sum_i Y_i^l \bigg]S$~:
\begin{equation*}
\begin{split}
S&\bigg[X_1,\sum_i Y_i^l \bigg]S\\
&=(1-q)SX_1 \bigg\{ Y_1 \sum_i Y_i^{l-1}-Y_1 \sum_{1<j_2}Y_{j_2} \sum_iY_i^{l-2} +Y_1\sum_{1<j_2<j_3}Y_{j_2}Y_{j_3}\sum_i Y_i^{l-3} -\cdots\\
&\quad \quad +\cdots +(-1)^lY_1\sum_{1<j_2 < \cdots <j_{l-1}}Y_{j_2} \cdots Y_{j_{l-1}} \sum_i Y_i +
(-1)^{l+1}lY_1\sum_{1<j_2 < \cdots <j_{l}}Y_{j_2} \cdots Y_{j_{l}}\bigg\}\\
&\quad \quad+(1-q^{l-1})SX_1Y_1^{l-1}\sum_i Y_iS + (1-q^{l-1})(q-1)SX_1Y_1^lS\\
& \quad \quad -(1-q^{l-2})SX_1Y_1^{l-2}\sum_{j_1<j_2}Y_{j_1}Y_{j_2}-(1-q^{l-2})(q-1)SX_1Y_1^{l-1}\sum_{1<j_2}Y_{j_2}S\\
&\quad \quad + \cdots\\
&\quad \quad +(-1)^l(1-q)SX_1Y_1\sum_{j_1< \cdots <j_{l-1}}Y_{j_1} \cdots Y_{j_{l-1}}S\\
&\quad \quad +(-1)^l(1-q)(q-1)SX_1Y_1^2\sum_{j_2< \cdots <j_{l-1}}Y_{j_2} \cdots Y_{j_{l-1}}S.
\end{split}
\end{equation*}
Collecting terms we get
\begin{equation*}
\begin{split}
S&\bigg[X_1,\sum_i Y_i^l \bigg]S\\
&=SX_1Y_1^lS\bigg\{ (1-q)+(1-q^{l-1})+(1-q^{l-1})(q-1)\bigg\}\\
&\quad \quad +SX_1Y_1^{l-1}\sum_{1<j_2}Y_{j_2}S\bigg\{(q-1)+(1-q^{l-1})-(1-q^{l-2})(q-1)-(1-q^{l-2})\bigg\}\\
&\quad \quad +SX_1Y_1^{l-2}\sum_{1<j_2<j_3}Y_{j_2}Y_{j_3}S\bigg\{(1-q)+(q^{l-2}-1)-(q^{l-3}-1)(q-1)-(q^{l-3}-1)\bigg\}\\
&\quad \quad + \cdots\\
&\quad \quad +SX_1Y_1^2\sum_{1<j_2 < \cdots <j_{l-1}}Y_{j_2} \cdots Y_{j_{l-1}}S \bigg\{(-1)^l\bigg((1-q)+(q^2-1)-(q-1)^2-(q-1)\bigg)\bigg\}\\
&\quad \quad +SX_1Y_1\sum_{1<j_2 < \cdots <j_{l}}Y_{j_2} \cdots Y_{j_{l}}S \bigg\{(1-q)\bigg((-1)^l(l-1)+(-1)^{l+1}l+(-1)^l\bigg)\bigg\}\\
&\quad \quad +\sum_{j>1}\sum_{t=1}^{l-2}SX_1Y_1Y_j^{l-t}\sum_{\underset{j_u \neq j}{1<j_2 < \cdots <j_{t}}}Y_{j_2} \cdots Y_{j_t} \bigg\{(1-q)\bigg((-1)^{t+1}+(-1)^{t+2}\bigg)\bigg\}\\
&=(1-q^l)SX_1Y_1^lS.
\end{split}
\end{equation*}
This concludes the proof of Proposition~\ref{P:1} in the first case.

\vspace{.2in}

\noindent
\textit{Case 2.} Let us deal with the situation when $l > n$. 
The method is very similar to the one used in the proof of Case 1 above. 
This time we use the following identity~:
\begin{equation*}
\begin{split}
\sum_i Y_i^l=\sum_i Y_i^{l-1} \cdot \sum_i Y_i-\sum_i Y_i^{l-2} \cdot 
\sum_{j_1<j_2} Y_{j_1}Y_{j_2}
+ \cdots  (-1)^{n-1}\sum_i Y_i^{l-n}Y_1 \cdots Y_n.
\end{split}
\end{equation*}
Based on this decomposition, we may write
\begin{equation}\label{E:Case2}
\begin{split}
S\bigg[ X_1,&\sum_i Y_i^l \bigg]S=\\
=& S \bigg[X_1, \sum_i Y_i \bigg]S\cdot\sum_i Y_i^{l-1} + \sum_i Y_i\cdot S \bigg[ X_1,\sum_i Y_i^{l-1}\bigg]S\\
&- S \bigg[X_1, \sum_{j_1<j_2} Y_{j_1}Y_{j_2} \bigg]S\cdot \sum_i Y_i^{l-2}  -\sum_{j_1<j_2} Y_{j_1}Y_{j_2}\cdot S \bigg[ X_1,\sum_i Y_i^{l-2}\bigg]S\\
&+\cdots\\
&+(-1)^{n-1} S\bigg[ X_1, Y_{1} \cdots Y_{n}\bigg]S\cdot  \sum_i Y_i^{l-n}\\
&+(-1)^{n-1} Y_1 \cdots Y_n  S\bigg[ X_1, \sum_i Y_i^{l-n} \bigg]S.
\end{split}
\end{equation}

Using the induction hypothesis and Lemma~\ref{L:44}, this simplifies to

\begin{equation*}
\begin{split}
S\bigg[& X_1,\sum_i Y_i^l \bigg]S=\\
=& (1-q)SX_1 Y_1 \sum_i Y_i^{l-1}S + (1-q^{l-1})S\sum_i Y_i X_1 Y_1^{l-1}S\\
&\quad \quad -(1-q)SX_1\sum_{1<j_2}Y_{j_2}Y_1\sum_i Y_i^{l-2}S -(1-q^{l-2})S\sum_{j_1<j_2}Y_{j_1}Y_{j_2}X_1Y_1^{l-2}S\\
&\quad \quad + \cdots\\
&\quad \quad (-1)^{n-1}\bigg( (1-q)SX_1Y_1 \cdots Y_n
\sum_i Y_i^{l-n}S+(1-q^{l-n})SY_1 \cdots Y_n X_1Y_1^{l-n}S\bigg)\\
=&\bigg\{(1-q)S X_1Y_1\sum_i Y_i^{l-1}S+(1-q^{l-1})SX_1\sum_i Y_iY_1^{l-1}S\\
& \qquad \qquad \qquad \qquad \qquad\qquad\qquad \qquad+(1-q^{l-1})(q-1)SX_1Y_1^lS\bigg\}\\
&\quad \quad -\bigg\{(1-q)SX_1\sum_{1<j_2}Y_{j_2}Y_1\sum_iY_i^{l-2}S+(1-q^{l-2})SX_1 \sum_{j_1<j_2}Y_{j_1}Y_{j_2}Y_1^{l-2}S\\
& \qquad \qquad \qquad \qquad\qquad \qquad \qquad\qquad+(1-q^{l-2})(q-1)SX_1\sum_{1<j_2}Y_{j_2}Y_1^{l-1}S\bigg\}\\
&\quad \quad + \cdots\\
&\quad \quad +(-1)^{n-1}\bigg\{(1-q)SX_1Y_1 \cdots Y_n \sum_{i}Y_i^{l-n}S + (1-q^{l-n})SX_1Y_1 \cdots Y_n Y_1^{l-n}S\\
&\qquad \qquad \qquad \qquad\qquad \qquad \qquad \qquad+(1-q^{l-n})(q-1)SX_1 Y_2 \cdots Y_n Y_1^{l-n+1}S\bigg\}
\end{split}
\end{equation*}

Gathering terms, we obtain
\begin{equation*}
\begin{split}
S\bigg[& X_1,\sum_i Y_i^l \bigg]S=\\
&=SX_1Y_1^lS\bigg\{(1-q)+(1-q^{l-1})+(q-1)(1-q^{l-1})\bigg\}\\
&\quad \quad + SX_1\sum_{1<j_2}Y_{j_2}Y_1^{l-1}S \bigg\{-(1-q)+(1-q^{l-1})-(1-q^{l-2})-(q-1)(1-q^{l-2})\bigg\}\\
&\quad \quad + \cdots \\
&\quad \quad +SX_1Y_2 \cdots Y_n Y_1^{l-n+1}S\bigg\{(-1)^{n-1}\bigg((1-q)-(1-q^{l-n+1})+(1-q^{l-n})\\
& \hspace{3.7in} +(q-1)(1-q^{l-n})\bigg)\bigg\}\\
&\quad \quad + \sum_{i>1}\sum_{t=1}^n SX_1 \sum_{\underset{j_u \neq i}{1<j_2 < \cdots < j_{t-1}}}Y_{j_2} \cdots Y_{j_{t-1}}Y_i^{l-t}S \bigg\{(-1)^{t+1}(q-1)+(-1)^t(q-1)\bigg\}\\
&= (1-q^l)SX_1Y_1^lS
\end{split}
\end{equation*}
as desired.
This concludes the proof of Case 2 as well as that of Proposition~\ref{P:1}. 
\qed

\vspace{.2in}

\section{Proof of Proposition~\ref{P:2}.}

\vspace{.1in}

We will start by giving a closed expression for the commutator $S \bigg[ \sum_i Y_i, X_1^lY_1^{-1}\bigg]S$.

\begin{lem}\label{L:51}  For any $l \geq 1$ we have
\begin{equation}\label{E:L51}
\begin{split}
S \bigg[ \sum_i Y_i, X_1^lY_1^{-1}\bigg]S=& (v^{2(1-n)}-1)S\bigg\{qX_1^{l-1} X_n+q^2X_1^{l-2}X_n^2 + \cdots + q^lX_n^l\bigg\}S\\
& +q^lSX_n^lS-SX_1^lS
\end{split}
\end{equation}
\end{lem}
\noindent
\textit{Proof.} First of all, by (\ref{E:L13})
\begin{equation*}
\begin{split}
Y_1X_1Y_1^{-1}&=qT_1 \cdots T_{n-2}T_{n-1}^2 T_{n-2} \cdots T_1X_1\\
&=qT_1 \cdots T_{n-2} T_{n-1}X_nT_{n-1}^{-1} \cdots T_1^{-1},
\end{split}
\end{equation*}
from which it follows that
\begin{equation}\label{E:517}
Y_1X_1^lY_1^{-1}=q^lT_1 \cdots T_{n-2} T_{n-1}X^l_nT_{n-1}^{-1} \cdots T_1^{-1}
\end{equation}
and hence that $SY_1X_1^lY_1^{-1}S=q^lSX_n^lS$. Now we compute
\begin{equation*}
\begin{split}
S\bigg[\sum_i Y_i, X_1^lY_1^{-1}\bigg]S&=SY_1X_1^lY_1^{-1}S-SX_1^lS+\sum_{m =2}^n S \bigg[Y_m,X_1^lY_1^{-1}\bigg]S\\
&=q^lSX_n^lS-SX_1^lS+\sum_{m =2}^n S \bigg[Y_m,X_1^lY_1^{-1}\bigg]S.
\end{split}
\end{equation*}
The lemma will thus be proved once we have shown that
\begin{equation}\label{E:L512}
S \bigg[Y_m,X_1^lY_1^{-1}\bigg]S=(1-v^2)v^{2(1-m)}S\bigg\{qX_1^{l-1}X_n+q^2X_1^{l-2}X_n^2 + \cdots + q^lX_n^l\bigg\}S.
\end{equation}
For this, we need some preparatory result

\begin{slem} For any $m$, the following identity holds~:
\begin{equation}\label{E:L513}
\begin{split}
T_{m-1}^{-1} \cdots T_2^{-1}T_1^{-1}T_2^{-1} \cdots T_{m-1}^{-1}T_1& \cdots T_{n-2}T_{n-1}^2 T_{n-2} \cdots T_1\\
&=T_1^{-1} \cdots T_{m-2}^{-1} T_m T_{m+1} \cdots T_{n-2}T_{n-1}^2 T_{n-2} \cdot T_1.
\end{split}
\end{equation}
\end{slem}
\noindent
\textit{Proof.} We argue by induction. The relation can easily be checked directly for $m=2$. Fix $m$ and assume that (\ref{E:L513}) holds for $m-1$. We have
\begin{equation*}
\begin{split}
T_{m-1}^{-1} &\cdots T_2^{-1}T_1^{-1}T_2^{-1} \cdots T_{m-1}^{-1}T_1 \cdots T_{n-2}T_{n-1}^2 T_{n-2} \cdots T_1\\
&=T_{m-1}^{-1} \cdots T_1^{-1}T_2^{-1}T_1^{-1} \cdots T_{m-1}^{-1}T_1 \cdots T_{n-2}T_{n-1}^2 T_{n-2} \cdots T_1\\
&=T_1^{-1}T_{m-1}^{-1} \cdots T_3^{-1}T_2^{-1}T_3^{-1} \cdots T_{m-1}^{-1}T_1^{-1}T_1 \cdots T_{n-2}T_{n-1}^2 T_{n-2} \cdots T_1\\
&=T_1^{-1}\bigg(T_{m-1}^{-1} \cdots T_3^{-1}T_2^{-1}T_3^{-1} \cdots T_{m-1}^{-1}T_2 \cdots T_{n-2}T_{n-1}^2 T_{n-2} \cdots T_2 \bigg) T_1.
\end{split}
\end{equation*}
Using the induction hypothesis applied to the set of indices $2, 3, \ldots, n$ instead of $1,2, \ldots, n$, we may simplify the expression in parenthesis to get
\begin{equation*}
\begin{split}
T_{m-1}^{-1} &\cdots T_2^{-1}T_1^{-1}T_2^{-1} \cdots T_{m-1}^{-1}T_1 \cdots T_{n-2}T_{n-1}^2 T_{n-2} \cdots T_1\\
&=T_1^{-1} \bigg(T_2^{-1} \cdots T_{m-2}^{-1}T_m T_{m+1} \cdots T_{n-2} T_{n-1}^2 T_{n-2} \cdots T_{2} \bigg) T_1
\end{split}
\end{equation*}
which proves (\ref{E:L513}) for the integer $m$. This finishes the induction step and the proof of the sublemma. \qed

We may now prove Lemma~\ref{L:51}. We argue once again by induction. Fix $m$ and set $u_l=Y_mX_1^lY_1^{-1}$.  We compute $u_1$ directly, using (\ref{E:L12}) and (\ref{E:L513})~:

\begin{equation*}
\begin{split}
u_1&=X_1Y_mY_1^{-1}+(v^{-1}-v)T_{m-1}^{-1} \cdots T_2^{-1}T_1^{-1}T_2^{-1} \cdots
T_{m-1}^{-1}Y_1X_1Y_1^{-1}\\
&=X_1Y_mY_1^{-1}+q(v^{-1}-v)T_{m-1}^{-1} \cdots T_1^{-1} \cdots T_{m-1}^{-1}T_1 \cdots T_{n-1}^2 \cdots T_1 X_1Y_1Y_1^{-1}\\
&=X_1Y_mY_1^{-1}+q(v^{-1}-v)T_1^{-1} \cdots T_{m-2}^{-1} T_m \cdots T_{n-1}^2 \cdots T_1 X_1\\
&=X_1Y_mY_1^{-1}+q(v^{-1}-v)T_1^{-1} \cdots T_{m-2}^{-1} T_m \cdots T_{n-1}X_n T_{n-1}^{-1} \cdots T_1^{-1}.
\end{split}
\end{equation*}

We will now prove, by induction on $l$ the following formula~:
\begin{equation}\label{E:L515}
u_l=X_1^lY_1^{-1}Y_m+(v^{-1}-v)T_1^{-1} \cdots T_{m-2}^{-1}T_m \cdots T_{n-1} e_l T_{n-1}^{-1} \cdots T_1^{-1}
\end{equation}
where
$$e_l=qX_nX_1^{l-1} + q^2X_n^2X_1^{l-2} + \cdots + q^lX_n^l.$$

The case $l=1$ is proved above. Let us assume that formula (\ref{E:L515}) holds for the integer $l$. We have
\begin{equation}\label{E:L516}
\begin{split}
u_{l+1}&=X_1u_l+q(v^{-1}-v)T_1^{-1} \cdots T_{m-2}^{-1}T_m \cdots T_{n-1}X_n T_{n-1}^{-1} \cdots T_1^{-1} Y_1 X_1^l Y_1^{-1}\\
&=X_1u_l+q^{l+1}(v^{-1}-v)T_1^{-1} \cdots T_{m-2}^{-1}T_m \cdots T_{n-1}X_n T_{n-1}^{-1} \cdots T_1^{-1}\times \\
& \hspace{2.6in} \times T_1 \cdots T_{n-1}X_n^l T_{n-1}^{-1} \cdots T_1^{-1}\\
&=X_1u_l+q^{l+1}(v^{-1}-v)T_1^{-1} \cdots T_{m-2}^{-1}T_m \cdots T_{n-1}X_n^{l+1}T_{n-1}^{-1} \cdots T_1^{-1}
\end{split}
\end{equation}
from which (\ref{E:L515}) follows for $l+1$ by the induction hypothesis.

Equation~(\ref{E:L512}) is obtained simply by multiplying (\ref{E:L515}) by $S$ on both sides. Lemma~\ref{L:51} is now proved.\qed

\vspace{.2in}

We can now start the proof of Proposition~\ref{P:2}. Let us form the generating series for
$S \bigg[ \sum_i Y_i, X_1^lY_1^{-1}\bigg]S$. By Lemma~\ref{L:51}, we find
\begin{equation*}
\begin{split}
\sum_{r\geq 1}S \bigg[ \sum_i Y_i, X_1^rY_1^{-1}\bigg]Su^r&=S\bigg\{-\frac{X_1u}{1-X_1u}+(v^{
2(1-n)}-1)\frac{\sum_{i \geq 1}q^iX_n^iu^i \cdot X_1u}{1-X_1u}\\
& \hspace{2in} +v^{2(1-n)}\frac{qX_nu}{1-X_nu}\bigg\}\\
&=S \frac{-X_1u+v^{2(1-n)}qX_nu}{(1-X_1u)(1-qX_nu)}S.
\end{split}
\end{equation*}

On the other hand, we have
\begin{equation*}
\begin{split}
exp\bigg(\sum_{r \geq 1}& \frac{(v^{-r}-v^r)(v^r-q^rv^{-r})}{r}\sum_{i}X_i^ru^r\bigg)\\
&=\frac{exp\bigg(\sum_{r \geq 1} \frac{1}{r}\sum_{i}X_i^ru^r\bigg)exp\bigg(\sum_{r \geq 1} \frac{q^r}{r}\sum_{i}X_i^ru^r\bigg)}{exp\bigg(\sum_{r \geq 1} \frac{v^{2r}}{r}\sum_{i}X_i^ru^r\bigg)exp\bigg(\sum_{r \geq 1} \frac{q^rv^{-2r}}{r}\sum_{i}X_i^ru^r\bigg)}\\
&=\prod_{i=1}^n\frac{exp\bigg(\sum_{r \geq 1} \frac{1}{r}X_i^ru^r\bigg)exp\bigg(\sum_{r \geq 1} \frac{q^r}{r}X_i^ru^r\bigg)}{exp\bigg(\sum_{r \geq 1} \frac{v^{2r}}{r}X_i^ru^r\bigg)exp\bigg(\sum_{r \geq 1} \frac{q^rv^{-2r}}{r}X_i^ru^r\bigg)}\\
&=\prod_{i=1}^n \frac{(1-v^2X_iu)(1-qv^{-2}X_iu)}{(1-qX_iu)(1-X_iu)}.
\end{split}
\end{equation*}

Hence we are reduced to proving the following relation~:

\begin{equation}\label{E:**}
\begin{split}
S\prod_{i=1}^n&(1-v^2X_iu)(1-v^{-2}qX_iu)S\\
=&S\prod_{i=1}^n(1-X_iu)(1-qX_iu)S\\
&+ \frac{(1-v^{2n})(1-v^{-2}q)}{1-q}S\bigg\{(X_1u-v^{2(1-n)}qX_nu)\prod_{2}^n (1-X_iu)\prod_{i=1}^{n-1}(1-qX_iu)\bigg\}S
\end{split}
\end{equation}
Of course we may, by homogeneity, drop the dummy variable $u$ in this formula. A brute force approach based on the equalities~:
$$SX_1S=\frac{1}{1+v^2}S(X_1+X_2)S={v^{-2}}SX_2S,$$
$$SX_1^2S=\frac{1}{1+v^2}S(X_1^2+X_2^2+(1-v^2)X_1X_2)S,$$
$$SX_2^2S=\frac{v^2}{1+v^2}S(X_1^2+X_2^2+(1-v^{-2})X_1X_2)S$$
allows one to check (\ref{E:**}) directly for $n=2$. We will now prove (\ref{E:**}) by induction on $n$.
So let us fix $n$ and assume that (\ref{E:**}) holds for the integer $n-1$, with $n-1 \geq 2$. For any subset $\{i_1, \ldots, i_r\}$ of $\{1, \ldots, n\}$ we denote by $S_{i_1, \ldots, i_r}$ the partial symmetrizer with respect to the indices $\{i_1, \ldots, i_r\}$.

Using the relation
$$S_{12}(X_1-v^{2(1-n)}qX_n)(1-X_2)S_{12}=v^{-2}S_{12}(X_2-v^{2(2-n)}qX_n)(1-v^2X_1)S_{12}$$
we get
\begin{equation*}
\begin{split}
&\frac{(1-v^{2n})(1-v^{-2}q)}{1-q}S\bigg\{(X_1-v^{2(1-n)}qX_n)\prod_{2}^n (1-X_i)\prod_{i=1}^{n-1}(1-qX_i)\bigg\}S\\
 &\quad= \frac{(1-v^{2n})(1-v^{-2}q)}{1-q}v^{-2}S\bigg\{(1-qX_1)(1-v^2X_1)(X_2-v^{2(2-n)}qX_n)\prod_{3}^n (1-X_i)\prod_{i=2}^{n-1}(1-qX_i)\bigg\}S\\
 &\quad=\frac{1-v^{2n}}{1-v^{2(n-1)}}v^{-2}S \bigg\{(1-qX_1)(1-v^2X_1)\bigg(\prod_{i=2}^n(1-v^2X_i)(1-v^{-2}qX_i)-\prod_{i=2}^n (1-X_i)(1-qX_i)\bigg)\bigg\}S
\end{split}
\end{equation*}

Next, we use the formulas
$$S(1-qX_1)\prod_{i=2}^n(1-v^{-2}qX_i)S=v^2\prod_{i=1}^n(1-v^{-2}qX_i)+(1-v^2)S\prod_{i=2}^n(1-v^{-2}qX_i)S,$$
$$S(1-v^2X_1)\prod_{i=2}^n(1-X_i)S=v^2\prod_{i=1}^n(1-X_i)+(1-v^2)S\prod_{i=2}^n(1-X_i)S$$
to simplify (\ref{E:**}) to the following relation~:
\begin{equation}\label{E:***}
\begin{split}
(v^{2n}&-v^{2(n-1)})S\bigg\{\prod_{i=1}^n (1-v^2X_i)(1-v^{-2}qX_i)-\prod_{i=1}^n(1-X_i)(1-qX_i)\bigg\}S\\
&=(1-v^{2n})(v^{-2}-1)\bigg( \prod_{i=1}^n(1-v^2X_i)S\prod_{i=2}^n (1-v^{-2}qX_i)S-\prod_{i=1}^n(1-qX_i)S\prod_{i=2}^n(1-X_i)S\bigg)
\end{split}
\end{equation}

Let as usual
$$m_{\lambda}(z_1, \ldots, z_t)=\sum_{\sigma \in \mathfrak{S}_t}z_{\sigma(1)}^{\lambda_1} \cdots z_{\sigma(t)}^{\lambda_t}$$
stand for the monomial symmetric function. The computation of the $(v)$-symmetrization of a monomial symmetric function $m_{(1^r)}$ is an easy exercise which we leave to the reader~:

\vspace{.1in}

\begin{slem}\label{C:*&}
For any $1 \leq r \leq n$, we have
$$S m_{(1^r)}(X_2, \ldots X_n)S=v^{2r}\frac{\begin{bmatrix} n-1\\ r \end{bmatrix}^{+}}{\begin{bmatrix} n\\ r \end{bmatrix}^+} Sm_{(1^r)}(X_1, \ldots, X_n)S.$$
\end{slem}

\vspace{.1in}

In the above we take the convention that $\begin{bmatrix} n-1\\ r \end{bmatrix}^{+}=0$ if $r=n$. Using Sublemma~\ref{C:*&}, we may now write down closed and symmetric expressions for all terms involved in (\ref{E:***})~:
$$\prod_{i=1}^n(1-v^2X_i)=\sum_{r=0}^{n}(-1)^rv^{2r}m_{(1^r)}(X_2, \ldots, X_n),$$
$$\prod_{i=1}^n(1-qX_i)=\sum_{r=0}^{n}(-1)^rq^{r}m_{(1^r)}(X_2, \ldots, X_n),$$
$$S\prod_{i=2}^n(1-X_i)S=\sum_{r=0}^{n}(-1)^rv^{2r}\frac{\begin{bmatrix} n-1\\ r \end{bmatrix}^{+}}{\begin{bmatrix} n\\ r \end{bmatrix}^+} Sm_{(1^r)}(X_1, \ldots, X_n)S,$$
$$S\prod_{i=2}^n(1-v^{-2}qX_i)S=\sum_{r=0}^{n}(-1)^rq^{r}\frac{\begin{bmatrix} n-1\\ r \end{bmatrix}^{+}}{\begin{bmatrix} n\\ r \end{bmatrix}^+} Sm_{(1^r)}(X_1, \ldots, X_n)S.$$

This allows us to write

\begin{equation*}
\begin{split}
\prod_{i=1}^n(1-v^2X_i)&S\prod_{i=2}^n (1-v^{-2}qX_i)S\\
=&\sum_{r,t}(-1)^r \bigg(\sum_{u=0}^r \begin{pmatrix} r\\ u \end{pmatrix} q^{u+t}\frac{\begin{bmatrix} n-1\\ u+t \end{bmatrix}^{+}}{\begin{bmatrix} n\\ u+t \end{bmatrix}^+} v^{2(r-u+t)}\bigg)Sm_{(1^r2^t)}(X_1, \ldots, X_n)S
\end{split}
\end{equation*}

and

\begin{equation*}
\begin{split}
\prod_{i=1}^n(1-qX_i)&S\prod_{i=2}^n (1-X_i)S\\
=&\sum_{r,t}(-1)^r \bigg(\sum_{u=0}^r \begin{pmatrix} r\\ u \end{pmatrix} v^{2(u+t)}\frac{\begin{bmatrix} n-1\\ u+t \end{bmatrix}^{+}}{\begin{bmatrix} n\\ u+t \end{bmatrix}^+} q^{r-u+t}\bigg)Sm_{(1^r2^t)}(X_1, \ldots, X_n)S
\end{split}
\end{equation*}

Hence,

\begin{equation}\label{E:wang2}
\begin{split}
\prod_{i=1}^n&(1-v^2X_i)S\prod_{i=2}^n (1-v^{-2}qX_i)S-\prod_{i=1}^n(1-qX_i)S\prod_{i=2}^n (1-X_i)S\\
=&\sum_{r,t}(-1)^r \sum_{u=0}^r \begin{pmatrix} r\\ u \end{pmatrix}\frac{\begin{bmatrix} n-1\\ u+t \end{bmatrix}^{+}}{\begin{bmatrix} n\\ u+t \end{bmatrix}^+}\bigg\{ q^{u+t}v^{2(r-u+t)}-v^{2(u+t)}q^{r-u+t}\bigg\}Sm_{(1^r2^t)}(X_1, \ldots, X_n)S
\end{split}
\end{equation}

while of course

\begin{equation}\label{E:wang1}
\begin{split}
S\bigg\{\prod_{i=1}^n (1-v^2X_i)(1-v^{-2}&qX_i)-\prod_{i=1}^n(1-X_i)(1-qX_i)\bigg\}S\\
&=\sum_{r,t}(-1)^r\sum_{u=0}^r \begin{pmatrix} r\\ u \end{pmatrix} (v^{2(r-2u)}-1)q^{u+t}m_{(1^r2^t)}(X_1, \ldots, X_n).
\end{split}
\end{equation}

Let $A$ denote the right-hand side of (\ref{E:wang2}) multiplied by $(1-v^{2n})(v^{-2}-1)$ and let $B$ stand for the right-hand side of (\ref{E:wang1}) multiplied by $(v^{2n}-v^{2(n-1)})$. Equation (\ref{E:***}) is simply that $A=B$. To show this, we check that the term $q^{u+t}m_{(1^r2^t)}(X_1, \ldots, X_n)$ appears in $A$ and $B$ with the same coefficient. For $B$ it is clearly equal to
$$(-1)^r \begin{pmatrix} r\\ u \end{pmatrix} (v^{2(r-2u)}-1)(v^{2n}-v^{2(n-1)}).$$
As far as $A$ is concerned, it is equal to

\begin{equation*}
\begin{split}
(1-v^{2n})&(v^{-2}-1)(-1)^r \begin{pmatrix} r\\ u \end{pmatrix} \bigg( v^{2(r-u+t)}\frac{\begin{bmatrix} n-1\\ u+t \end{bmatrix}^{+}}{\begin{bmatrix} n\\ u+t \end{bmatrix}^+}-v^{2(r-u+t)}\frac{\begin{bmatrix} n-1\\ r-u+t \end{bmatrix}^{+}}{\begin{bmatrix} n\\ r-u+t \end{bmatrix}^+}\bigg)\\
&\\
&=(1-v^{2n})(v^{-2}-1)(-1)^r \begin{pmatrix} r\\ u \end{pmatrix} v^{2(r-u+t)}\bigg(\frac{v^{2(n-r+u-t)}-v^{2(n-u-t)}}{1-v^{2n}}\bigg)\\
&=(v^{-2}-1)v^{2n}(1-v^{2(r-2u)})(-1)^r \begin{pmatrix} r\\ u \end{pmatrix}\\
&=(v^{2n}-v^{2(n-1)})(v^{2(r-2u)}-1)(-1)^r \begin{pmatrix} r\\ u \end{pmatrix}
\end{split}
\end{equation*}
as wanted. Equation (\ref{E:***}) and Proposition~\ref{P:2} are (finally) proved ! \qed

\vspace{.2in}

\section{Proof of Theorem~\ref{T:EisHeck}}

\vspace{.2in}

\paragraph{\textbf{C.1.}} We begin with equation (\ref{E:EisHeck1}). Since
$\eE_r(z_2)=\eE_r^{vec}(z_2)\eE_0(\nu^{2r}z_2)$ and since
$[\eE_0(z_1),\eE_0(\nu^{2r}z_2)]=0$, the relation (\ref{E:EisHeck1}) is equivalent to
\begin{equation}\label{E:Appc1}
[T_{(0,1)},\eE^{vec}_r(z)]=\nu \# \E(\mathbb{F}_l)\frac{\nu^{-2r}-1}{\nu^{-2}-1}z^{-1}\eE^{vec}_r(z).
\end{equation}
We prove (\ref{E:Appc1}) by showing that for any $x \in \E(\mathbb{F}_l)$,
$$[\mathbf{1}_{\mathcal{O}_x},\mathbf{1}^{vec}_{r,d}]=\frac{\nu^{-r}-\nu^r}{\nu^{-2}-1}\mathbf{1}^{vec}_{r,d+1},$$
where $\mathcal{O}_x$ is the structure sheaf at $x$. Indeed, we have
$$\mathbf{1}^{vec}_{r,d} \cdot\mathbf{1}_{\mathcal{O}_x}=\nu^{-r} \sum_{\substack{\mathcal{F}\;vec.\;bdle\\ \overline{\mathcal{F}}=(r,d)}} \mathbf{1}_{\mathcal{F} \oplus \mathcal{O}_x},$$
whereas, since every nonzero map to $\mathcal{O}_x$ is onto,
$$\mathbf{1}_{\mathcal{O}_x}\cdot\mathbf{1}^{vec}_{r,d}=
\nu^{r} \left\{\sum_{\substack{\mathcal{G}\;vec.\;bdle\\
\overline{\mathcal{G}}=(r,d+1)}} 
\frac{\#\text{Hom}(\mathcal{G},\mathcal{O}_x)-1}
{\nu^{-2}-1}\mathbf{1}_{\mathcal{G}}+
\sum_{\substack{\mathcal{F}\;vec.\;bdle\\ 
\overline{\mathcal{F}}=(r,d)}} 
\#\text{Hom}(\mathcal{F},\mathcal{O}_x)\mathbf{1}_{\mathcal{F} 
\oplus \mathcal{O}_x}\right\}.$$
We conclude using $dim\;\text{Hom}(\mathcal{G},\mathcal{O}_x)=dim\;\text{Hom}(\mathcal{F},\mathcal{O}_x)={r}$.

\vspace{.2in}

\paragraph{\textbf{C.2.}} We now turn to the proof of (\ref{E:EisHeck2}). We begin with a Lemma~:

\begin{lem} For any torsion sheaf $\mathcal{T}$, the series $\eE_r(z)$ is an eigenvector for the adjoint action of $\mathbf{1}_{\mathcal{T}}$.
\end{lem}
\noindent
\textit{Proof.} It is essentially the same as for the above case of $\mathcal{T}=\mathcal{O}_x$ (see {C.1.}). It suffices to notice that the number $\text{Surj}(\mathcal{G},\mathcal{T})$ of surjective maps
from a vector bundle $\mathcal{G}$ of rank $r$ to $\mathcal{T}$ is 
independent of the choice (and of the degree) of $\mathcal{G}$. 
This last statement is clear when $\mathcal{T}$ is stable, 
and may be proved in general by induction using the formula
$$\#\text{Hom}(\mathcal{G},\mathcal{T})=\sum_{\mathcal{T}'\subseteq \mathcal{T}} \# \text{Surj}(\mathcal{G},\mathcal{T}').$$
\qed

\vspace{.1in}

From the above Lemma and from the formula $\eE_0(z)=exp\left( \sum_r \frac{T_{(0,r)}}{[r]}z^r\right)$ we deduce that there exists a series $\mathcal{E}_r(z_1/z_2) \in \C[[z_1/z_2]]$ such that
\begin{equation}\label{E:appC2}
\eE_0(z_1)\eE_r(z_2)=\mathcal{E}_r\left(\frac{z_1}{z_2}\right)\eE_r(z_2)\eE_0(z_1).
\end{equation}
Let us first determine $\mathcal{E}_1(z_1/z_2)$. Relation (\ref{E:appC2}) for $r=1$ is equivalent to
\begin{equation}\label{E:appC3}
\eE_0(z_1)\eE^{vec}_1(z_2)=\mathcal{E}_1\left(\frac{z_1}{z_2}\right)\eE^{vec}_1(z_2)\eE_0(z_1).
\end{equation}
Thus, in order to compute $\mathcal{E}_1(z_1/z_2)$ it is enough to consider the restriction of
$\eE_0(z_1)\eE_1(z_2)$ to line bundles of degree, say, zero. If $\mathcal{L}$ is such a line bundle then for any $d>0$ we have
\begin{equation*}
\begin{split}
\mathbf{1}_{(0,d)} \cdot \mathbf{1}_{(1,-d)}(\mathcal{L})&=\nu^d \sum_{\mathcal{L}_{-d} \in Pic^{-d}(\E)} \frac{\#\text{Hom}(\mathcal{L}_{-d},\mathcal{L})-1}{\nu^{-2}-1}\\
&=\nu^d\#\E(\mathbb{F}_l) \frac{\nu^{-2d}-1}{\nu^{-2}-1}
\end{split}
\end{equation*}
wherefrom we get
\begin{equation*}
\begin{split}
\eE_0(z_1)\eE_1(z_2)(\mathcal{L})&=1+\sum_{d>0}\nu^{-d} \left(\frac{z_1}{z_2}\right)^d \mathbf{1}_{(0,d)}\mathbf{1}_{(1,-d)}(\mathcal{L})\\
&=1+\frac{\#\E(\mathbb{F}_l)}{\nu^{-2}-1}\sum_{d>0}\left(\frac{z_1}{z_2}\right)^d (\nu^{-2d}-1)\\
&=1+\frac{z_1}{z_2} \frac{\#\E(\mathbb{F}_l)}{(1-\frac{z_1}{z_2})(1-\nu^{-2}\frac{z_1}{z_2})}\\
&=\zeta\left( \frac{z_1}{z_2}\right).
\end{split}
\end{equation*}

This shows that $\mathcal{E}_1(z_1/z_2)=\zeta(z_1/z_2)$. Finally, to determine $\mathcal{E}_r(z_1/z_2)$, observe that by the coproduct formulas in Proposition~\ref{P:Eiscop},
\begin{equation*}
\begin{split}
\mathcal{E}_r\left(\frac{z_1}{z_2}\right)&\Delta_{1, \ldots, 1}(\eE_r(z_2)\eE_0(z_1))=\Delta_{1, \ldots, 1}(\eE_0(z_1)\eE_r(z_2))\\
&=\eE_0(z_1)\eE_1(z_2) \otimes \eE_0(z_1)\eE_1(\nu^{2}z_2) \otimes \cdots \otimes \eE_0(z_1)\eE_1(\nu^{2(r-1)}z_2)\\
&=\prod_{i=0}^{r-1}\zeta\left(\nu^{-2i}\frac{z_1}{z_2}\right)\eE_1(z_2) \eE_0(z_1) \otimes \cdots \otimes
\eE_1(\nu^{2(r-1)}z_2)\eE_0(z_1)\\
&=\prod_{i=0}^{r-1}\zeta\left(\nu^{-2i}\frac{z_1}{z_2}\right)\Delta_{1, \ldots, 1}(\eE_r(z_2)\eE_0(z_1)).
\end{split}
\end{equation*}
It follows that $\mathcal{E}_r(z_1/z_2)=\prod_{i=0}^{r-1}\zeta(\nu^{-2i}z_1/z_2)$ as desired. Theorem~\ref{T:EisHeck} is proved. \qed

\vspace{.2in}

\section{Proof of Lemma~\ref{L:strangeid}.}

\vspace{.2in}

Let us denote by $F(\a_1, \ldots, \a_r)$ the r.h.s. of (\ref{E:strangeid}). Each $\zeta\left( \frac{\a_i}{\a_j}\right)$ is a rational function in degree zero and leading coefficient one in both $\a_i$ and $\a_j$. From this and the expression (\ref{E:strangeid}) we see that $F(\a_1, \ldots, \a_r)$ is a rational function of degree one in each of the variables $\a_1, \ldots, \a_r$, and whose leading coefficient in any of these variables is also equal to one. Next, since $\zeta(z)$ has a simple pole at $z=1$ and $z=v^2$, the function $F(\a_1, \ldots, \a_r)$ has at most simple poles and these are located along the hyperplanes
$\a_i=\a_j$ and $\a_i=v^2\a_j$. We claim that the residues on each of these hyperplanes in fact vanish, so that $F(\a_1, \ldots, \a_r)$ is a polynomial in $\a_1, \ldots, \a_r$.

Indeed, the residues along hyperplanes $\a_i=\a_j$ vanish because $F(\a_1, \ldots, \a_r)$ is symmetric in $\a_1, \ldots, \a_r$; as for the hyperplanes $a_i=v^2\a_j$, we compute
\begin{equation*}
\begin{split}
\text{Res}_{v^2\a_j-\a_i}&\;F(\a_1, \ldots, \a_r)\\
&=\frac{1}{v^{-2}-1}\bigg[\prod_{\substack{l \neq i\\l \neq j} }\zeta \left( \frac{\a_l}{\a_j}\right) \cdot \text{Res}_{v^2\a_j-\a_i}\zeta\left( \frac{\a_i}{\a_j}\right) \cdot \bigg(\sum_{\substack{l \neq i\\l \neq j}} \a_l+v^2\a_j\bigg)\\
&\qquad \qquad\;\;\;-\prod_{\substack{l \neq i\\l \neq j} }\zeta \left( \frac{v^2\a_j}{\a_l}\right) \cdot \text{Res}_{v^2\a_j-\a_i}\zeta\left( \frac{\a_i}{\a_j}\right) \cdot \bigg(\sum_{\substack{l \neq i\\l \neq j}} \a_l+\a_j\bigg)\bigg]\\
&\quad+\prod_{\substack{l \neq i\\l \neq j} }\zeta \left( \frac{v^2\a_j}{\a_l}\right) \cdot \text{Res}_{v^2\a_j-\a_i}\zeta\left( \frac{\a_i}{\a_j}\right) v^2\a_j.
\end{split}
\end{equation*}
Using the relation $\zeta\left(\frac{v^2\a_j}{\a_l}\right)=\zeta\left(\frac{\a_l}{\a_j}\right)$ we simplify this to
\begin{equation*}
\begin{split}
\text{Res}_{v^2\a_j-\a_i}\;F(\a_1, \ldots, \a_r)&=\prod_{\substack{l \neq i\\l \neq j} }\zeta \left( \frac{\a_l}{\a_j}\right) \cdot\text{Res}_{v^2\a_j-\a_i}\zeta\left( \frac{\a_i}{\a_j}\right) \cdot \bigg\{ \frac{v^2\a_j-\a_j}{v^{-2}-1}-v^2\a_j\bigg\}\\
&=0
\end{split}
\end{equation*}
as wanted. Combining all the information we have on the function 
$F(\a_1, \ldots, \a_r)$ we see that necessarily 
$F(\a_1, \ldots, \a_r)=\a_1+\ldots + \a+r +u$ for some $u \in \mathcal{K'}$. 
It remains to observe that (for instance) we have $F(1,\ldots, 1)=r$. 
We are done.\qed

\vspace{.2in}

\centerline{\textbf{Acknowledgments.}}

\vspace{.1in}

We would like to thank Iain Gordon, Mark Haiman and Fran\c cois Bergeron for interesting discussions.
O.S. would  like to thank Pavel Etingof for a crucial suggestion made a few years ago concerning spherical DAHA's.

\newpage

\small{}

\vspace{4mm}

\noindent
O. Schiffmann, \texttt{schiffma@dma.ens.fr},\\
DMA, \'Ecole Normale Sup\'erieure, 45 rue d'Ulm, 75230 Paris Cedex 05, FRANCE.

\vspace{.1in}

\noindent
E. Vasserot, \texttt{vasserot@math.jussieu.fr},\\
D\'epartement de Math\'ematiques, Universit\'e de Paris 7, 175 rue du Chevaleret, 75013 Paris, FRANCE.

\end{document}